\documentclass[final,hidelinks,onefignum,onetabnum]{siamart220329}

\usepackage{lipsum}
\usepackage{amsfonts}
\usepackage{amssymb}
\usepackage{graphicx}
\usepackage{epstopdf}
\usepackage{subfigure}
\ifpdf
\DeclareGraphicsExtensions{.eps,.pdf,.png,.jpg}
\else
\DeclareGraphicsExtensions{.eps}
\fi

\usepackage{booktabs}
\usepackage{stmaryrd}
\SetSymbolFont{stmry}{bold}{U}{stmry}{m}{n}
\usepackage{bm}
\usepackage{caption}
\usepackage{subcaption}
\usepackage{multirow}
\usepackage{enumitem}
\usepackage{geometry}

\usepackage{cleveref}

\newsiamthm{assumption}{Assumption}

\renewcommand{\bf}{\textbf}
\newcommand{\fn}{\mathfrak{n}}
\newcommand{\fp}{\mathfrak{p}}

\newtheorem{thm}{Theorem}[section]

\newtheorem{lem}{Lemma}[section]
\newtheorem{coro}{Corollary}[section]
\newtheorem{rem}{Remark}[section]
\newtheorem{prop}[thm]{Proposition}
\newtheorem{assum}{Assumption}[section]

\newcommand{\bsub}{\begin{subequations}}
\newcommand{\esub}{\end{subequations}$\!$}

\newtheorem{example}{\textbf{Example}}[section]
\numberwithin{equation}{section}

\newsiamremark{remark}{Remark}
\newsiamremark{hypothesis}{Hypothesis}
\crefname{hypothesis}{Hypothesis}{Hypotheses}
\newsiamthm{claim}{Claim}
\newsiamremark{exmp}{Example}
\usepackage{algorithmic} 

\makeatletter
\newenvironment{breakablealgorithm}
{
		\begin{center}
			\refstepcounter{algorithm}
			\kern2pt \hrule height.8pt depth0pt \kern2pt
			\renewcommand{\caption}[2][\relax]{
				{\raggedright\textbf{\ALG@name~\thealgorithm} ##2\par}%
				\ifx\relax##1\relax 
				\addcontentsline{loa}{algorithm}{\protect\numberline{\thealgorithm}##2}%
				\else 
				\addcontentsline{loa}{algorithm}{\protect\numberline{\thealgorithm}##1}%
				\fi
				\kern2pt\hrule\kern2pt
			}
		}{
		\kern1pt \hrule\kern2pt
	\end{center}
}
\makeatother

\headers{}{}

\title{A second-order dynamical low-rank mass-lumped finite element method for the Allen-Cahn equation
}

\author{Jun Yang\thanks{School of Mathematics and Computational Science, Xiangtan University, Xiangtan 411105, P.R.China \email{yangjun@smail.xtu.edu.cn})}
\and Nianyu Yi\thanks{Hunan Key Laboratory for Computation and Simulation in Science and Engineering; School of Mathematics and Computational Science, Xiangtan University, Xiangtan 411105, P.R.China (\email{yinianyu@xtu.edu.cn}).}
\and Peimeng Yin\thanks{Corresponding author. Department of Mathematical Sciences, The University of Texas at El Paso, El Paso, TX 79968, USA (\email{pyin@utep.edu}).}} 

\usepackage{amsopn}

\ifpdf
\hypersetup{
  pdftitle={},
  pdfauthor={}
}
\fi

\usepackage{comment}

\begin{document}

\maketitle
\begin{abstract}
In this paper, we propose a novel second-order dynamical low-rank mass-lumped finite element method for solving the Allen-Cahn (AC) equation, a semilinear parabolic partial differential equation. The matrix differential equation of the semi-discrete mass-lumped finite element scheme is decomposed into linear and nonlinear components using the second-order Strang splitting method. The linear component is solved analytically within a low-rank manifold, while the nonlinear component is discretized using a second-order augmented basis update \& Galerkin (BUG) integrator, in which the $S$-step matrix equation is solved by the explicit 2-stage strong stability-preserving Runge-Kutta method. The algorithm has lower computational complexity than the full-rank mass-lump finite element method.
The dynamical low-rank finite element solution is shown to conserve mass up to a truncation tolerance for the conservative Allen-Cahn equation. Meanwhile, the modified energy is dissipative up to a high-order error and is hence stable. Numerical experiments validate the theoretical results. Symmetry-preserving tests highlight the robustness of the proposed method for long-time simulations and demonstrate its superior performance compared to existing methods.
\end{abstract}



\section{Introduction}
We focus on developing high-order numerical methods for solving gradient flows using low-dimensional surrogates that effectively capture their essential features. Dynamical low-rank approximation (DLRA) \cite{koch2007dynamical} is an emerging method that provides significant computational savings for solving many high-dimensional dynamical systems (e.g., \cite{grasedyck2013literature}). Gradient flows are inherently high-dimensional but typically evolve toward an equilibrium exhibiting a low-rank structure, such as phase separations or the disappearance of patterns. To leverage this property, we propose a dynamical low-rank finite element method for solving the Allen-Cahn (AC) equation, a gradient flow associated with the $L^2$ inner product.

The AC equation \cite{allen1979microscopic} was originally developed to model the motion of anti-phase boundaries in crystalline solids, and the classical AC equation is given by
\begin{equation}\label{1.1}
\left\{ \begin{aligned}
w_t &= \epsilon^2 \Delta w +f(w), & & x \in \Omega, \quad t \in(0,T],\\
w(&x,0)= w_0,& & x \in \Omega,
\end{aligned}
\right.
\end{equation}
subject to periodic or homogeneous Neumann boundary condition, $\partial_{\mathbf{n}}w =0$, on $\partial \Omega$. 
For simplicity of the presentation, we take $\Omega  \subset \mathbb{R}^{2}$ to be a convex, bounded domain.
Here, $w$ represents the concentration of one of the two metallic components of an alloy, and the parameter $\epsilon$ denotes the interfacial width. The nonlinear function $f(w)$ is given by $f(w) = -F^{\prime}(w)$, where $F(w) = \frac{1}{4}\left(1-w^{2}\right)^{2}$ is the double-well potential function.

The mass of solution in $\eqref{1.1}$ is not conserved, i.e., $\frac{d}{dt} \int_\Omega w dx \neq 0$. 
A nonlocal Lagrange multiplier was introduced in \cite{rubinstein1992nonlocal} to construct the conservative AC equation
\begin{equation}\label{1.2}
w_t = \epsilon^2 \Delta w + \bar{f}(w),
\end{equation}
where the modified nonlinear term $\bar{f}(w)$ in $\eqref{1.2}$ is defined as $\bar{f}(w):=f(w)-\lambda(w)$, and $\lambda(w)$ is the Lagrange multiplier enforcing mass conservation. Two primary forms of $\lambda(w)$ have been proposed, as described by Rubinstein and Sternberg \cite{rubinstein1992nonlocal}:
\begin{equation}\label{1.3}
\text{RSLM: } \lambda(w) = \frac{1}{|\Omega|} \int_{\Omega} f(w)dx,
\end{equation}
and the nonlocal multiplier proposed by Brassel and Bretin \cite{brassel2011modified}:
\begin{equation}\label{1.4}
\text{BBLM: } \lambda(w) = \frac{\int_{\Omega} f(w) d\mathbf{x}}{\int_{\Omega} \sqrt{4F(w)} d\mathbf{x}} \sqrt{4F(w)}.
\end{equation}
Both nonlocal Lagrange multipliers have been widely used in conservative AC equations \cite{ward1996metastable,kim2014conservative,lee2022high,zhang2022up,yang2023high}. Since the RSLM multiplier depends only on time, it has limitations in preserving small features. Brassel and Bretin suggested that the space-time dependent Lagrange multiplier \eqref{1.4} offered better mass conservation properties. 

The AC equation $\eqref{1.1}$ can be viewed as the $L^2$-gradient flow of the following Ginzburg-Landau energy functional:
\begin{equation}\label{Oener}
E(w) = \int_{\Omega} \left( \frac{\epsilon^2}{2} | \nabla w|^2 + F(w)  \right )d\mathbf{x}.
\end{equation}
Meanwhile, the conservative AC equation \eqref{1.2} can be interpreted as the $L^2$-gradient flow of the reformulated energy functional:
\begin{equation}\label{energyBBLM}
E(w) = \int_{\Omega} \left( \frac{\epsilon^2}{2} |\nabla w|^2 + F(w) + \int_{0}^{w} \lambda(s)  ds \right) d\mathbf{x}.
\end{equation}
The solutions of both AC equation \eqref{1.1} and its conservative AC equation \eqref{1.2} satisfy the energy dissipation law,
\begin{equation}\label{1.5}
\frac{d}{dt} E(w) = \left(\frac{\delta E(w)}{\delta u} , w_t \right) = - \left\| \frac{\partial w }{\partial t}\right\|^2 \leq 0, \qquad \forall t>0,
\end{equation}
where $E(w)$ is given by \eqref{Oener} for the classical AC equation \eqref{1.1}, and by \eqref{energyBBLM} for the conservative AC equation \eqref{1.2}.
Moreover, the conservative form \eqref{1.2} conserves mass, i.e.,
\begin{equation}\label{1.6}
\frac{d}{dt} \int_{\Omega} w d\mathbf{x} = 0.
\end{equation}

Preserving the physical properties, \(\eqref{1.5}\), and \(\eqref{1.6}\) in the conservative case, at the discrete level for the AC equation \(\eqref{1.1}\) and its conservative form \(\eqref{1.2}\), is a challenging and compelling research topic.
Many efficient numerical methods have been developed to preserve these properties, including the invariant energy quadratization (IEQ) methods \cite{yang2016linear,yang2017efficient, chen2024unconditionally}, scalar auxiliary variable (SAV) methods \cite{gong2020arbitrarily,shen2018scalar,shen2019new}, splitting methods \cite{li2022stability,li2022stability2}, integrating factor Runge-Kutta (IFRK) methods \cite{isherwood2018strong,ju2021maximum,zhang2021numerical}, and exponential time differencing (ETD) methods \cite{cox2002exponential,du2019maximum,du2021maximum,wang2016efficient}.
However, preserving the physical properties, even within acceptable tolerances, while reducing the computational complexity of these methods is a worthwhile endeavor.

The DLRA methods \cite{koch2007dynamical}, which can be traced back to the Dirac–Frenkel–McLachlan variational principle from the 1930s \cite{dirac1930note, frenkel1934wave}, have emerged as powerful tools for efficiently modeling high-dimensional systems whose dynamics evolve on lower-dimensional manifolds. 
More recently, the methods have been widely used in various fields, including weakly compressible flows \cite{einkemmer2019low}, radiation transport equations \cite{einkemmer2018low,peng2020low,einkemmer2021efficient}, hyperbolic problems \cite{kusch2022dynamical}, advection-diffusion equations \cite{nakao2023reduced}, and neural network applications \cite{schotthofer2022low}. Further reviews can be found in \cite{einkemmer2024review} and the references therein.

The DLRA methods evolve a dynamical system on the Riemannian manifold of fixed-rank matrices by projecting the right-hand side of a matrix differential equation onto the tangent space of the manifold, resulting in a set of differential equations that govern the factors of an SVD-type decomposition.
A projector-splitting integrator for factorized differential equations in DLRA was introduced in \cite{lubich2014projector}. Later on, Ceruti and Lubich proposed a robust,  \textit{basis-update \& Galerkin} (BUG) integrator, also known as an unconventional integrator. This alternative integrator eliminates the backward time integration substep inherent in the projector-splitting approach, which can be unstable for dissipative problems.  Additionally, the BUG integrator facilitates greater parallelism in its substeps.
However, the BUG integrator relies on a fixed rank, which limits its applicability when the optimal rank is either unknown or evolves dynamically during the computation. To overcome this limitation and enhance the efficiency of solving sub-differential equations in DLRA, Ceruti et al. introduced a rank-adaptive integrator \cite{ceruti2022rank}, which mitigates unwanted projection errors in the S-step of the BUG integrator by employing augmented bases. More recently, a second-order augmented BUG integrator was introduced, offering improved accuracy and robustness \cite{ceruti2024robust}.

Both the finite element method (FEM) and the discontinuous Galerkin (DG) method are high-order numerical approaches offering significant advantages, including high-order accuracy on compact stencils, compatibility with \textit{hp}-adaptivity, and flexibility in handling complex geometries. The integration of dynamical low-rank approximation (DLRA) with the DG method was explored for a homogeneous kinetic equation in \cite{yin2023semi}. However, applying DLRA to the DG bilinear form for the Laplace term in the AC equations introduces additional complexity. The formulation on cell interfaces leads to an approximation involving a sum of four additional SVD-like matrices with differing bases, necessitating the development of specialized solvers for the resulting DLRA system.
To address this complexity, we adopt FEM in the DLRA framework, where the bilinear form is approximated using a single SVD-like matrix. This approach significantly reduces the computational complexity of solving the DLRA while maintaining the high-order accuracy of FEM. 

Integrating DLRA with numerical methods typically requires the basis functions in different directions to be orthogonal. Unlike the discontinuous Galerkin (DG) method \cite{yin2023semi}, this orthogonality condition is not inherently satisfied by the FEM. To address this, we employ the mass-lumped finite element method \cite{chen1985lumped}, which substitutes the $L^2$ inner products in the variational formulation with piecewise fixed-point Gauss-Lobatto quadratures. While this approach introduces a numerical integration error relative to the original variational formulation, it does not compromise the accuracy of the finite element approximation.
Moreover, the mass-lumped finite element formulation provides essential benefits for integration with DLRA:
(i) The basis functions in different directions are orthogonal, making it well-suited for the integration.
(ii) The coefficient matrices for polynomial-type nonlinear terms can be explicitly expressed using the Hadamard product, facilitating the low-rank approximation of the nonlinear terms.

To accurately approximate the solution to the AC problems, we primarily focus on second-order time discretization, with only a brief mention of first-order discretization. Specifically, we employ the second-order Strang splitting method to decouple the linear and nonlinear components of the AC equations. For the linear part, we solve it analytically within a low-rank manifold. For the nonlinear part, we design a second-order augmented BUG integrator, inspired by the approach in \cite{ceruti2024robust}. Specifically, we will employ the explicit 2-stage strong stability preserving Runge-Kutta method for the S-step matrix differential equation. Building on these techniques, we develop the dynamical low-rank mass-lumped finite element method (DLR-MLFEM). This method retains the structure of the full-rank mass-lumped finite element method (FR-MLFEM) but computes the solution in subspaces of the finite element space. The computational complexity of the proposed second-order augmented BUG integrator or the DLR-MLFEM is $\mathcal{O}((m+n)r^{4})$, where $m$ and $n$ are the total number of unknowns in each direction, and $r$ is the rank of the manifold. 
The dynamical low-rank finite element solution is shown to conserve the mass up to a truncation tolerance. The modified energy is dissipative up to a high-order error of $\mathcal{O}(\tau(\tau^2+h^{k+1}))$, where $\tau$ is the time step, $h$ is the mesh size, and $k$ is the degree of the polynomials. Consequently, the energy stability property remains valid.
Numerical examples are provided to validate the theoretical findings, with symmetry-preserving tests underscoring the robustness of the proposed method for long-time simulations.

The remainder of the paper is organized as follows. In \Cref{Sec2}, we introduce the formulation of the matrix differential equations and propose a second-order FR-MLFEM. In \Cref{Sec3}, we develop a second-order DLR-MLFEM based on the Strang-splitting method, incorporating an exact integrator for the linear component and a second-order augmented BUG integrator for the nonlinear component. Additionally, we analyze the conservation properties of the proposed method. Numerical examples are provided in \Cref{Sec4}, and concluding remarks are presented in \Cref{Sec5}.

\section{Matrix differential equation}\label{Sec2}
Consider the AC-type semilinear parabolic equation 
\begin{equation}\label{2.1}
\begin{aligned}
&w_t=\epsilon^2 \Delta w + \mathcal{N}(w), & &\text { in } \Omega \times(0, T), \\
&w({x}, t=0)=w_0({x}), & &\text { in } \Omega \times\{0\}, \\ 
&\partial_{\mathbf{n}} w=0,  & &\text { on } \partial \Omega \times(0, T),
\end{aligned}
\end{equation}
where the $\mathcal{N}(w)=f(w)$ for the classical AC equation \eqref{1.1}, and $\mathcal{N}(w)=\bar{f}(w)$ for the conservative AC equation \eqref{1.2}. 

\subsection{1D semi-discrete finite element scheme}
Let $\Omega = I = [a,b]$ and $a=x_1 < x_2 < \ldots < x_{Mk+1} = b$ be the Gauss-Lobatto points, where $x_{ik+1}$, $i=0,\ldots, M$ are the nodal points. Consider a uniform mesh size $h = x_{ik+1} - x_{(i-1)k+1} =(b-a)/M$. We define the finite element space 
\begin{equation}\label{2.2}
V_{I,h}^k = \left\lbrace v \in H^{1}(I): v|_{I_i} \in \mathbb{P}^k, i =1,...,M \right\rbrace,
\end{equation}
where $I_i = [x_{(i-1)k+1},x_{ik+1}]$ and $\mathbb{P}^k$ denotes the space of polynomials of degree no more than $k$. 
Let $x_{(i-1)k+j}$ and $\omega_j$, $j=1, \ldots, k+1$ be the quadrature points and weights of $(k+1)$-point Gauss-Lobatto quadrature on the subinterval $I_i$, and denote 
\begin{equation}\label{GLweights}
\bar{\omega}_{(i-1)k+j} = 
\left\lbrace
\begin{aligned}
\omega_j, & \quad 2\leq j \leq k, \\
2\omega_j, & \quad j=1,k+1.
\end{aligned}
\right.
\end{equation}
Given functions $f,g$, the $L^2$ inner product
\[(f, g) = \int_{\Omega} f g dx\]
is approximated by the piecewise Gauss-Lobatto quadrature
\begin{equation}\label{GaussLobatto}
(f,g)_h = \sum\limits_{j=1}^{Mk+1} \bar{\omega}_j f(x_j)g(x_j).
\end{equation}
The norm associated with the inner product in \eqref{GaussLobatto} is defined by
\[\|f\|_h = \sqrt{(f,f)_h}, \qquad \forall f \in V^k_h.\]

The variational formulation of the equation $\eqref{2.1}$ is to find $w \in L^2((0,T];H^1(\Omega))$ such that:
\begin{equation}\label{2.5}
(w_t,v) = -\epsilon^2 (\partial_x w, \partial_x v) + (\mathcal{N}(w),v),  \quad \forall v \in H^1(\Omega),
\end{equation}
which upon using the Gauss-Lobatto quadrature gives
\begin{equation}\label{2.5+}
(w_t,v_h)_h = -\epsilon^2 (\partial_x \Pi_h w,  \partial_x v_h) + (\mathcal{N}(w),v_h)_h + \mathcal{E}(v_h), \quad \forall v_h \in V_h^k,
\end{equation}
where $\Pi_h: C(\bar{\Omega}) \to V^k_h$ is the Lagrange interpolation operator and 
\[\mathcal{E}(v_h) =  (w_t,v_h)_h - (w_t,v_h) + (\mathcal{N}(w),v_h)-(\mathcal{N}(w),v_h)_h + (\partial_x (\Pi_h w - w), \partial_x v_h).\]
\begin{rem}
The error term, $\mathcal{E}(v_h)$, which includes both quadrature and interpolation errors, satisfies the following error bound \cite{li2020arbitrarily}:
\[|\mathcal{E}(v_h)| \leq \mathcal{O}(h^k),\]
and studies imply that this error does not impact the convergence rate of the finite element solution.
\end{rem}
Note that the $(k+1)$-point Gauss-Lobatto quadrature is exact for polynomials of degree up to $2k-1$. For $w_h \in V_h^k$, it holds
\begin{equation}\label{hinnereq}
    (\partial_x \Pi_h w_h,  \partial_x v_h)=(\partial_x w_h,  \partial_x v_h)=(\partial_x w_h,  \partial_x v_h)_h.
\end{equation}
The spatially semi-discrete mass-lumped finite element method for $\eqref{2.1}$, derived using \eqref{2.5+} and \eqref{hinnereq}, is to find $w_h \in L^2((0,T]; V^k_h)$ such that
\begin{subequations}\label{2.7-}
    \begin{align}
        & (w_{h,t},v_h)_h = -\epsilon^2 (\partial_x w_h, \partial_x v_h)_h + (\mathcal{N}(w_h),v_h)_h, \quad \forall v_h \in V_h^k, \\
        & (w_h(x,0), v_h) = (w_0, v_h), \quad \forall v_h \in V_h^k.
    \end{align}
\end{subequations}

\subsection{Matrix differential equation in 2D} 
We extend the implementation using the same setting as one dimension to describe the proposed spatially semi-discrete scheme on a 2D rectangular domain $\Omega =[a,b] \times[c,d]$. 
To simplify the notation, we denote
\begin{equation}\label{2.8}
m = Mk+1, \quad n = Nk+1.
\end{equation}
Let $a = x_1 < x_2 < \cdots < x_{m}=b$ represent the Gauss-Lobatto points on the interval $\Omega_x= [a,b]$, where $x_{ik+1}$, $i=0,\ldots, M$, are the nodal points.
The uniform mesh size in the $x$-direction is $h_x = x_{ik+1} - x_{(i-1)k+1} = (b-a)/M$ for $i=1,\ldots, M$. Similarly, let $c = y_1 < y_2 < \cdots < y_{n}=d$ be the Gauss-Lobatto points on the interval $\Omega_y = [c,d]$, where
$y_{jk+1}$, $j=0,\ldots,N$, are the nodal points. The uniform mesh size in the $y$-direction is $h_y = y_{jk+1} - y_{(j-1)k+1} = (d-c)/N$ for $j=1,\ldots, N$. 
The domain $\Omega = \Omega_x \times \Omega_y$ is now divided into $M \times N$ subrectangles, collectively denoted by by $\mathcal{K}$, using the grid points $(x_{ik+1}, y_{jk+1})$, with $0 \leq i \leq M$ and $0 \leq j \leq N$. The mesh size of the partition $\mathcal{K}$ is given by $h = \max \{h_x, h_y\}$.
We define the finite element method space in each direction as:
\begin{align*}
V_{x,h}^k = & \left\lbrace \phi \in H^{1}(\Omega_x): \phi|_{\Omega_{x,i}} \in \mathbb{P}^k, i =1,...,M \right\rbrace,\\
V_{y,h}^k = & \left\lbrace \psi \in H^{1}(\Omega_y): \psi|_{\Omega_{y,j}} \in \mathbb{P}^k, j =1,...,N \right\rbrace,
\end{align*}
where $\Omega_{x,i} = [x_{(i-1)k+1},x_{ik+1}]$ and $\Omega_{y,j} = [y_{(j-1)k+1},y_{jk+1}]$. 
The $H^1$-conforming tensor-product finite element space for $\Omega$ is given by
\[Q_h^k = V_{x,h}^{k} \otimes V_{y,h}^k  =  \{ v \in H^1(\Omega): v|_{K} \in \mathcal{Q}_k, \quad \forall K \in \mathcal{K} \}.\]
where $\mathcal{Q}_k$ be space of polynomials in the variables $x,y$ with real coefficients, and of degree at most $k$ in each direction on $K$.
Similar to the 1D case, the $L^2$  inner product in 2D is approximated using the tensor product of the 1D Gauss-Lobatto quadrature: 
\[(f,g)_h = \sum\limits_{i=1}^{m}\sum\limits_{j=1}^{n} w_i w_j f(x_i,y_j) g(x_i,y_j).\]
With this Gauss-Lobatto quadrature approximation, the 2D spatially semi-discrete mass-lumped finite element method for $\eqref{2.2}$ is to find $w_{h,t} \in Q_h^k$ such that
\begin{subequations}\label{2.9}
    \begin{align}
        & (w_{h,t},v_h)_h =  -\epsilon^2 (\nabla w_h, \nabla v_h) + (\mathcal{N}(w_h),v_h)_h, \quad \forall v_h \in Q_h^k,\label{2.9a}\\
        & (w_h|_{t=0}, v_h) =  (w_0, v_h), \quad \forall v_h \in Q_h^k.\label{2.9b}
    \end{align}
\end{subequations}

Let $\{ \phi_i(x)\}_{i=1}^{m}$ and $\{ \psi_j(y)\}_{i=1}^{n}$ be basis functions for the finite elements spaces in the $x$- and $y$- directions, respectively. For any  $1\leq p \leq m$, $1\leq q \leq n$, let $x_p$ and $y_q$ represent the corresponding Gauss-Lobatto quadrature points. Then it follows
\begin{equation}\label{basisdelta}
   \phi_i(x_l) = \delta_{il}, \quad \psi_j(y_q) = \delta_{jq},
\end{equation}
where $\delta_{ij}$ is Kronecker delta function.
Moreover, for any $s>0$, it holds
\begin{equation}\label{basis01}
(\phi_i(x_p))^s = \phi_i(x_p), \quad \text{ and } \quad (\psi_j(y_q))^{s} = \psi_j(y_q).
\end{equation}

Denoted by
\begin{equation}\label{2.10}
\Phi(x) = [\phi_1(x), \phi_2(x), \cdots, \phi_{m}(x)]^{\top} \text{ and } \Psi(y) = [\psi_1(y), \psi_2(y), \cdots, \psi_{m}(y)]^{\top}.
\end{equation}
Then we can define the lumped mass matrices
\begin{equation}\label{2.11}
M_x = (\Phi(x),\Phi(x)^{\top})_{\Omega_x} = [(\phi_i, \phi_j)_h ]_{m \times m}, \quad M_y = (\Psi(y),\Psi(y)^{\top})_{\Omega_y} = [(\psi_i, \psi_j)_h ]_{n \times n},
\end{equation}
which are sparse and diagonal.
We can similarly introduce the stiff matrices 
\begin{equation}\label{2.12}
\begin{aligned}
A_x =  (\partial_x \Phi(x), \partial_x \Phi(x)^{\top})_{\Omega_x} = [(\partial_x \phi_i ,\partial_x \phi_j )_h ]_{m \times m}, \quad 
A_y =  (\partial_y \Phi(y), \partial_y \Phi(y)^{\top})_{\Omega_y} = [(\partial_y \psi_i ,\partial_y \psi_j )_h ]_{n \times n},
\end{aligned}
\end{equation}
which are sparse, symmetric, and banded with a bandwidth $2k+1$.

Given a function $w_h \in Q^k_h$, its expansion in terms of these bases can be written as:
\begin{equation}\label{2.13}
w_h = \sum\limits_{i=1}^{m} \sum\limits_{j=1}^{n} W_{i,j}(t) \phi_i(x) \psi_j(y) = \Phi(x)^{\top} W(t) \Psi(y),
\end{equation}
The matrix $W = [W_{i,j}] \in \mathbb{R}^{m \times n}$ is called the coefficient matrix of $w_h$, with respect to the bases $\{ \phi_i(x)\}_{i=1}^{m}$ and $\{ \psi_j(y)\}_{i=1}^{n}$. 

\begin{definition}\label{frobdef}
Given matrices $A,B \in \mathbb{R}^{m \times n}$ with entries $A_{ij}$ and $B_{ij}$, their {Frobenius inner product} is defined as
$
(A,B)_{\rm{F}}
=\text{tr}(A^\top B) 
=\sum_{i=1}^{m} \sum_{j=1}^{n} A_{ij}B_{ij}.
  $
The {Frobenius norm} of $A$ is 
$
\|A\|_{\rm{F}} = \sqrt{( A,A)_{\rm{F}}}.
$
The $M$-weighted Frobenius inner product is $(A,B)_{\rm M} = (M_x A M_y,B)_\mathrm{F}$ with $M_x$ and $M_y$ given in \eqref{2.11}.
The $\rm M$-weighted Frobenius norm of $A$ is $\|A\|_{\rm M}:= \sqrt{(A,A)_{\rm M}}$.
\end{definition}

\begin{definition}\label{hadadef}
Given matrices \( A, B \in \mathbb{R}^{m \times n} \) with entries \( A_{ij} \) and \( B_{ij} \), the {Hadamard product} 
\[
A \odot B \in \mathbb{R}^{m \times n}, \quad \text{with entries} \quad (A \odot B)_{ij} = A_{ij} B_{ij}.
\]
For any integer \( s \geq 1 \), the {Hadamard power} \( A^{\circ s} \in \mathbb{R}^{m \times n} \) is defined as  
\[
A^{\circ s} = \underbrace{A \odot A \odot \cdots \odot A}_{s \text{ times}}, \quad \text{with entries} \quad (A^{\circ s})_{ij} = A_{ij}^s.
\]
Additionally, for \( C \in \mathbb{R}^{m \times n} \), the distributive property of the Hadamard product holds:  
\[
A \odot (B + C) = A \odot B + A \odot C.
\]
\end{definition}

Then the following result holds.

\begin{lem}\label{lem2.2}
For any $z_h, w_h \in Q_h^k$, let $Z, W \in \mathbb{R}^{m \times n}$ be their coefficient matrices, respectively. Then,
\begin{equation}\label{2.14}
(z_h,w_h)_h = (M_x Z M_y,W)_{\rm F} = (Z,W)_{\rm M},
\end{equation} 
and for any integer $s>0$,
\begin{equation}\label{2.15}
(z_h^s,w_h)_h = (M_x Z^{\circ s} M_y,W)_{\rm F} = (Z^{\circ s},W)_{\rm M}.
\end{equation} 
Moreover, 
\begin{equation}\label{2.16}
\begin{aligned}
(\nabla z_h, \nabla w_h)_h = & (A_x Z M_y + M_x Z A^{\top}_y,W)_{\rm F} = -(L_x Z + Z L^{\top}_y,W)_{\rm M},
\end{aligned}
\end{equation} 
where the matrices $L_x =  -M^{-1}_x A_x \in \mathbb{R}^{m\times m}, \ L_y = -M^{-1}_y A_y \in \mathbb{R}^{n \times n}$.
\end{lem}
\begin{proof}
Note that $z_h = \Phi(x)^{\top} Z \Psi(y)$, $w_h = \Phi(x)^{\top} W\Psi(y)$. Then,  
\begin{equation}
\begin{aligned}
(z_h,w_h)_h & = \left( \sum\limits_{i=1}^{m}\sum\limits_{j=1}^{n} Z_{ij} \phi_i(x) \psi_j(y), \sum\limits_{k=1}^{m}\sum\limits_{l=1}^{n} W_{kl}\phi_k(x) \psi_l(y) \right)_h  = \sum\limits_{ijkl} Z_{ij} W_{kl} (\phi_i,\phi_k)_h (\psi_j,\psi_l)_h \\
& = \sum\limits_{kl} (M_x Z M_y^{\top})_{kl} W_{kl} = (M_x Z M_y^{\top}, W) = (M_x Z M_y, W)_{\rm F}.
\end{aligned}
\end{equation}
where $M_x$ and $M_y$ are defined by $\eqref{2.11}$.
For any  $1\leq p \leq m$, $1\leq q \leq n$, let $x_p$ and $y_q$ represent the Gauss-Lobatto quadrature points in $x$- and $y$- directions, respectively,
\[\begin{aligned}
(z_h^s, w_h)_h &= \left( \left(\sum\limits_{i=1}^{m} \sum\limits_{j=1}^{n} Z_{ij} \phi_i(x) \psi_j(y)\right)^{s}, \sum\limits_{\nu=1}^{m} \sum\limits_{l=1}^{n} W_{\nu l} \phi_\nu(x) \psi_l(y)\right)_h \\
&=\sum\limits_{p=1}^{m} \sum\limits_{q=1}^{n} w_p w_q \sum\limits_{i=1}^{m} \sum\limits_{j=1}^{n} Z^s_{ij} (\phi_i(x_p))^{s} (\psi_j(y_q))^{s} \sum\limits_{\nu=1}^{m} \sum\limits_{l=1}^{n} W_{\nu l} \phi_\nu(x_p) \psi_l(y_q) \\
&= \sum\limits_{p=1}^{m} \sum\limits_{q=1}^{n} w_p w_q \sum\limits_{i=1}^{m} \sum\limits_{j=1}^{n} Z^s_{ij} \phi_i(x_p) \psi_j(y_q) \sum\limits_{\nu=1}^{m} \sum\limits_{l=1}^{n} W_{\nu l} \phi_\nu(x_p) \psi_l(y_q) \quad \text{(by \eqref{basis01})} \\
& = \left( \sum\limits_{i=1}^{m} \sum\limits_{j=1}^{n} Z^s_{ij} \phi_i(x) \psi_j(y) , \sum\limits_{\nu=1}^{m} \sum\limits_{l=1}^{n} W_{\nu l} \phi_\nu(x) \psi_l(y) \right)_h \\
& = (M_x Z^{\circ s} M_y , W )_{\rm F} =(Z^{\circ s} , W )_{\rm M}.
\end{aligned}\] 

Similarly, 
\[
\begin{aligned}
(\nabla z_h, \nabla w_h)_h&= \left( \nabla \left(\sum\limits_{i=1}^{m}\sum\limits_{j=1}^{n} Z_{ij} \phi_i(x) \psi_j(y) \right), \nabla \left( \sum\limits_{\nu=1}^{m}\sum\limits_{l=1}^{n} W_{\nu l} \phi_\nu(x) \psi_l(y) \right) \right)_h  \\
& = \sum\limits_{ij\nu l} Z_{ij}W_{\nu l} (\partial_x \phi_i, \partial_x \phi_k)_h  (\psi_\nu, \psi_l)_h +\sum\limits_{ij\nu l} Z_{ij}W_{\nu l} (\partial_y \psi_i, \partial_y \psi_\nu)_h  (\phi_j, \phi_l)_h\\
& = \left( A_x Z M_y, W \right)_{\rm F} + \left( M_x Z A^{\top}_y, W \right)_{\rm F}.
\end{aligned}
\]
where $A_x$ and $A_y$ are defined by $\eqref{2.12}$. 
\end{proof}

We introduce the coefficient matrix associated with the function \(1 \in Q_h^k\), which is the all-ones matrix \(\mathbf{I} = q_m q_n^\top \in \mathbb{R}^{m \times n}\), where each entry is equal to 1. Here, \(q_m \in \mathbb{R}^m\) and \(q_n \in \mathbb{R}^n\) are all-ones vectors.
\begin{coro}\label{intoin}
For any $w_h \in Q_h^k$, let $W \in \mathbb{R}^{m \times n}$ be their coefficient matrices,
\begin{equation}\label{nlmat}
(\mathcal{N}(w_h),v_h)_h = (\mathcal{N}(W), V)_{\rm M},
\end{equation}
where $V$ is the coefficient matrix of $v_h$, and
\begin{subequations}\label{matrix-nonlinear}
    \begin{align}
        &\mathcal{N}(W) := W-W^{\circ 3}, & & \text{(for \eqref{1.1})}, \label{matrix-nonlineara}\\
        &\mathcal{N}(W) := W-W^{\circ 3} - \frac{1}{|\Omega|} (W-W^{\circ 3}, \mathbf{I}), & & \text{(for \eqref{1.2} with RSLM),}\\
&\mathcal{N}(W) := W-W^{\circ 3} - \frac{(W-W^{\circ 3},\mathbf{I})_{\rm M}}{((\mathbf{I} - W^{\circ 2}),\mathbf{I})_{\rm M}} (\mathbf{I}-W^{\circ 2}), & & \text{(for \eqref{1.2} with BBLM)},\label{matrix-nonlinearb}
    \end{align}
\end{subequations}
where the matrix $\mathbf{I}\in \mathbb{R}^{m \times n}$ is an all-ones matrix.

Moreover, the conservative AC equation \eqref{1.2} with RSLM and BBLM holds
\begin{equation}\label{2.27}
(\mathcal{N}(w_h),1)_h = (\mathcal{N}(W), \mathbf{I})_{\rm M} = 0.
\end{equation}
\end{coro}
By \Cref{lem2.2} and \Cref{intoin}, we have the following result.
\begin{coro}\label{fullfemmat}
Let $W \in \mathbb{R}^{m \times n}$ be the coefficient matrix of the finite element solution $w_h \in Q_{h}^{k}$, and $V \in \mathbb{R}^{m \times n}$ be the coefficient matrix of any function $v_h \in Q_{h}^{k}$. Then, the semi-discrete \eqref{2.9} is equivalent to the following problem: Find $W(t) \in \mathbb{R}^{m \times n}$ such that
\begin{subequations}\label{matvar}
    \begin{align}
        &(\partial_t W(t), V)_{\rm M} = (\mathcal{L}(W(t)), V)_{\rm M} + (\mathcal{N}(W(t)), V)_{\rm M}.\\
        &W(t_0) = W_0, 
    \end{align}
\end{subequations}
where $W_0 \in \mathbb{R}^{m \times n}$ is the coefficient matrix of $w_h(x,y,0)$ in \eqref{2.9b}. Here, the linear term $\mathcal{L}(W) = \epsilon^2( L_x W + W L^{\top}_y)$ 
and the nonlinear term $\mathcal{N}(W(t))$ is given by $\eqref{matrix-nonlinear}$
\end{coro}

By \Cref{fullfemmat}, the 2D semi-discrete finite element scheme \eqref{2.9} can be rewritten as the matrix differential equation
\begin{equation}\label{2.19}
\begin{aligned}
&\partial_t W = \mathcal{L}(W(t)) + \mathcal{N}(W(t)), \\
&W(t_0) = W_0.
\end{aligned}
\end{equation}

\subsection{Second-order FR-MLFEM}
Before presenting the dynamical low-rank approximation for solving the matrix differential equation \eqref{2.19}, we introduce a second order full-rank mass-lumped finite element method (FR-MLFEM).

For $\mathfrak{n}\geq 0$, let $w_h^\mathfrak{n} = w_h(x,y,t_\mathfrak{n}) \in Q^k_h$ be an approximation of $w(x,y,t_\mathfrak{n})$, where $t_\mathfrak{n}=\mathfrak{n}\tau$ and $\tau>0$ is a specified time step. We further denote the coefficient matrix of $w_h^\mathfrak{n}$ by $W_\mathfrak{n} \in \mathbb{R}^{m\times n}$, namely, $w_h^\mathfrak{n} = \Phi(x)^{\top} W_\mathfrak{n} \Psi(y)$.

Given $w_h^\mathfrak{n}$ or its coefficient matrix $W_h^{\mathfrak{n}}$, to design a second-order numerical integrator that computes $W_{\mathfrak{n}+1}$ from the nonlinear matrix differential equation \eqref{2.19}, we consider splitting the matrix equation in \eqref{2.19} into two sub-equations associated with its linear and nonlinear terms,
\begin{subequations}\label{2.20}
\begin{align}
&\begin{aligned}
X_t  = \mathcal{L} (X),
\end{aligned} \label{2.20a}\\
&\begin{aligned}
Z_t = \mathcal{N}(Z).
\end{aligned} \label{2.20b}
\end{align}
\end{subequations}

We denote by $S_{\tau}^{\mathcal{L}} X_{\mathfrak{n}}$ and $S_{\tau}^{\mathcal{N}} Z_{\mathfrak{n}}$ the solution for $\eqref{2.20a}$ and $\eqref{2.20b}$ at $t_{\mathfrak{n}+1} = t_{\mathfrak{n}}+\tau$ with initial conditions $X_{\mathfrak{n}}$ and $Z_{\mathfrak{n}}$, respectively. 
We first solve the linear equation $\eqref{2.20a}$. 
For any integer $r>0$ and matrix $Z \in \mathbb{R}^{r \times r}$, we define the matrix exponential as 
$$e^Z := \sum_{l=0}^\infty \frac{Z^l}{l!} = I_r + \sum_{l=1}^\infty \frac{Z^l}{l!} \in \mathbb{R}^{r \times r}.$$
Here and throughout what follows, $I_r$ represents the identity matrix of dimension $r$.

For any matrix $W \in \mathbb{R}^{m \times n}$ and $s \in \mathbb{R}$, we define the operator $e^{s \mathcal{L}}$ by 
\[e^{s \mathcal{L}} W : = e^{s L_x}  W  e^{s L^{\top}_y}, \]
where $L_x$ and $L_y$ are defined in \eqref{2.14}.
It can be verified that the inverse operator of $e^{s \mathcal{L}}$ is given by $e^{-s\mathcal{L}}$, namely,
\[ e^{s \mathcal{L}} e^{-s \mathcal{L}} W  =e^{-s \mathcal{L}} e^{s \mathcal{L}} W  = e^{-s L_x}e^{s L_x}  W  e^{s L^{\top}_y}e^{-s L^{\top}_y} = W. \] 

The linear equation $\eqref{2.20a}$, which a linear ordinary differential equation in matrix form, could be solved analytically, and the solution is given by
\begin{equation}\label{2.21}
X_{\mathfrak{n}+1} = S_\tau^{\mathcal{L}} X_{\mathfrak{n}}= e^{\epsilon^2 \tau \mathcal{L}} X_{\mathfrak{n}}.
\end{equation}

Next, for the nonlinear problem $\eqref{2.20b}$, we use the explicit 2-stage strong stability preserving Runge-Kutta (SSP-RK2) to approximate the solution at $t_{\mathfrak{n}+1}$ 
with the initial condition $Z_\mathfrak{n}$ by
\begin{equation}\label{2.22}
\left\lbrace
\begin{aligned}
&Z_{\mathfrak{n},1} = Z_\mathfrak{n} + \tau \mathcal{N}(Z_\mathfrak{n}),\\
&Z_{\mathfrak{n}+1} = Z_\mathfrak{n} + \frac{\tau}{2} \mathcal{N}(Z_\mathfrak{n}) + \frac{\tau}{2} \mathcal{N}(Z_{\mathfrak{n},1}),
\end{aligned}
\right.
\end{equation}
which can be represented in an abstract form as:
\begin{equation}\label{2.23}
Z_{\mathfrak{n}+1} = S_\tau^{\mathcal{N}} Z_{\mathfrak{n}}.
\end{equation}

The FR-MLFEM for \eqref{2.19} formulated in matrix form proceeds as follows: Given $W_\mathfrak{n} \in \mathbb{R}^{m \times n}$, we find $W_{\mathfrak{n}+1} \in \mathbb{R}^{m \times n}$ using a second-order Strang splitting approach, relying on the solvers outlined in $\eqref{2.21}$ and $\eqref{2.23}$, 
\begin{equation}\label{2.24}
\left\lbrace
\begin{aligned}
&W_{\mathfrak{n},1} = e^{\frac{\tau}{2} \epsilon^2 \mathcal{L}} W_\mathfrak{n} ,\\
&W_{\mathfrak{n},2} = W_{\mathfrak{n},1} + \tau \mathcal{N}(W_{\mathfrak{n},1}),\\
&W_{\mathfrak{n},3} = \frac{1}{2} W_{\mathfrak{n},1} + \frac{1}{2} W_{\mathfrak{n},2} + \frac{1}{2} \tau \mathcal{N}(W_{\mathfrak{n},2}),\\
&W_{\mathfrak{n}+1} = e^{\frac{\tau}{2} \epsilon^2 \mathcal{L}} W_{\mathfrak{n},3},
\end{aligned}
\right.
\end{equation}
which can be written in the abstract form:
\begin{equation}\label{2.25}
W_{\mathfrak{n}+1} = S_{\frac{\tau}{2}}^{\mathcal{L}} \circ S_{\tau}^\mathcal{N} \circ S_{\frac{\tau}{2}}^{\mathcal{L}} W_{\mathfrak{n}}.
\end{equation}

The FR-MLFEM $\eqref{2.24}$ or \eqref{2.25} can be reformulated in the finite element form: Given $w_h^\mathfrak{n} = \Phi(x)^{\top} W_{\mathfrak{n}} \Phi(y)
\in Q_h^k$, we find $w_h^{\mathfrak{n}+1} = \Phi(x)^{\top} W_{\mathfrak{n}+1} \Phi(y) \in Q_h^k$ following
\begin{equation}\label{finiteform}
\left\lbrace
\begin{aligned}
&w_h^{\mathfrak{n},1} = e^{\frac{\tau}{2} \epsilon^2\mathcal{L}_h} w_h^\mathfrak{n},\\
&(w_h^{\mathfrak{n},2},v_h)_h = (w_h^{\mathfrak{n},1},v_h)_h + \tau (\mathcal{N}(w_h^{\mathfrak{n},1}) ,v_h)_h, & &\forall v_h \in Q_h^k,\\
&(w_h^{\mathfrak{n},3},v_h)_h = \frac{1}{2} (w_h^{\mathfrak{n},1},v_h)_h + \frac{1}{2} (w_h^{\mathfrak{n},2} + \frac{\tau}{2} \mathcal{N}(w_h^{\mathfrak{n},2}),v_h)_h,& &\forall v_h \in Q_h^k,\\
&w_h^{\mathfrak{n}+1} = e^{\frac{\tau}{2} \epsilon^2\mathcal{L}_h} w_h^{\mathfrak{n},3},
\end{aligned}
\right.
\end{equation}
where $ e^{\frac{\tau}{2} \epsilon^2 \mathcal{L}_h} w_h= \Phi(x)^{\top} e^{\frac{\tau}{2} \epsilon^2 \mathcal{L}} W \Phi(y)$, and $w_h^{\mathfrak{n},l} =\Phi(x)^{\top} W_{\mathfrak{n},l} \Phi(y) \in Q_h^k$ for $l=1,2,3$.

\begin{rem}\label{fullcpx}
    The computational complexity of the FR-MLFEM \eqref{2.24} or \eqref{finiteform} is $\mathcal{O}(mn)$.
\end{rem}

To present the main properties of the solution obtained from the system \eqref{2.21}, we introduce the following result.
\begin{lem}\label{lem2.3}
For any $\tau >0$, it holds
\[(e^{\tau \epsilon^2 \mathcal{L}} Z, \mathbf{I})_{\rm M}= (Z,\mathbf{I})_{\rm M}, \quad \forall Z \in \mathbb{R}^{m \times n}.\]
\end{lem}
The proof of \Cref{lem2.3} is presented in \Cref{lem2.2proof}. 
Then, we have the following results.
\begin{lem}\label{lem2.4}
Let $W_{\mathfrak{n}} \in \mathbb{R}^{m\times n}$ be the solution to the system $\eqref{2.24}$ corresponding to the conservative AC equation, and $w_h^{\mathfrak{n}} = \Phi(x)^\top W_{\mathfrak{n}} \Psi(y) \in Q_h^k$ be the finite element solution. 
Then the solution conserves the mass
\[(w_h^{\mathfrak{n}},1)_h  = (W_{\mathfrak{n}}, \mathbf{I})_{\rm M} = (W_{0}, \mathbf{I})_{\rm M} = (w_h^0,1)_h. \]
\end{lem}

\begin{proof}
We prove this by mathematical induction. Assume that $(w_h^\mathfrak{n},1)_h = (w_h^0, 1)_h$. 
Next, we take the $\rm M$-weighted Frobenius inner product on both sides of  $\eqref{2.24}$ with $\mathbf{I}$. Using $\eqref{2.27}$ and Lemma \ref{lem2.3} yield
\[
\begin{aligned}
(w^{\mathfrak{n}+1}_h,1)_h & = (W_{\mathfrak{n}+1}, \mathbf{I})_{\rm M} = ( e^{\tau \epsilon^2 \mathcal{L}}W_{\mathfrak{n},3},\mathbf{I})_{\rm M} = (W_{\mathfrak{n},3},\mathbf{I})_{\rm M} \\
&= \frac{1}{2} (W_{\mathfrak{n},1},\mathbf{I})_{\rm M} + \frac{1}{2} (W_{\mathfrak{n},2},\mathbf{I})_{\rm M} +\frac{1}{2} \tau(\mathcal{N}(W_{\mathfrak{n},2}),\mathbf{I})_{\rm M}  =  (W_{\mathfrak{n},1},\mathbf{I})_{\rm M}+\frac{1}{2} \tau(\mathcal{N}(W_{\mathfrak{n},1}),\mathbf{I})_{\rm M}  \\
& = (W_{\mathfrak{n},1},\mathbf{I})_{\rm M} = (e^{\frac{\tau}{2} \epsilon^2 \mathcal{L}} W_{\mathfrak{n}}, \mathbf{I})_{\rm M} = (W_{\mathfrak{n}}, \mathbf{I})_{\rm M} = (w_h^{\mathfrak{n}},1)_h.
\end{aligned} 
\]
This completes the proof.
\end{proof}
Similar to \cite{li2022stability,li2022stability2}, we introduce an $O(\tau)$ modification to the original energy \eqref{Oener}, leading to the following modified energy
\begin{equation}\label{2.28}
\widetilde{E}^{\fn} = \frac{1}{2\tau} ((e^{-\tau \epsilon^2 \mathcal{L}_h}-1)w_h^{\mathfrak{n},1}, w_h^{\mathfrak{n},1})_h + (G(w_h^{\mathfrak{n},1}),1)_h,
\end{equation}
where $G(w_h^{\mathfrak{n},1}) = F(0) + \int_{0}^{w_h^{\mathfrak{n},1}} g(s) \, ds$ with
\begin{equation}\label{gweqn}
g(s) = - \frac{1}{2} \mathcal{N}(s) - \frac{1}{2} \mathcal{N}(s+\tau \mathcal{N}(s)).
\end{equation}
For $g(s)$ defined in \eqref{gweqn}, it holds
$$
|g^{\prime}(s)| \leq C_g,
$$
where the constant 
\[C_g = \max\limits_{\|s\| \leq 1} |g^{\prime}(s)| \leq 2 \max\limits_{\|s\| \leq 1} |\mathcal{N}^{\prime}(s)|. \]
For the classical AC equation, $C_g = \frac{1}{2}$, and for the conservative AC equation, $C_g \leq C_{|\Omega|}$, where $C_{|\Omega|}$ is a constant depending on $|\Omega|$ \cite{jiang2022unconditionally}.
\begin{lem}\label{lem2.5}
For both the classical AC equation \eqref{1.1} and the conservative AC equation \eqref{1.2}, if the time step $\tau$ satisfies $\tau \in (0,\frac{1}{C_{g}}]$, the solution of system $\eqref{2.21}$ or \eqref{finiteform} satisfies the following energy law
\[\widetilde{E}^{\mathfrak{n}+1} \leq \widetilde{E}^{\mathfrak{n}}.\]
\end{lem}
The proof of \Cref{lem2.4} is presented in \Cref{lem2.5proof}.

\section{Dynamical low-rank mass-lumped finite element method}\label{Sec3}

Let $\mathcal{M}_r \subset \mathbb{R}^{m \times n}$ be the manifold of rank-$r$ matrices ($r \leq \min(m, n)$). 
For convenience, let $Z \in \mathbb{R}^{m \times n}$ represent the coefficient matrix of $z_h \in Q_h^k$. In the following, we use script letters to denote low-rank matrix approximations (e.g., $\mathcal{Z}$ represents the low-rank approximation of $Z$). 
Similar to the full-rank strategy, we continue to use the spitting \eqref{2.20} to solve the matrix differential equation \eqref{2.19} in the low-rank approximation. We begin with the low-rank approximation for the linear matrix differential equation \eqref{2.20a}.

\subsection{An exact low-rank integrator for the linear problem \eqref{2.20a}} 
Given $\tau>0$ and a rank-$r$ approximation 
$\mathcal{X}_{\mathfrak{n}} = U_{\mathfrak{n}}S_{\mathfrak{n}}V_{\mathfrak{n}}^{\top}$ at $t_\mathfrak{n}$ with the factors $U_{\mathfrak{n}} \in \mathbb{R}^{m\times r}$ and $V_{\mathfrak{n}} \in \mathbb{R}^{n\times r}$ satisfying
\[U^{\top}_{\mathfrak{n}} M_x U_{\mathfrak{n}} = V^{\top}_{\mathfrak{n}} M_y V_{\mathfrak{n}} = I_r,\]
then the linear matrix differential equation \eqref{2.20a} at $t_{\mathfrak{n}+1}=t_\mathfrak{n} + \tau$ can be solved analytically by \eqref{2.21} as the rank-$r$ matrix
\begin{equation}\label{linlrn1}
\mathcal{X}_{\mathfrak{n}+1} = e^{\tau \epsilon^2 L_x} \mathcal{X}_{\mathfrak{n}} e^{\tau \epsilon^2 L^{\top}_y}= e^{\tau \epsilon^2 L_x} U_{\mathfrak{n}} S_{\mathfrak{n}} V_{\mathfrak{n}}^{\top} e^{\tau \epsilon^2 L^{\top}_y}=U_{\mathfrak{n}+1} S_{\mathfrak{n}+1} V_{\mathfrak{n}+1}^{\top}, 
\end{equation}
where the last equality is obtained by \Cref{lralg}, and 
\begin{equation}\label{morth}
U_{\mathfrak{n}+1}^\top M_x U_{\mathfrak{n}+1} = V_{\mathfrak{n}+1}^\top M_y V_{\mathfrak{n}+1}= I_r.
\end{equation}

\begin{breakablealgorithm}
    \normalsize
    \medskip
    \caption{An exact low-rank integrator to the linear matrix differential equation \eqref{2.20a}. }\label{lralg}
    \begin{algorithmic}
    \textbf{Input:} $U_{\fn}, S_{\fn}, V_{\fn}$, $\epsilon, \tau$, $M_x, M_y$.\\
    \textbf{Output:} $U_{\fn+1},S_{\fn+1},V_{\fn+1}$.\\
    \begin{itemize}
         \item \textbf{Step 1:} Update $U^{\fn} \rightarrow  U^{\fn+1} \in \mathbb{R}^{m \times {r}}$ and $V_{\fn} \rightarrow  V_{\fn+1}\in \mathbb{R}^{n\times {r}}$ in parallel: perform a generalized QR (GQR) decomposition \cite{yin2023semi} with weight $M_x$ and $M_y$, respectively,
            \[
            \begin{aligned}
                \left[ U_{\mathfrak{n}+1}, R \right] &= GQR\left( e^{\tau \epsilon^2 L_x} U_{\mathfrak{n}}, M_x \right), \\
                \left[ V_{\mathfrak{n}+1}, P \right] &= GQR\left( e^{\tau \epsilon^2 L_y} V_{\mathfrak{n}}, M_y \right),
            \end{aligned}
            \]
        where $U_{\mathfrak{n}+1} $ and $ V_{\mathfrak{n}+1}$ satisfy \eqref{morth}.\\
        \item{\textbf{Step 2:}} Update $S_{\mathfrak{n}} \rightarrow S_{\mathfrak{n}+1} \in \mathbb{R}^{r\times r}$ by
             \begin{equation}\label{linS}
                 S_{\mathfrak{n}+1} = R S_\mathfrak{n} P^{\top}.
             \end{equation}\\
    \end{itemize}
\end{algorithmic}
\end{breakablealgorithm}

\begin{lem}\label{lem3.1}
Given low-rank matrix $\mathcal{X}_{\mathfrak{n}} = U_{\mathfrak{n}}S_{\mathfrak{n}}V_{\mathfrak{n}}^{\top}$ at $t_\mathfrak{n}$, the low-rank approximation  $\mathcal{X}_{\mathfrak{n}+1} = U_{\mathfrak{n}+1} S_{\mathfrak{n}+1} V_{\mathfrak{n}+1}^{\top}$, obtained in \eqref{linlrn1} for the linear matrix differential equation \eqref{2.20a} at $t_{\mathfrak{n}+1}$, satisfies (i) the rank of matrix $\mathcal{X}_{\mathfrak{n}+1}$ is no more than that of $\mathcal{X}_{\mathfrak{n}}$, and (ii)
$(U_{\mathfrak{n}+1} S_{\mathfrak{n}+1} V_{\mathfrak{n}+1}^{\top}, \mathbf{I})_{\rm M} = (U_{\mathfrak{n}} S_{\mathfrak{n}} V_{\mathfrak{n}}^{\top}, \mathbf{I})_{\rm M}$.
\end{lem}
\begin{proof} 
(i) follows from \eqref{linS}. For (ii), it holds
\[(U_{\mathfrak{n}+1} S_{\mathfrak{n}+1} V_{\mathfrak{n}+1}^{\top}, \mathbf{I})_{\rm M} = ( U_{\mathfrak{n}+1} (R S_{\mathfrak{n}} P^{\top}) V_{\mathfrak{n}+1}^{\top}, \mathbf{I})_{\rm M} = (e^{\tau \epsilon^2 L_x} U_{\mathfrak{n}} S_{\mathfrak{n}}(e^{\tau \epsilon^2 L_y} V_{\mathfrak{n}})^{\top}, \mathbf{I})_{\rm M} = (U_{\mathfrak{n}} S_{\mathfrak{n}} V_{\mathfrak{n}}^{\top}, \mathbf{I})_{\rm M},
\]
where we have used \Cref{lem2.3} for the last equality.
\end{proof}

\subsection{DLRA for nonlinear problem \eqref{2.20b}}
The Dynamical Low-Rank Approximation (DLRA) is traditionally formulated by evolving the matrix differential equation through a minimization problem using the Frobenius norm on the tangent space of $\mathcal{M}_r$ at $Z$ \cite{koch2007dynamical}. 
In the following, to maintain formulation equivalence between the Galerkin equation of the DLRA and the matrix variational problem \eqref{matvar}, we introduce a weighted DLRA, in which the minimization problem is formulated using the $\rm M$-weighted Frobenius norm in \Cref{frobdef}. A similar weighted DLRA was investigated in \cite{yin2023semi}.
 
\begin{definition}\label{define3.1}
The weighted DLRA to $\eqref{2.20b}$ is defined as the solution $\mathcal{Z} \in \mathcal{M}_r$ (where $\mathcal{Z}$ is a rank-r approximation of $Z$) to the differential equation
\begin{equation}\label{3.1}
\partial_t \mathcal{Z} = \mathop{\arg\min}\limits_{\delta \mathcal{Z} \in \mathcal{T}_{\mathcal{Z} } \mathcal{M}_r} \|\delta \mathcal{Z} - \mathcal{N}(\mathcal{Z}) \|_{\rm M},
\end{equation}
where $\mathcal{N}(\mathcal{Z})$ is the nonlinear term defined in $\eqref{2.20b}$, with the initial condition $\mathcal{Z}(0) = \mathcal{Z}_0$ being a rank-$r$ approximation of $Z_\mathfrak{n}$.
$\mathcal{T}_{\mathcal{Z}} \mathcal{M}_r$ is the tangent space of $\mathcal{M}_r$ at $\mathcal{Z}$.
\end{definition}

Assume that matrix $\mathcal{Z} \in \mathcal{M}_r$ has a rank-$r$ decomposition
\begin{equation}\label{3.2}
\mathcal{Z} = U S V^{\top}, \text{ where } U^{\top} M_x U = V^{\top}M_y V = I_r,
\end{equation}
with $U \in \mathbb{R}^{m \times r}, S \in \mathbb{R}^{r\times r} $and $V \in \mathbb{R}^{n \times r}$. Then the tangent space of $\mathcal{M}_r$ at $\mathcal{Z}$ is given by \cite{koch2007dynamical}:
\begin{equation}\label{3.3}
\mathcal{T}_{\mathcal{Z}} \mathcal{M}_r = \{ \delta U S V^{\top} + U \delta S V^{\top} + U S \delta V^{\top} : U^{\top} M_x\delta U  =0, \; V^{\top} M_y \delta V=0.\},
\end{equation}
where $\delta U \in \mathbb{R}^{m \times r}, \delta S \in \mathbb{R}^{r\times r} $and $\delta V \in \mathbb{R}^{n \times r}$. Under the gauge conditions $U^{\top} M_x \delta U = V^{\top} M_y\delta V =0$, each matrix $\delta \mathcal{Z} \in \mathcal{T}_{\mathcal{Z}} \mathcal{M}_r$ has a unique decomposition
\begin{equation}\label{3.4}
\delta \mathcal{Z} = \delta U S V^{\top} + U \delta S V^{\top} + US\delta V^{\top} = P_U^{\perp} M_x \delta \mathcal{Z} M_y P_V + P_U M_x \delta \mathcal{Z} M_y P_V + P_U M_x \delta \mathcal{Z}^{\top} M_y P_V^{\perp},
\end{equation}
where 
\begin{equation}\label{3.5}
\delta U = P_U^{\perp} M_x\delta \mathcal{Z} M_y V S^{-1}, \delta S = U^{\top} M_x \delta \mathcal{Z}M_y V, \text{ and } \delta V = P_V^{\perp} M_y \delta \mathcal{Z}^{\top}M_x U S^{-\top},
\end{equation}
with the symmetric matrices
\begin{equation}\label{3.6}
\begin{aligned}
P_U = U U^{\top}, &\quad P_{U}^{\perp} =M_x^{-1} - P_U,\\
P_V = V V^{\top},&\quad P_{V}^{\perp} = M_y^{-1} - P_V.  
\end{aligned}
\end{equation}
Here, the matrix $P_U M_x$ (resp. $P_V M_y$) is the orthogonal projection onto the column space of $U$ (resp. $V$)  with respect to the inner product on $\mathbb{R}^m$ (resp. $\mathbb{R}^n$) with weight $M_x$ (resp. $M_y$).

\begin{prop}\label{prop:DLRA}
Suppose the solution $\mathcal{Z} = USV^{\top} \in \mathcal{M}_r$ satisfies Definition \ref{define3.1}, with initial data $\mathcal{Z}(0) = U_0 S_0 V_0^{\top} \in \mathcal{M}_r$ where $U_0^{\top} M_x U_0 = V_0^{\top} M_y V_0 = I_r$. Then, the following formulations are equivalent to the weighed DLRA defined in \Cref{define3.1} \cite{koch2007dynamical, yin2023semi}.
\begin{itemize}
\item[(i)] Galerkin condition: $\partial_t \mathcal{Z} \in \mathcal{T}_{\mathcal{Z}} \mathcal{M}_r$ is the solution of the Galerkin condition
\begin{equation}\label{3.7}
\begin{aligned}
&(\partial_t \mathcal{Z} - \mathcal{N}(\mathcal{Z}),\delta \mathcal{Z})_{\rm M} = 0 ,\quad \forall \delta \mathcal{Z} \in \mathcal{T}_{\mathcal{Z}}\mathcal{M}_r, \\
&\mathcal{Z}(0) = U_{0} S_{0} V_{0}^{\top}.
\end{aligned}
\end{equation}
\item[(ii)] Equations of motion: the factors of $\mathcal{Z}$ satisfy
\begin{equation}\label{3.8}
\begin{aligned}
&\dot{U} = P_{U}^{\perp} M_x \mathcal{N}(\mathcal{Z}) M_y V S^{-1}, \ \dot{S} = U^{\top} M_x \mathcal{N}(\mathcal{Z}) M_y V, \ \dot{V} = P_{V}^{\perp} M_y \mathcal{N}(\mathcal{Z})^{\top}  M_x U S^{-\top}, \\
&U(0) = U_{0},\; S(0) = S_{0},\; V(0) = V_{0},\;
\end{aligned}
\end{equation}
where $P_U$ and $P_V$ are defined in $\eqref{3.6}$.

\item[(iii)] Coupled ordinary differential equations (ODEs), referred to as the KLS system: the matrices $K = US \in \mathbb{R}^{m \times r}$, $L = V S^{\top} \in \mathbb{R}^{n \times r}$, and $S\in \mathbb{R}^{r \times r}$ satisfy
\begin{equation}\label{3.9}
\begin{aligned}
&\dot{K} =  \mathcal{N}(KV^{\top}) M_y V, \quad \dot{L} = \mathcal{N}(UL^{\top})^{\top} M_x U, \quad \dot{S} = U^{\top} M_x \mathcal{N}(USV^{\top}) M_y V, \\
&K(0) = U_{0}S_{0},\; L(0) = V_{0}S_{0}^{\top},\; S(0) = S_{0}.
\end{aligned}
\end{equation}
\end{itemize}
\end{prop}
The proof of \Cref{prop:DLRA} is similar to that of \cite[Proposition 3.5]{yin2023semi}.

\subsection{Low complexity strategy for the nonlinear and nonlocal terms}

To reduce the computational cost of integrating the differential equation \eqref{3.1} in DLRA, or its equivalent formulations \eqref{3.7}-\eqref{3.9}, for the nonlinear matrix differential equation \eqref{2.20b}, it is essential to minimize the computational complexity of both the nonlinear terms and any nonlocal terms, if present, in $\mathcal{N}(W)$. The term $\mathcal{N}(W)$ includes a Hadamard power of $W$ up to degree $3$ as described in \eqref{matrix-nonlinear}. Directly forming $W = USV^\top$ and then computing $W^{\circ 3}$ significantly increases the computational complexity, making it comparable to that of the full-rank scheme. Therefore, an efficient low-rank technique for handling this nonlinear term is crucial.

The coupled ODE system \eqref{3.9} implies that we can compute the product involving $\mathcal{N}(W)$ without explicitly computing $\mathcal{N}(W)$ itself. For example, in the differential equation of $K$, we can compute the product $\mathcal{N}(W) M_y V$ for any $V \in \mathbb{R}^{n \times r}$ without explicitly computing the elements of $W$. 
We only consider the nonlinear term $W^{\circ 3} M_y V$ for simplicity.
Recall that $K = U S \in \mathbb{R}^{m \times r}$. Then, it holds
\begin{equation}
W = USV^{\top} = K V^{\top} = \sum\limits_{i=1}^{r} K_i V_i^{\top},    
\end{equation}
where $K_i$, $V_i$ are the i-th column of $K$ and $V$, respectively. 
In addition to the properties in \Cref{hadadef}, the following property holds for the Hadamard product,
\[(K_i V_i^{\top}) \odot (K_j V_j^{\top}) = (K_i\odot K_j)(V_i \odot V_j)^{\top}, \quad \text{for} 1\leq i,j\leq r.\]
Hence, the cubic nonlinearity which appears in the nonlinear term and $V \in \mathbb{R}^{n \times r}$ it holds:
\begin{equation}\label{computationcost}
\begin{aligned}
(W^{\circ 3}) M_y V & = \left( \left( \sum\limits_{i=1}^{r} K_i V_i^{\top}\right) \odot \left( \sum\limits_{j=1}^{r} K_j V_j^{\top}\right) \odot \left( \sum\limits_{l=1}^{r} K_l V_l^{\top}\right)\right) M_y  V \\
& = \left(  \sum\limits_{i,j,l=1}^{r} (K_i \odot K_j \odot K_l) \left( V_i \odot V_j \odot V_l\right)^{\top} \right) M_y V \\
& =  \sum\limits_{i,j,l=1}^{r} (K_i \odot K_j \odot K_l)  \left(\left( V_i \odot V_j \odot V_l\right)^{\top} M_y  V\right). 
\end{aligned}
\end{equation}
Therefore, the computation of $(W^{\circ 3}) M_y V$, when implemented following \eqref{computationcost}, results in a computational complexity of $\mathcal{O}((m+nr)r^{3})$.
Analogously, the computational complexity of $(W^{\circ 3})^{\top} M_x U$ with $U \in \mathbb{R}^{m \times r}$ is $\mathcal{O}((mr+n)r^{3})$.
By \eqref{computationcost}, the dominant term in the product $U^{\top} M_x \mathcal{N}(USV^{\top}) M_y V$ has the following reformulation
\begin{equation}
U^{T} M_x (W^{\circ 3}) M_y V =  \sum\limits_{i,j,l=1}^{r} \left(U^{T} M_x(K_i \odot K_j \odot K_l) \right)  \left(\left( V_i \odot V_j \odot V_l\right)^{\top} M_y  V\right), 
\end{equation}
which has the computational complexity $\mathcal{O}((m+n)r^{4})$. 

Next, we consider the low-rank strategy for the nonlocal term $(W^{\circ 3},\mathbf{I})_{\rm M}$ in $\eqref{matrix-nonlinearb}$ and calculate its computational complexity. Recall that the all-ones matrix \(\mathbf{I} = q_m q_n^\top \in \mathbb{R}^{m \times n}\), where \(q_m \in \mathbb{R}^m\) and \(q_n \in \mathbb{R}^n\) are all-ones vectors. 
We can reduce its computational complexity by avoiding direct computation and using the following reformulation instead:
\[(W^{\circ 3},\mathbf{I})_{\rm M}=(M_x W^{\circ 3} M_y, q_m q_n^{\top})_{\rm F} = (M_x W^{\circ 3} M_y q_n, q_m  )_{\rm F}.\]
Since \(M_x\) is a diagonal matrix, the computational complexity of computing \(M_x W^{\circ 3} M_y q_n\), which dominates the computational complexity of $(W^{\circ 3},\mathbf{I})_{\rm M}$, is equal to that of computing \(W^{\circ 3} M_y q_n\). The decomposition of \(W^{\circ 3} M_y q_n\) follows the same as \eqref{computationcost}, with \(V\) replaced by \(q_n\), and its computational complexity is \(\mathcal{O}((m+n)r^{3})\), which is also the computational complexity of $(W^{\circ 3},\mathbf{I})_{\rm M}$.
Similarly, the computational complexity for the nonlocal term $(W^{\circ 2},\mathbf{I})_{\rm M}$ in $\eqref{matrix-nonlinearb}$ is $\mathcal{O}((m+n)r^{2})$.

Based on the discussion above, we summarize the following result.
\begin{lem}\label{lem3.2}
    The computational complexity of evaluating all the products $\mathcal{N}(KV^{\top}) M_y V$, $\mathcal{N}(UL^{\top})^{\top}\\ M_x U$, and $U^{\top} M_x \mathcal{N}(USV^{\top}) M_y V$ in the coupled system \eqref{3.12}, corresponding to the nonlinear matrix differential equation \eqref{2.20b} with either nonlinear term in \eqref{matrix-nonlinear}, is \(\mathcal{O}((m+n)r^{4})\).
\end{lem}

\subsection{Second-order integrator}
Similar to the Strang splitting approach for the FR-MLFEM \eqref{2.25}, we propose a second-order low-rank integrator for the matrix differential equation \eqref{2.19} by solving the split linear and nonlinear matrix differential equations \eqref{2.20a} and \eqref{2.20b}. For the linear component \eqref{2.20a}, we use the exact low-rank integrator in \Cref{lralg}. 

For the nonlinear component \eqref{2.20b}, we employ the augmented BUG integrator \cite{ceruti2022rank, ceruti2024robust} for the KLS system \eqref{3.9}, an equivalent formulation of DLRA to the nonlinear matrix differential equation $\eqref{2.20b}$. 
This method extends the bases used in the S-step of the BUG integrator \cite{ceruti2022unconventional, ceruti2024parallel} to eliminate unwanted projection errors. 
The BUG integrator, viewed as a splitting method, is applied to the KLS system $\eqref{3.9}$, where the $K$ and $L$ equations are updated in parallel, followed by an update with the $S$ equation. 

The second-order augmented BUG integrator \cite{ceruti2024robust} builds upon the first-order augmented BUG integrator \cite{ceruti2022unconventional}. We start by outlining the key aspects of the first-order integrator for the matrix differential equation $\eqref{2.20b}$.
Given a factored rank-\( r \) matrix \(\mathcal{Z}_\mathfrak{n} = U_\mathfrak{n} S_\mathfrak{n} V_\mathfrak{n}^{\top}\) at $t_\mathfrak{n}$, with factors $U_{\mathfrak{n}} \in \mathbb{R}^{m\times r}$ and $V_{\mathfrak{n}} \in \mathbb{R}^{n\times r}$ satisfying  
\[
U_{\mathfrak{n}}^{\top} M_x U_{\mathfrak{n}} = V_{\mathfrak{n}}^{\top} M_y V_{\mathfrak{n}} = I_r,
\]
the one-step first-order augmented BUG integrator gives a low-rank approximation at the next time step \(t_{\mathfrak{n}+1} = t_\mathfrak{n} + \tau\), namely,
\begin{equation}\label{3.10}
\mathcal{Z}_{\mathfrak{n}+1} = U_{\mathfrak{n}+1} S_{\mathfrak{n}+1} V_{\mathfrak{n}+1}^{\top},
\end{equation}
where the updated bases $U_{\mathfrak{n}+1} \in \mathbb{R}^{m \times \hat{r}}$ and $V_{\mathfrak{n}+1} \in \mathbb{R}^{n \times \hat{r}}$ satisfy  
\[
U^{\top}_{\mathfrak{n}+1} M_x U_{\mathfrak{n}+1} = V^{\top}_{\mathfrak{n}+1} M_y V_{\mathfrak{n}+1} = I_{\hat{r}}.
\]
The specifics of the first-order augmented BUG integrator are outlined in \Cref{firstBUG}.

\begin{breakablealgorithm}
    \caption{One-step first-order augmented BUG integrator for \eqref{2.20b}. }\label{firstBUG}
    \begin{algorithmic}
    \textbf{Input:} {$U_{\fn}, S_{\fn}, V_{\fn}$, $\tau$, $M_x, M_y$.}\\
    \textbf{Output:}{$U_{\fn+1},S_{\fn+1},V_{\fn+1}$. }
    \begin{itemize}
    \item{\textbf{Step 1:}} Update $U_{\fn} \rightarrow  U_{\fn+1} \in \mathbb{R}^{m \times \hat{r}}$ and $V_{\fn} \rightarrow  V_{\fn+1}\in \mathbb{R}^{n\times \hat{r}}$ in parallel, with $(\hat{r} = 2r)$:\\
        \begin{itemize}
            \item{\texttt{$K$-step}}: 
            \begin{itemize}
                \item Solve $K_{\fn+1}$ from the $m\times {r} $ matrix equation
                    \begin{equation}\label{3.12}
                    \dfrac{K_{\mathfrak{n}+1}-K_{\mathfrak{n}}}{\tau} = \mathcal{N}(K_\mathfrak{n}V_\mathfrak{n}^{\top})M_y V_\mathfrak{n}, \quad \quad K_\mathfrak{n} = U_\mathfrak{n} S_\mathfrak{n}.
                    \end{equation}\\
                \item Set $\widetilde{K}_{\mathfrak{n}+1} = [K_{\mathfrak{n}+1}, U_\mathfrak{n}]$ and compute \([U_{\mathfrak{n}+1}, R_K] = GQR(\widetilde{K}_{\mathfrak{n}+1}, M_x)\).\\
                \item Compute the $\hat{r} \times r$ matrix $M = U_{\mathfrak{n}+1}^{\top} M_x U_\mathfrak{n}.$\\
            \end{itemize}
        \end{itemize}
        \begin{itemize}
            \item{\texttt{$L$-step}}: 
            \begin{itemize}
                \item Solve $L_{\fn+1}$ from the $n \times {r}$ matrix equation
                \begin{equation}\label{3.13}
                \dfrac{L_{\mathfrak{n}+1}-L_{\mathfrak{n}}}{\tau} = \mathcal{N}(U_\mathfrak{n} L_\mathfrak{n}^{\top})^{\top} M_x U_\mathfrak{n}, \quad L_\mathfrak{n} = V_\mathfrak{n} S_\mathfrak{n}^{\top}.
                \end{equation}
                \item Set $\widetilde{L}_{\mathfrak{n}+1} = [L_{\mathfrak{n}+1}, V_\mathfrak{n}]$ and compute $[V_{\mathfrak{n}+1}, R_L]=GQR(\widetilde{L}_{\mathfrak{n}+1}, M_y)$.
                \item Compute the $\hat{r} \times r$ matrix $N = V_{\mathfrak{n}+1}^{\top}M_y V_\mathfrak{n}.$
            \end{itemize}
        \end{itemize}
    \item{\textbf{Step 2:}} Update $S_{\fn}\rightarrow S_{\fn+1} \in \mathbb{R}^{\hat{r} \times \hat{r}}$:
        \begin{itemize}
        \item{\texttt{$S$-step}}:
            \begin{itemize}
                \item Set $S_{\mathfrak{n},*} = M S_\mathfrak{n} N^{\top}.$
                \item Solve the $\hat{r} \times \hat{r}$ matrix equation 
                    \begin{equation}\label{3.14}
                    \dfrac{S_{\mathfrak{n}+1} - S_{\mathfrak{n},*}}{\tau} = U_{\mathfrak{n}+1}^{\top} M_x \mathcal{N}(U_{\mathfrak{n}+1} S_{n,*} V_{\mathfrak{n}+1}^{\top}) M_yV_{\mathfrak{n}+1}.
                    \end{equation}
            \end{itemize}
        \end{itemize}
    \end{itemize}
\end{algorithmic}
\end{breakablealgorithm}
\begin{rem}
Since the first-order augmented BUG integrator in \Cref{firstBUG} is a component of the second-order integrator for the nonlinear matrix differential equation \eqref{2.20b}, the truncation step is deferred to the main algorithm for the matrix differential equation \eqref{2.20b}.
\end{rem}

With the preparation above, we present the second-order integrator for the matrix differential equations in \eqref{2.20}.
The exact low-rank integrator from \Cref{lralg} is employed for the linear matrix differential equation \eqref{2.20a}, while the second-order augmented BUG integrator \cite{ceruti2024robust} is applied to solve the nonlinear matrix differential equation \eqref{2.20b}. In particular, we will use SSP-RK2 to solve the matrix equation in the S-step.
This approach utilizes a second-order Strang splitting method, similar to \eqref{2.25} or \eqref{2.24}, as implemented in FR-MLFEM, to integrate the linear and nonlinear integrators. The details are outlined in \Cref{lrnonlinear2}.

\begin{breakablealgorithm}
    \caption{A second-order integrator for the split matrix differential equations in \eqref{2.20}. }\label{lrnonlinear2}
    \begin{algorithmic}
    \textbf{Input:} $U_{\fn}, S_{\fn}, V_{\fn}$, $\tau$, $M_x, M_y$.\\
    \textbf{Output:} $U_{\fn+1},S_{\fn+1},V_{\fn+1}$. 
    \begin{itemize}
        \item{\textbf{Linear Step:}} Perform \Cref{lralg} with a time step $\frac{\tau}{2}$ and initial values $U_{\fn}, S_{\fn}, V_{\fn}$ to obtain 
        $$
        \mathcal{W}_{\mathfrak{n},1} = U_{\mathfrak{n},1} S_{\mathfrak{n},1} V_{\mathfrak{n},1}^{\top}.
        $$
    	\item{\textbf{Nonlinear Step}:}
    	\begin{itemize}
        \item{\texttt{Trapezoidal Approximation}}:
        \begin{itemize}
        \item  Perform \Cref{firstBUG} with a time step $\tau$ and initial values $U_{\mathfrak{n},1}, S_{\mathfrak{n},1}, V_{\mathfrak{n},1}$ to obtain the approximation of rank $\hat{r} \leq 2r$,
        \begin{equation}\label{3.20} 
        \mathcal{W}_{\fn,2} = U_{\mathfrak{n},2} S_{\mathfrak{n},2} V_{\mathfrak{n},2}^{\top}.   
\end{equation}
        \end{itemize}
         \item{\texttt{Galerkin Step}}:
         \begin{itemize}
        \item Assemble the augmented matrices $\mathbb{U}$ and $\mathbb{V}$ with rank $\bar{r} \leq 4r+1$:
\[
\begin{aligned}
&\mathbb{U} = [q_m,U_{\mathfrak{n},1}, \tau \mathcal{N}(\mathcal{W}_{\mathfrak{n},1})V_{\mathfrak{n},1}, \tau \mathcal{N}(\mathcal{W}_{\mathfrak{n},2})V_{\mathfrak{n},2}], \\
&\mathbb{V} = [q_n,V_{\mathfrak{n},1}, \tau \mathcal{N}(\mathcal{W}_{\mathfrak{n},1})^{\top} U_{\mathfrak{n},1}, \tau \mathcal{N}(\mathcal{W}_{\mathfrak{n},2})^{\top} U_{\mathfrak{n},2}].
\end{aligned}
\]

\item 
Compute
\begin{equation}\label{3.21}
\begin{aligned}
\overline{U}_{\fn,3} = \text{GQR}(\mathbb{U},M_x), \quad 
\overline{V}_{\fn,3} = \text{GQR}(\mathbb{V},M_y).
\end{aligned}
\end{equation}
\item Compute the $\bar{r} \times r$ matrices
$M = \overline{U}_{\fn,3}^{\top} M_x U_{\mathfrak{n},1}$ and $N = \overline{V}^{\top}_{\fn,3} M_y V_{\fn,1}$, and $\bar{r} \times \bar r$ matrix
\begin{equation}\label{Sn1proj}
    \overline{S}_{\mathfrak{n},1} = M S_{\mathfrak{n},1} N^{\top}.
\end{equation}
\item Solve the $\bar{r} \times \bar r$ matrix equations by the SSP-RK2:
\begin{subequations}\label{3.22}
\begin{align}
\overline{S}_{\mathfrak{n},2} &= \overline{S}_{\mathfrak{n},1} + \tau \overline{U}^{\top}_{\fn,3} M_x \mathcal{N}(\overline{U}_{\fn,3}  \overline{S}_{\mathfrak{n},1} \overline{V}_{\fn,3}^{\top})M_y \overline{V}_{\fn,3}, \label{3.22a}\\
\overline{S}_{\mathfrak{n},3} &= \frac{1}{2} \overline{S}_{\mathfrak{n},1} + \frac{1}{2} \overline{S}_{\fn,2} + \frac{\tau }{2} \overline{U}_{\fn,3}^{\top} M_x\mathcal{N}(\overline{U}_{\fn,3} \overline{S}_{\mathfrak{n},2} \overline{V}_{\fn,3}^{\top})M_y \overline{V}_{\fn,3}.\label{3.22b}
\end{align}
\end{subequations}
\end{itemize}
\end{itemize} 
    \begin{itemize}
    \item{\texttt{Truncation Step}}:
    \begin{itemize} 
\item Compute the SVD decomposition of $\overline{S}_{\mathfrak{n}, 3} = \overline{R}  \; \overline{\Sigma} \; \overline{P}^{\top}$ with $\overline{\Sigma} = diag(\sigma_i)$. 
\item Truncate $\overline{\Sigma}$ to $S_{\mathfrak{n},3}$ with either the original rank $r$ or a new rank $r=\tilde r$, determined by prescribed a truncation tolerance $\eta$, such that
\begin{equation}\label{3.23}
\left( \sum\limits_{i=r +1}^{\bar{r}} \sigma_i^2 \right)^{\frac{1}{2}} \leq \eta.
\end{equation}

\item Set $R \in \mathbb{R}^{\bar{r} \times r}$ and $P \in \mathbb{R}^{\bar{r} \times r}$ containing the first $r$ columns of $\overline{R}$ and $\overline{P}$, respectively. 
\item Set $U_{\fn,3} = \overline{U}_{\fn,3} R \in \mathbb{R}^{m \times r}$ and $V_{\fn,3} = \overline{V}_{\fn,3} P \in \mathbb{R}^{n \times r}.$
\end{itemize}
\end{itemize}
\item{\textbf{Linear Step:}} Perform \Cref{lralg} with a time step $\frac{\tau}{2}$ and initial values $U_{\mathfrak{n},3}, S_{\mathfrak{n},3}, V_{\mathfrak{n},3}$ \\to obtain 
 $$
        \mathcal{W}_{\mathfrak{n}+1} = U_{\mathfrak{n}+1} S_{\mathfrak{n}+1} V_{\mathfrak{n}+1}^{\top}.
        $$
\end{itemize}
\end{algorithmic}
\end{breakablealgorithm}

The computational complexity of \Cref{lrnonlinear2} is primarily determined by the nonlinear matrix products, as described in \Cref{lem3.2}. Consequently, the following result holds.
\begin{lem}\label{lowcpx}
The computational complexity of \Cref{lrnonlinear2} is $\mathcal{O}((m+n)r^{4})$.
\end{lem}

\begin{rem}
\Cref{lowcpx}, combined with \Cref{fullcpx}, demonstrates that \Cref{lrnonlinear2} is computationally advantageous when $r^4 < \min\{m,n\}$, a condition typically satisfied when the solution exhibits a low-rank structure.
In future work, we will explore strategies, such as the discrete empirical interpolation method \cite{naderi2023adaptive}, to further reduce the computational complexity of the low-rank integrator.
\end{rem}

\begin{rem}\label{lr1st}
A first-order integrator can also be proposed by using the Lie-Trotter splitting \cite{xiao2017highly}, which combines \Cref{lralg} and \Cref{firstBUG} analogous to \Cref{lrnonlinear2}. 
The computational complexity for the first-order Algorithm is also $\mathcal{O}((m+n)r^{4})$.
We omit the details here. 
\end{rem}

\subsection{Second-order DLR-MLFEM}
Recall that script letters denote low-rank matrix approximations (e.g., $\mathcal{Z}$ represents the low-rank approximation of the matrix $Z$). 
From this subsection, we adopt bold italic font,  such as $\boldsymbol{z}_h=\Phi(x)^{\top} \mathcal{Z} \Psi(y) \in Q_h^k$, to represent the low-rank function associated with the low-rank matrix $\mathcal{Z}$.
First, we defined the following subspaces of the finite element space $Q_h^k$,
\begin{subequations}\label{3.24+}
\begin{align}
&\begin{aligned}
Q_{0}^{\fn,k} = \left\{ v ~| ~v(x,y) = \Phi(x)^{\top} U_\mathfrak{n} S V_\mathfrak{n}^{\top} \Psi(y), \quad \forall S \in \mathbb{R}^{r \times r} \right\},
\end{aligned} \label{3.24+a}\\
&\begin{aligned}
Q_{1}^{\fn,k} = \left\{ v ~| ~v(x,y) = \Phi(x)^{\top} U_{\mathfrak{n},1} S V_{\mathfrak{n},1}^{\top} \Psi(y), \quad \forall S \in \mathbb{R}^{r \times r} \right\},
\end{aligned} \label{3.24+b}\\
&\begin{aligned}
\overline{Q}^{\fn,k}_3 = \left\{ v ~| ~v(x,y) = \Phi(x)^{\top} \overline{U}_{\mathfrak{n},3} S \overline{V}_{\mathfrak{n},3}^{\top} \Psi(y), \quad \forall S \in \mathbb{R}^{\bar{r} \times \bar{r}} \right\},
\end{aligned} \label{3.24+d}\\
&\begin{aligned}
Q_{3}^{\fn,k} = \left\{ v ~| ~v(x,y) = \Phi(x)^{\top} U_{\mathfrak{n},3} S V_{\mathfrak{n},3}^{\top} \Psi(y), \quad \forall S \in \mathbb{R}^{\tilde r \times \tilde r} \right\},
\end{aligned} \label{3.24+c}
\end{align}
\end{subequations}
where $r$, $\bar{r}$, and $\tilde{r}$ are given in \Cref{lrnonlinear2}.
\begin{rem}
    We have used the adaptive rank $\tilde r$ for $S$ matrix in the subspace $Q_{3}^{\fn,k}$ in \eqref{3.24+c}, and it can be changed to fixed rank $r$ as stated in Truncation Step in \Cref{lrnonlinear2}. 
\end{rem}
By \eqref{3.21} in \Cref{lrnonlinear2}, it can be observed that 
\begin{equation}\label{Q1nk}
    Q_{1}^{\fn,k} \subset \overline{Q}^{\fn,k}_3.
\end{equation}

Similar to the full-rank finite element method, \Cref{lrnonlinear2} provides a second-order low-rank finite element approximation for the matrix differential equations in \eqref{2.20}. Notably, the low-rank finite element formulations in \Cref{lrnonlinear2} are mathematically equivalent to the full-rank finite element scheme \eqref{finiteform}, differing only in that the approximations are computed in subspaces of the original finite element space. Specifically, the following statement holds.
\begin{lem}\label{lrfem}
\Cref{lrnonlinear2} can be reformulated as the DLR-MLFEM: Given $\boldsymbol{w}_h^\mathfrak{n} = \Phi(x)^{\top} U_{\mathfrak{n}} S_{\mathfrak{n}} V_{\mathfrak{n}}^{\top} \Phi(y)
\in Q_{0}^{\fn,k}$, we find $ \boldsymbol{w}_h^{\mathfrak{n}+1} = \Phi(x)^{\top} U_{\mathfrak{n}+1} S_{\mathfrak{n}+1} V_{\mathfrak{n}+1}^{\top} \Phi(y) \in Q_{0}^{\fn+1,k}$ such that
\begin{subequations}\label{finiteform2}
\begin{align}
&\begin{aligned}
\boldsymbol{w}_S^{\mathfrak{n},1} = e^{\frac{\tau}{2} \epsilon^2 \mathcal{L}_h} \boldsymbol{w}_h^\mathfrak{n}
\end{aligned} \label{3.25+a}\\
&\begin{aligned}
(\boldsymbol{w}_{\overline{S}}^{\mathfrak{n},2},\boldsymbol{w}_2)_h = (\boldsymbol{w}_{S}^{\mathfrak{n},1},\boldsymbol{w}_2)_h + \tau (\mathcal{N}(\boldsymbol{w}_{S}^{\mathfrak{n},1}) ,\boldsymbol{w}_2)_h, & &\forall \boldsymbol{w}_2 \in \overline{Q}_3^{\fn,k},
\end{aligned} \label{3.25+b}\\
&\begin{aligned}
(\boldsymbol{w}_{\overline{S}}^{\mathfrak{n},3},\boldsymbol{w}_3)_h = \frac{1}{2} (\boldsymbol{w}_{S}^{\mathfrak{n},1},\boldsymbol{w}_3)_h + \frac{1}{2} (\boldsymbol{w}_{\overline{S}}^{\mathfrak{n},2},\boldsymbol{w}_3)_h + \frac{\tau}{2} (\mathcal{N}(w_{\overline{S}}^{\mathfrak{n},2}),\boldsymbol{w}_3)_h,& &\forall \boldsymbol{w}_3 \in \overline{Q}_3^{\fn,k},
\end{aligned} \label{3.25+c}\\
&\begin{aligned}
\boldsymbol{w}_{h}^{\mathfrak{n}+1} = e^{\frac{\tau}{2} \epsilon^2\mathcal{L}_h} \boldsymbol{w}_S^{\mathfrak{n},3}.
\end{aligned} \label{3.25+d}
\end{align}
\end{subequations}
Here, the operator $e^{\frac{\tau}{2} \epsilon^2 \mathcal{L}_h} $ is defined in $\eqref{finiteform}$, and the low-rank finite element approximations
\[
\begin{aligned}
\boldsymbol{w}_{S}^{\mathfrak{n},1} &= \Phi(x)^{\top} U_{\mathfrak{n},1} S_{\mathfrak{n},1} V_{\mathfrak{n},1}^{\top} \Phi(y) \in Q_{1}^{\fn,k}, \\
\boldsymbol{w}_{\overline{S}}^{\mathfrak{n},2} &= \Phi(x)^{\top} \overline{U}_{\mathfrak{n},3} \overline{S}_{\mathfrak{n},2} \overline{V}_{\mathfrak{n},3}^{\top} \Phi(y) \in \overline{Q}_3^{\fn,k}, \\
\boldsymbol{w}_{\overline{S}}^{\mathfrak{n},3} &= \Phi(x)^{\top} \overline{U}_{\mathfrak{n},3} \overline{S}_{\mathfrak{n},3} \overline{V}_{\mathfrak{n},3}^{\top} \Phi(y) \in \overline{Q}_3^{\fn,k}, \\
\boldsymbol{w}_{S}^{\mathfrak{n},3} &= \Phi(x)^{\top} U_{\mathfrak{n},3} S_{\mathfrak{n},3} V_{\mathfrak{n},3}^{\top} \Phi(y) \in Q_3^{\fn,k},
\end{aligned}
\]
where $\boldsymbol{w}_{S}^{\mathfrak{n},3}$ is the truncated low-rank finite element approximations of $\boldsymbol{w}_{\overline{S}}^{\mathfrak{n},3}$.
\end{lem}

\begin{proof}
It can be verified that the $\eqref{3.25+a}$ and $\eqref{3.25+d}$ are equivalent to the two linear steps in \Cref{lrnonlinear2}, respectively. Below, we demonstrate the equivalence for \eqref{3.25+a}; the proof for \eqref{3.25+d} follows a similar argument. By \Cref{lralg}, it holds
\begin{equation}\label{3.26+}
\begin{aligned}
\boldsymbol{w}_{S}^{n,1} & = \Phi(x)^{\top} U_{\mathfrak{n},1} S_{\mathfrak{n},1} V_{\mathfrak{n},1}^{\top} \Phi(y) = \Phi(x)^{\top} U_{\mathfrak{n},1} R S_{\mathfrak{n}} P^{\top} V_{\fn,1}^{\top} \Phi(y)\\
&= \Phi(x)^{\top} e^{\frac{\tau}{2} \epsilon^2 L_x} U _{\fn} S_{\mathfrak{n}} V_{\fn}^{\top}e^{\frac{\tau}{2} \epsilon^2 L_y^{\top}} \Phi(y) = e^{\frac{\tau}{2} \epsilon^2 \mathcal{L}_h} \boldsymbol{w}_h^{\fn}.
\end{aligned}
\end{equation}
In this following, we focus on $\eqref{3.25+b}$ and $\eqref{3.25+c}$.
Let \[\boldsymbol{w}_{\overline{S}}^{\fn,1} = \Phi(x)^{\top} \overline{U}_{\fn,3} \overline{S}_{\fn,1} \overline{V}^{\top}_{\fn,3} \Psi(y) \in \overline{Q}_3^{\fn,k},\]
be the projection of $\boldsymbol{w}_{S}^{\mathfrak{n},1} \in Q_{1}^{\fn,k}$ onto $\overline{Q}_3^{\mathfrak{n},k}$, namely,
\[(\boldsymbol{w}_{\overline{S}}^{\mathfrak{n},1}, \boldsymbol{w}_h)_h = (\boldsymbol{w}_{S}^{\mathfrak{n},1}, \boldsymbol{w}_h)_h, \quad \forall \boldsymbol{w}_h \in \overline{Q}_3^{\mathfrak{n},k},\]
which, by \Cref{lem2.2}, is equivalent to the matrix equation
\begin{equation}\label{3.26++}
\left(M_x \overline{U}_{\mathfrak{n},3} \overline{S}_{\mathfrak{n},1} \overline{V}_{\mathfrak{n},3}^{\top}M_y, \overline{U}_{\mathfrak{n},3} S \overline{V}_{\mathfrak{n},3}^{\top}\right)_{\rm F} = \left(M_x \overline{U}_{\mathfrak{n},3} S_{\mathfrak{n},1} \overline{V}_{\mathfrak{n},3}^{\top}M_y, \overline{U}_{\mathfrak{n},3} S \overline{V}_{\mathfrak{n},3}^{\top}\right)_{\rm F}, \quad  \forall S \in \mathbb{R}^{\bar{r} \times \bar{r}}.
\end{equation}
From \eqref{Q1nk}, it also holds
$$
\boldsymbol{w}_{\overline{S}}^{\mathfrak{n},1} = \boldsymbol{w}_{S}^{\mathfrak{n},1} \in \overline{Q}_3^{\mathfrak{n},k},
$$
By \Cref{lem2.2} and \Cref{intoin}, the low-rank finite element formulation \eqref{3.25+b} is equivalent to 
\[\left( M_x \overline{U}_{\mathfrak{n},3} (\overline{S}_{\mathfrak{n},2} - S_{\mathfrak{n},1} )\overline{V}_{\mathfrak{n},3}^{\top} M_y, \overline{U}_{\mathfrak{n},3} S \overline{V}_{\mathfrak{n},3}^{\top} \right)_{\rm F} = \tau \left(M_x \mathcal{N}(U_{\mathfrak{n},1} S_{\mathfrak{n},1} V_{\mathfrak{n},1}^{\top}) M_y, \overline{U}_{\mathfrak{n},3} S \overline{V}_{\mathfrak{n},3}^{\top}\right)_{\rm F},\quad  \forall S \in \mathbb{R}^{\bar{r} \times \bar{r}}.\]
Using \eqref{3.26++} and the properties of the Frobenius inner product gives
\begin{equation}\label{3.28+}
\left( \overline{S}_{\mathfrak{n},2},S \right)_{\rm F} = \left(\overline{S}_{\mathfrak{n},1},S \right)_{\rm F} + \tau \left(\overline{U}_{\mathfrak{n},3}^{\top} M_x \mathcal{N}(U_{\mathfrak{n},1} S_{\mathfrak{n},1} V_{\mathfrak{n},1}^{\top})M_y \overline{V}_{\mathfrak{n},3}, S\right)_{\rm F}, \quad  \forall S \in \mathbb{R}^{\bar{r} \times \bar{r}}.
\end{equation}
From $\eqref{3.21}$,  $\overline{U}_{\fn,3}\overline{U}_{\fn,3}^{\top} M_x$ is the orthogonal projection onto the range of $\overline{U}_{\fn,3}$, which by definition equals the range of $\mathbb{U}$. In particular, the columns of $U_{\mathfrak{n},1}$ lie in range of $\bar U_{\fn,3}$, and hence $\overline{U}_{\fn,3}\overline{U}_{\fn,3}^{\top} M_x U_{\mathfrak{n},1}=\overline{U}_{\fn,3}M = U_{\mathfrak{n},1} $. Similarly, $\overline{V}_{\fn,3}\overline{V}_{\fn,3}^{\top} M_y V_{\mathfrak{n},1} =\overline{V}_{\fn,3}N = V_{\mathfrak{n},1}$. Here, $M$ and $N$ were given in \Cref{lrnonlinear2}.
Therefore, it holds
\begin{equation}\label{3.29+}
U_{\mathfrak{n},1} S_{\mathfrak{n},1} V_{\mathfrak{n},1}^{\top} = \overline{U}_{\fn,3} M S_{\mathfrak{n},1} N^\top \overline{V}^{\top}_{\fn,3} = \overline{U}_{\fn,3} \overline{S}_{\fn,1} \overline{V}^{\top}_{\fn,3}.
\end{equation}
Then, $\eqref{3.28+}$ can be written as
\[\left( \overline{S}_{\mathfrak{n},2}, S\right)_{\rm F} = \left(\overline{S}_{\mathfrak{n},1},S \right)_{\rm F} + \tau \left(\overline{U}_{\mathfrak{n},3}^{\top} M_x \mathcal{N}(\overline{U}_{\mathfrak{n},3} \overline{S}_{\mathfrak{n},1} \overline{V}_{\mathfrak{n},3}^{\top})M_y \overline{V}_{\mathfrak{n},3}, S\right)_{\rm F}, \quad  \forall S \in \mathbb{R}^{\bar{r} \times \bar{r}},\]
which gives \eqref{3.22a} by the arbitrariness of $S$. Similarly, $\eqref{3.25+c}$ is equivalent to \eqref{3.22b}.
\end{proof}

\begin{rem}
In \Cref{lrfem}, we presented the finite element formulations for the {$S$-step} but omitted those for the {$K$-step} and {$L$-step}, as the finite element solution for the {$S$-step} is equivalent to the low-rank finite element solution.
\end{rem}

\begin{thm}
Given $\boldsymbol{w}_h^\mathfrak{n} 
\in Q_{0}^{\fn,k}$, the DLR-MLFEM solution $ \boldsymbol{w}_h^{\mathfrak{n}+1} \in Q_{0}^{\fn+1,k}$ obtained from \Cref{lrnonlinear2} for the conservative AC equation \eqref{1.2} conserves the mass up to the truncation tolerance 
\[|(\boldsymbol{w}_h^{\mathfrak{n}+1} - \boldsymbol{w}_h^{\mathfrak{n}},1)_h| = (\mathcal{W}_{\fn+1},\mathbf{I})_{\rm M} -(\mathcal{W}_{\fn},\mathbf{I})_{\rm M}| = |(U_{\fn+1} S_{\fn+1} V^{\top}_{\fn+1},\mathbf{I})_{\rm M} - (U_{\fn} S_{\fn} V^{\top}_{\fn},\mathbf{I})_{\rm M} | \leq C \eta \;.\]
\end{thm}

\begin{proof}
The matrices $U_{\mathfrak{n},1}, S_{\mathfrak{n},1}, V_{\mathfrak{n},1}$ are obtained from \Cref{lralg} with initial values $\{ U_{\mathfrak{n}},S_{\mathfrak{n}},V_{\mathfrak{n}} \}$ and the time step $\frac{\tau}{2}$. Then, by Lemma \ref{lem3.1}, it follows
\begin{equation}\label{3.24}
(U_{\mathfrak{n},1} S_{\mathfrak{n},1} V_{\mathfrak{n},1}^{\top}, \mathbf{I})_{\rm M} =(U_{\mathfrak{n}}S_{\mathfrak{n}}V_{\mathfrak{n}},\mathbf{I})_{\rm M}.
\end{equation}  
Recall that the all-ones matrix \(\mathbf{I} = q_m q_n^\top \in \mathbb{R}^{m \times n}\), where \(q_m \in \mathbb{R}^m\) and \(q_n \in \mathbb{R}^n\) are all-ones vectors. 
Similar to $\eqref{3.29+}$, $ \overline{U}_{\fn,3}\overline{U}_{\fn,3}^{\top} M_x q_m = q_m$, and $\overline{V}_{\fn,3}\overline{V}_{\fn,3}^{\top} M_y q_n = q_n$. Then, for any matrix $W \in \mathbb{R}^{m \times n}$, 
\begin{equation}\label{3.31+}
\begin{aligned}
 \left(\overline{U}_{\fn,3} \overline{U}_{\fn,3}^{\top} M_x \mathcal{N}(W) M_y \overline{V}_{\fn,3} \overline{V}_{\fn,3}^{\top} , \mathbf{I} \right)_{\rm M} 
& = \left( M_x \overline{U}_{\fn,3}\overline{U}_{\fn,3}^{\top} M_x \mathcal{N}(W) M_y \overline{V}_{\fn,3} \overline{V}_{\fn,3}^{\top} M_y, \mathbf{I} \right)_{\rm F} \\
& = \left( M_x \mathcal{N}(W) M_y , \mathbf{I} \right)_{\rm F} = \left(\mathcal{N}(W), \mathbf{I} \right)_{\rm M}.
\end{aligned}
\end{equation}
Therefore, $\overline{U}_{\mathfrak{n},3}, \overline{S}_{\mathfrak{n},3}, \overline{V}_{\fn,3}$ obtained from \Cref{lrnonlinear2} with initial values $U_{\mathfrak{n},1}, S_{\mathfrak{n},1},V_{\mathfrak{n},1}$ satisfy
\begin{equation}
\begin{aligned}
&(\overline{U}_{\fn,3}\overline{S}_{\mathfrak{n},3}\overline{V}_{\fn,3}^{\top},\mathbf{I})_{\rm M}  = \frac{1}{2}\left( \overline{U}_{\fn,3} \left[ \overline{S}_{\mathfrak{n},1} + \overline{S}_{\mathfrak{n},2} + \tau \overline{U}_{\fn,3}^{\top} M_x \mathcal{N}(\overline{U}_{\fn,3}\; \overline{S}_{\mathfrak{n},2} \overline{V}_{\fn,3}^{\top})M_y \overline{V}_{\fn,3}\right]\overline{V}_{\fn,3}^{\top}, \mathbf{I} \right)_{\rm M}\\
& \quad = \left( \overline{U}_{\fn,3} \;\overline{S}_{\mathfrak{n},1} \overline{V}_{\fn,3}^{\top}, \mathbf{I} \right)_{\rm M} + \frac{\tau}{2} \left(\mathcal{N}(\overline{U}_{\fn,3} \;\overline{S}_{\mathfrak{n},1} \overline{V}_{\fn,3}^{\top}), \mathbf{I}\right)_{\rm M} + \frac{\tau}{2} \left( \mathcal{N}(\overline{U}_{\fn,3}\; \overline{S}_{\mathfrak{n},2} \overline{V}_{\fn,3}^{\top}), \mathbf{I} \right)_{\rm M}\\
& \quad = \left( \overline{U}_{\fn,3} \;\overline{U}_{\fn,3}^{\top} U_{\mathfrak{n},1}S_{\mathfrak{n},1}V_{\mathfrak{n},1}^{\top} \overline{V}_{\fn,3} \; \overline{V}_{\fn,3}^{\top}, \mathbf{I}  \right)_{\rm M} = \left( U_{\mathfrak{n},1} S_{\mathfrak{n},1} V^{\top}_{\mathfrak{n},1}, \mathbf{I} \right)_{\rm M},
\end{aligned}
\end{equation}
where we have used \eqref{2.27}.
The Truncation step in \Cref{lrnonlinear2} implies
\begin{equation}\label{3.27}
\begin{aligned}
\left| (U_{\fn,3} S_{\fn,3} V_{\fn,3}^{\top}, \mathbf{I})_{\rm M} - (\overline{U}_{\fn,3} \overline{S}_{\fn,3} \overline{V}_{\fn,3}^{\top}, \mathbf{I})_{\rm M}\right| &\leq \left\|U_{\fn,3} S_{\fn,3} V_{\fn,3}^{\top} - \overline{U}_{\fn,3} \overline{S}_{\fn,3} \overline{V}_{\fn,3}^{\top}\right\|_{\rm M} \|\mathbf{I}\|_{\rm M} \leq C \eta .
\end{aligned}
\end{equation}
Similar to $\eqref{3.24}$, the second linear step implies
\[(U_{\mathfrak{n}+1} S_{\mathfrak{n}+1} V_{\mathfrak{n}+1}, \mathbf{I})_{\rm M} = (U_{\fn,3} S_{\fn,3} V_{\fn,3}^{\top}, \mathbf{I})_{\rm M},\]
Therefore, it holds
\[\left|(U_{\mathfrak{n}+1} S_{\mathfrak{n}+1} V_{\mathfrak{n}+1}^{\top}, \mathbf{I})_{\rm M} - (U_{\mathfrak{n}} S_{\mathfrak{n}} V_\mathfrak{n}^{\top}, \mathbf{I})_{\rm M}  \right| \leq C \eta.\]
which implies the conclusion.
\end{proof}

To give the discrete energy law for the low-rank solution, we introduce the following assumptions.
\begin{assum}\label{assumRl}
(i) $\mathcal{N}$ is Lipschitz-continuous and bounded: for all $X, Y \in \mathbb{R}^{m \times n}$,
\[ \|\mathcal{N}(X) - \mathcal{N}(Y)\|_{M} \leq C_{N} \|X-Y\|_{\rm M}.\]
(ii) For coefficient matrix $\mathcal{W}(t) \in \mathbb{R}^{m \times n}$ to the low rank solution $\boldsymbol{w}_h \in \overline{Q}_{3}^{n,k}$, 
it holds projection error
\[\left\| \overline{U}_{\fn,3}\; \overline{U}_{\fn,3}^{\top} M_x \mathcal{N}(\mathcal{W}) M_y \overline{V}_{\fn,3} \overline{V}_{\fn,3}^{\top} - \mathcal{N}(\mathcal{W}) \right\|_{\rm M} \leq  \mu,\]
where the bases are given in \Cref{lrnonlinear2}.
\end{assum}

Then we are ready to state the following modified energy law.
\begin{thm}\label{energythm}
For both the classical AC equation \eqref{1.1} and the conservative AC equation \eqref{1.2}, if the time step $\tau$ satisfies the conditions stated in $\Cref{lem2.5}$, the DLR-MLFEM solution $\boldsymbol{w}_h^{\mathfrak{n}+1} = \Phi(x)^{\top} U_{\mathfrak{n}+1} S_{\mathfrak{n}+1} V^{\top}_{\mathfrak{n}+1}\Psi(y)$, obtained via \Cref{lrnonlinear2} or \eqref{finiteform2}, satisfies the following modified energy law 
\begin{equation}\label{lrenglaw}
    \widetilde{E}^{\mathfrak{n}+1}\leq \widetilde{E}^{\mathfrak{n}} - \beta^2 + \left(\frac{\eta}{\tau} +\frac{(1+C_N \tau) \mu}{2}\right)\beta,
\end{equation}
where $\beta = \left\|\boldsymbol{w}_h^{\mathfrak{n}+1,1} - \boldsymbol{w}_h^{\mathfrak{n},1} \right\|_h$, and the modified energy
\begin{equation}\label{3.28}
\widetilde{E}^{\fn}=\widetilde{E}\left(\boldsymbol{w}_h^{\mathfrak{n},1} \right) = \frac{1}{2\tau} \left[ (e^{-\tau \epsilon^2 \mathcal{L}_h}-1)  \boldsymbol{w}_h^{\mathfrak{n},1},\boldsymbol{w}_h^{\mathfrak{n},1})_h \right] + (G(\boldsymbol{w}_h^{\mathfrak{n},1}),1)_h,
\end{equation}
with $G$ being defined in \eqref{2.28}.
\end{thm}
\begin{proof}
From $\eqref{3.22}$, it holds
\begin{equation}\label{3.29}
\left\lbrace
\begin{aligned}
\overline{U}_{\fn,3} \overline{S}_{\mathfrak{n},2} \overline{V}_{\fn,3}^{\top} &= \overline{U}_{\fn,3} \overline{S}_{\mathfrak{n},1} \overline{V}_{\fn,3}^{\top} + \tau \overline{U}_{\fn,3} \overline{U}_{\fn,3}^{\top} M_x \mathcal{N} (\overline{U}_{\fn,3} \overline{S}_{\mathfrak{n},1} \overline{V}_{\fn,3}^{\top}) M_y \overline{V}_{\fn,3}\overline{V}_{\fn,3}^{\top}. \\
\overline{U}_{\fn,3} \overline{S}_{\mathfrak{n},3} \overline{V}_{\fn,3}^{\top} &= \frac{1}{2}\overline{U}_{\fn,3} \overline{S}_{\mathfrak{n},1} \overline{V}_{\fn,3}^{\top} +\frac{1}{2} \overline{U}_{\fn,3} \overline{S}_{\mathfrak{n},2} \overline{V}_{\fn,3}^{\top} + \frac{\tau}{2} \overline{U}_{\fn,3} \overline{U}_{\fn,3}^{\top} M_x \mathcal{N} (\overline{U}_{\fn,3} \overline{S}_{\mathfrak{n},2} \overline{V}_{\fn,3}^{\top})M_y \overline{V}_{\fn,3} \overline{V}_{\fn,3}^{\top}.
\end{aligned}
\right.
\end{equation}
Based on $\eqref{3.29+}$, we define
\begin{equation}\label{3.30}
\overline{\mathcal{W}}_{\fn,3} = \overline{U}_{\fn,3} \overline{S}_{\mathfrak{n},3} \overline{V}_{\fn,3}^{\top}.   
\end{equation}
and 
\begin{equation}\label{3.31}
\begin{aligned}
\overline{\mathcal{W}}_{\mathfrak{n},2} = \overline{U}_{\fn,3} \overline{S}_{\mathfrak{n},2} \overline{V}_{\fn,3}^{\top}  &= \overline{U}_{\fn,3} \overline{S}_{\mathfrak{n},1} \overline{V}_{\fn,3}^{\top} + \tau \overline{U}_{\fn,3}\overline{U}_{\fn,3}^{\top} M_x \mathcal{N} (\overline{U}_{\fn,3} \overline{S}_{\mathfrak{n},1} \overline{V}_{\fn,3}^{\top}) M_y \overline{V}_{\fn,3}\overline{V}_{\fn,3}^{\top} \\
& = U_{\mathfrak{n},1} S_{\mathfrak{n},1} V_{\mathfrak{n},1}^{\top} + \tau \overline{U}_{\fn,3}\overline{U}_{\fn,3}^{\top} M_x \mathcal{N} (U_{\mathfrak{n},1} S_{\mathfrak{n},1} V_{\mathfrak{n},1}^{\top}) M_y \overline{V}_{\fn,3}\overline{V}_{\fn,3}^{\top} \\
& = \mathcal{W}_{\mathfrak{n},1} + \tau \overline{U}_{\fn,3}\overline{U}_{\fn,3}^{\top} M_x \mathcal{N}(\mathcal{W}_{\mathfrak{n},1})M_y \overline{V}_{\fn,3} \overline{V}_{\fn,3}^{\top}.
\end{aligned}
\end{equation}
Thus, it can be concluded from $\eqref{3.29}$, $\eqref{3.30}$ and $\eqref{3.31}$:
\begin{equation}\label{df31}
\overline{\mathcal{W}}_{\mathfrak{n},3}  - \mathcal{W}_{\mathfrak{n},1} = \frac{\tau}{2} \overline{U}_{\fp} \mathcal{N}(\mathcal{W}_{\mathfrak{n},1}) \overline{V}_{\fp} + \frac{\tau}{2} \overline{U}_{\fp} \mathcal{N}(\overline{\mathcal{W}}_{\mathfrak{n},2}) \overline{V}_{\fp} ,    
\end{equation}
where $\overline{U}_{\fp} = \overline{U}_{\fn,3} \overline{U}_{\fn,3}^{\top} M_x$ and $\overline{V}_{\fp} = \overline{V}_{\fn,3} \overline{V}_{\fn,3}^{\top} M_y$ denote the orthogonal projector on to the range of $\overline{U}_{\fn,3}$ and $\overline{V}_{\fn,3}$, respectively. 
\eqref{df31} can be reformulated as
\begin{equation}\label{3.32}
\begin{aligned}
&\frac{1}{\tau} (\overline{\mathcal{W}}_{\mathfrak{n},3} -  \mathcal{W}_{\mathfrak{n}+1,1} )+ \frac{1}{\tau}( \mathcal{W}_{\mathfrak{n}+1,1} - \mathcal{W}_{\mathfrak{n},1}) = \frac{1}{2} \overline{U}_{\fp} \mathcal{N}(\mathcal{W}_{\mathfrak{n},1}) \overline{V}_{\fp} + \frac{1}{2} \overline{U}_{\fp}\mathcal{N}(\overline{\mathcal{W}}_{\mathfrak{n},2}) \overline{V}_{\fp}.
\end{aligned}
\end{equation}

The first term of the left side of $\eqref{3.32}$, taking ${\rm M}$-weighted inner product with $\mathcal{W}_{\mathfrak{n}+1,1}- \mathcal{W}_{\mathfrak{n},1}$, gives
\begin{equation}\label{3.33}
\begin{aligned}
&\quad \frac{1}{\tau}(\overline{\mathcal{W}}_{\mathfrak{n},3} -  \mathcal{W}_{\mathfrak{n}+1,1},  \mathcal{W}_{\mathfrak{n}+1,1}- \mathcal{W}_{\mathfrak{n},1})_{\rm M} \\
&=\frac{1}{\tau}(\overline{\mathcal{W}}_{\mathfrak{n},3} -\mathcal{W}_{\fn,3} ,  \mathcal{W}_{\mathfrak{n}+1,1} -\mathcal{W}_{\mathfrak{n},1})_{\rm M} +\frac{1}{\tau}(\mathcal{W}_{\fn,3} - \mathcal{W}_{\mathfrak{n}+1,1},  \mathcal{W}_{\mathfrak{n}+1,1} -\mathcal{W}_{\mathfrak{n},1})_{\rm M}.
\end{aligned}
\end{equation}
By \eqref{3.23}, the first term on the right hand side of $\eqref{3.33}$ implies
\[
\begin{aligned}
\frac{1}{\tau}(\overline{\mathcal{W}}_{\mathfrak{n},3} -\mathcal{W}_{\fn,3} , \mathcal{W}_{\mathfrak{n}+1,1} -\mathcal{W}_{\mathfrak{n},1})_{\rm M} &\leq \frac{1}{\tau}  \left\| \overline{\mathcal{W}}_{\mathfrak{n},3} -\mathcal{W}_{\fn,3}\right\|_{\rm M} \left\|\mathcal{W}_{\mathfrak{n}+1,1} -\mathcal{W}_{\mathfrak{n},1}\right\|_{\rm M}\leq \frac{\eta}{\tau} \left\|\mathcal{W}_{\mathfrak{n}+1,1} -\mathcal{W}_{\mathfrak{n},1}\right\|_{\rm M}. 
\end{aligned}\]
For linear step in \Cref{lrnonlinear2}, 
\[\mathcal{W}_{\fn,3}= e^{-\tau \epsilon^2 \mathcal{L}} \mathcal{W}_{\fn+1,1}.\]
Then the second term on the right hand side of $\eqref{3.33}$ holds
\begin{equation}\label{3.34}
\begin{aligned}
& \frac{1}{\tau}(\mathcal{W}_{\fn,3} - \mathcal{W}_{\mathfrak{n}+1,1},\mathcal{W}_{\mathfrak{n}+1,1} -\mathcal{W}_{\mathfrak{n},1})_{\rm M}=\frac{1}{\tau}(e^{-\tau \epsilon^2 \mathcal{L}} \mathcal{W}_{\mathfrak{n}+1,1} - \mathcal{W}_{\mathfrak{n}+1,1},\mathcal{W}_{\mathfrak{n}+1,1} -\mathcal{W}_{\mathfrak{n},1})_{\rm M} \\
&\quad=\frac{1}{2\tau} \left[ (e^{-\tau  \epsilon^2\mathcal{L}}  \mathcal{W}_{\mathfrak{n}+1,1} - \mathcal{W}_{\mathfrak{n}+1,1},\mathcal{W}_{\mathfrak{n}+1,1})_{\rm M} - (e^{-\tau  \epsilon^2\mathcal{L}}  \mathcal{W}_{\mathfrak{n},1}  - \mathcal{W}_{\mathfrak{n},1}, \mathcal{W}_{\mathfrak{n},1})_{\rm M} \right. \\
&\qquad \qquad +\left.(e^{-\tau  \epsilon^2\mathcal{L}}  (\mathcal{W}_{\mathfrak{n}+1,1}-\mathcal{W}_{\mathfrak{n},1}), \mathcal{W}_{\mathfrak{n}+1,1}-\mathcal{W}_{\mathfrak{n},1})_{\rm M} - \left\| \mathcal{W}_{\mathfrak{n}+1,1}-\mathcal{W}_{\mathfrak{n},1}\right\|^2_{\rm M}  \right],
\end{aligned}
\end{equation}
where we use the identity
\[(A-B,2A)_M =\|A\|_M^2 - \|B\|_M^2 + \|A-B\|_M^2.\]

The right side of $\eqref{3.32}$, taking ${\rm M}$-weighted inner product with $\mathcal{W}_{\mathfrak{n}+1,1}- \mathcal{W}_{\mathfrak{n},1}$, gives
\begin{align*}
&\quad \frac{1}{2} \left( \overline{U}_{\fp}  \mathcal{N}(\mathcal{W}_{\mathfrak{n},1})\overline{V}_{\fp} + \overline{U}_{\fp}\mathcal{N}(\overline{\mathcal{W}}_{\mathfrak{n},2}) \overline{V}_{\fp} , \mathcal{W}_{\fn+1,1}-\mathcal{W}_{\fn,1}\right)_{\rm M} \\
& =\frac{1}{2} \left( \overline{U}_{\fp}  \mathcal{N}(\mathcal{W}_{\mathfrak{n},1})\overline{V}_{\fp} + \overline{U}_{\fp}\mathcal{N}(\overline{\mathcal{W}}_{\mathfrak{n},2}) \overline{V}_{\fp}-\left(\mathcal{N}(\mathcal{W}_{\mathfrak{n},1}) +  \mathcal{N}(\overline{U}_{\fn,3} \overline{S}_{\mathfrak{n},2} \overline{V}_{\fn,3}^{\top})\right), \mathcal{W}_{\fn+1,1}-\mathcal{W}_{\fn,1}\right)_{\rm M} \\
& \qquad \qquad \qquad+ \frac{1}{2}\left( \left(\mathcal{N}(\mathcal{W}_{\mathfrak{n},1}) +  \mathcal{N}(\overline{U}_{\fn,3} \overline{S}_{\mathfrak{n},2} \overline{V}_{\fn,3}^{\top})\right), \mathcal{W}_{\fn+1,1}-\mathcal{W}_{\fn,1}\right)_{\rm M} \\
& \leq \frac{1}{2} (\mathcal{N}(\mathcal{W}_{\mathfrak{n},1}) +  \mathcal{N}(\overline{U}_{\fn,3} \overline{S}_{\mathfrak{n},2} \overline{V}_{\fn,3}^{\top}),  \mathcal{W}_{\mathfrak{n}+1,1} -\mathcal{W}_{\mathfrak{n},1})_{\rm M} +\frac{\mu}{2}\left\|\mathcal{W}_{\fn+1,1}-\mathcal{W}_{\fn,1}\right\|_{\rm M} \\
& \leq \frac{1}{2}( \mathcal{N}(\mathcal{W}_{\mathfrak{n},1})+\mathcal{N}(\mathcal{W}_{\mathfrak{n},1} + \tau \mathcal{N}(\mathcal{W}_{\mathfrak{n},1})),\mathcal{W}_{\mathfrak{n}+1,1} - \mathcal{W}_{\mathfrak{n},1})_{\rm M} + \frac{\mu}{2}\left\|\mathcal{W}_{\fn+1,1}-\mathcal{W}_{\fn,1}\right\|_{\rm M}\\
&\qquad \qquad \qquad + \frac{1}{2}(\mathcal{N}(\overline{U}_{\fn,3} \overline{S}_{\mathfrak{n},2} \overline{V}_{\fn,3}^{\top}) - \mathcal{N}(\mathcal{W}_{\mathfrak{n},1} + \tau \mathcal{N}(\mathcal{W}_{\mathfrak{n},1})), \mathcal{W}_{\mathfrak{n}+1,1} - \mathcal{W}_{\fn,1})_{\rm M}. \\
&\leq \frac{1}{2}( \mathcal{N}(\mathcal{W}_{\mathfrak{n},1})+\mathcal{N}(\mathcal{W}_{\mathfrak{n},1} + \tau \mathcal{N}(\mathcal{W}_{\mathfrak{n},1})),\mathcal{W}_{\mathfrak{n}+1,1} - \mathcal{W}_{\mathfrak{n},1})_{\rm M}  + \frac{\mu}{2}\left\|\mathcal{W}_{\fn+1,1}-\mathcal{W}_{\fn,1}\right\|_{\rm M}\\
&\qquad \qquad \qquad + \frac{\tau C_N}{2} \left\| \overline{U}_{\fn,3} \overline{U}_{\fn,3}^{T} M_x \mathcal{N}(\mathcal{W}_{\mathfrak{n},1}) M_y \overline{V}_{\fn,3} \overline{V}_{\fn,3}^{T}-  \mathcal{N}(\mathcal{W}_{\mathfrak{n},1})\right\|_{\rm M} \left\| \mathcal{W}_{\mathfrak{n}+1,1}-\mathcal{W}_{\mathfrak{n},1}\right\|_{\rm M}  \\
&\leq -( g(\mathcal{W}_{\mathfrak{n},1}),\mathcal{W}_{\mathfrak{n}+1,1} - \mathcal{W}_{\mathfrak{n},1})_{\rm M}  + \frac{\mu \tau C_N}{2}\left\| \mathcal{W}_{\mathfrak{n}+1,1}-\mathcal{W}_{\mathfrak{n},1}\right\|_{\rm M}+ \frac{\mu}{2}\left\|\mathcal{W}_{\fn+1,1}-\mathcal{W}_{\fn,1}\right\|_{\rm M} \\
&\leq -( g(\mathcal{W}_{\mathfrak{n},1}),\mathcal{W}_{\mathfrak{n}+1,1} - \mathcal{W}_{\mathfrak{n},1})_{\rm M}  + \frac{(1+C_N\tau)\mu}{2} \left\| \mathcal{W}_{\mathfrak{n}+1,1}-\mathcal{W}_{\mathfrak{n},1}\right\|_{\rm M},
\end{align*}
where $g$ is defined in \eqref{gweqn}, and \Cref{assumRl} has been applied in the first, third, and fourth inequalities.

Combining the results above, we obtain
\begin{equation}\label{3.36}
\begin{aligned}
& \widetilde{E}(\boldsymbol{w}^{\mathfrak{n}+1,1}_h) -\widetilde{E}(\boldsymbol{w}^{\mathfrak{n},1}_h) \\
&=\frac{1}{2\tau} \left[ (e^{-\tau \mathcal{L}_h}-1) \boldsymbol{w}_h^{\fn+1,1}, \boldsymbol{w}_h^{\mathfrak{n}+1,1})_h - (e^{-\tau \mathcal{L}_h}-1) \boldsymbol{w}_h^{\mathfrak{n},1}, \boldsymbol{w}_h^{\fn,1})_h  +(G(\boldsymbol{w}_{h}^{\mathfrak{n}+1,1}),1)_h - G(\boldsymbol{w}_{h}^{\mathfrak{n},1}),1)_h\right] \\
&\leq -\left(\frac{1}{2\tau} - \frac{g^{\prime}(\xi)}{2}\right) \left\|\boldsymbol{w}_h^{\mathfrak{n}+1,1} - \boldsymbol{w}_h^{\mathfrak{n},1} \right\|^2_h + \frac{\eta}{\tau} \left\|\boldsymbol{w}_h^{\mathfrak{n}+1,1}- \boldsymbol{w}_h^{\mathfrak{n},1}\right\|_h + \frac{(1+C_N\tau)\mu}{2}\left\|\boldsymbol{w}_h^{\mathfrak{n}+1,1} - \boldsymbol{w}_h^{\mathfrak{n},1}\right\|_h,
\end{aligned}
\end{equation}
which yields the energy dissipation law \eqref{lrenglaw}.
\end{proof}

\begin{rem}
The energy dissipation law \eqref{lrenglaw} for the DLR-MLFEM solution has an additional term, $(\frac{\eta}{\tau}  + (1+C\tau) \mu) \beta-\beta^2$, which accounts for the truncation error and projection error in \Cref{assumRl}.
To ensure that \Cref{lrnonlinear2} achieves second-order accuracy in time, it is reasonable to assume
\[\eta = \mathcal{O}(\tau^2), \quad \text{and} \quad \beta \sim \|\boldsymbol{w}_h^{\fn+1}- \boldsymbol{w}_h^{\fn}\|_h = \mathcal{O}(\tau(\tau^2+h^{k+1})).\]
Under these assumptions, \eqref{lrenglaw} reduces to
\[\widetilde{E}^{\fn+1} \leq \widetilde{E}^{\fn}  + \mathcal{O}(\tau(\tau^2+h^{k+1})),\]
which implies that the last two terms in \eqref{lrenglaw} are higher-order error terms. Consequently, the energy stability property remains valid.
\end{rem}

\section{Numerical experiment}\label{Sec4}
In this section, we present some numerical experiments to demonstrate the advantages of the proposed method.
\subsection{Convergence test}
\begin{example}
We first verify the spatial and temporal convergence of the proposed schemes on classical Allen-Cahn equation with a smooth initial condition
\begin{equation}\label{4.1}
u(x,y,t=0) = sin(\pi x) sin(\pi y),  \quad (x,y) \in [0,1] \times [0,1].
\end{equation}
\end{example}
We set $\epsilon = 0.01$ and terminal time $T = 1.0.$ Since the exact solution is unavailable, the temporal error is computed by fixing the spatial mesh size $h=1/128$ and comparing the numerical solutions with a reference solution obtained using a very small time step ($\tau = 0.0001$). Similarly, the spatial error is computed by fixing $\tau = 0.0001$ and comparing the numerical solutions with a reference solution obtained using a refined mesh (with $h = 1/512$). The numerical errors are computed as:
\[
L^2 \text{ error} = \|w_{ref}-w_h\|_h.
\]
where $w_{ref}$ is the reference solution. 

\begin{figure}
\centering
\subfigure[FR-MLFEM]{
\begin{minipage}[t]{0.22\linewidth}
\centering
\includegraphics[width=0.95\textwidth]{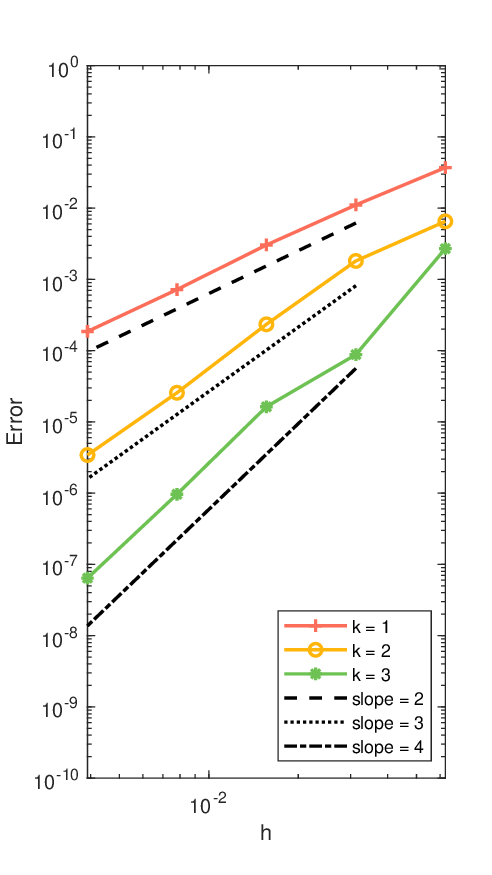}
\end{minipage}
}
\subfigure[DLR-MLFEM, $k = 1$]{
\begin{minipage}[t]{0.22\linewidth}
\centering
\includegraphics[width=0.95\textwidth]{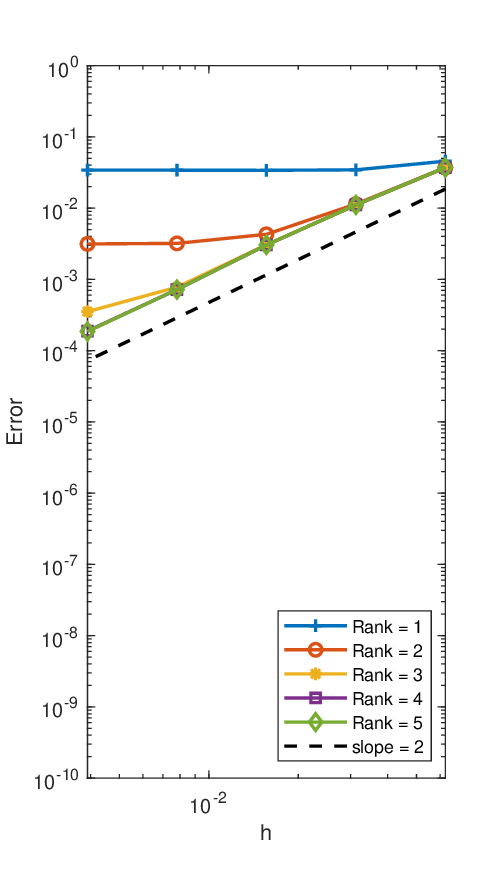}
\end{minipage}
}
\subfigure[DLR-MLFEM, $k = 2$]{
\begin{minipage}[t]{0.22\linewidth}
\centering
\includegraphics[width=0.95\textwidth]{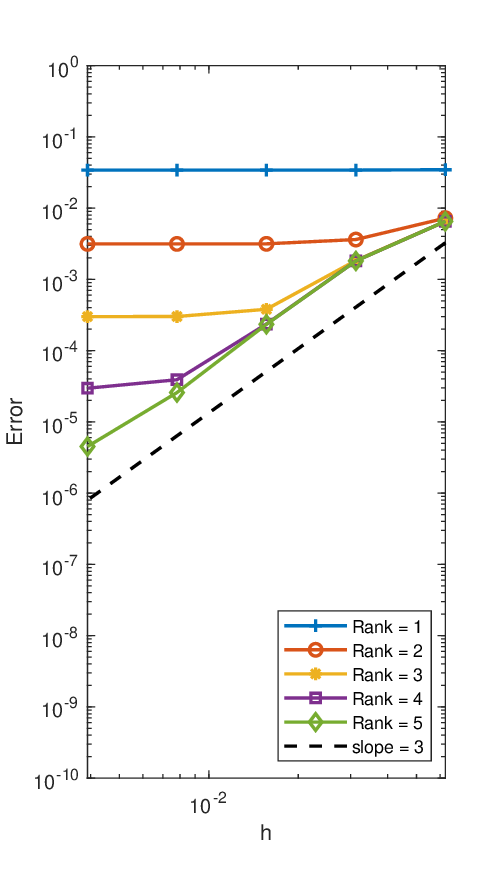}
\end{minipage}
}
\subfigure[DLR-MLFEM, $k = 3$]{
\begin{minipage}[t]{0.22\linewidth}
\centering
\includegraphics[width=0.95\textwidth]{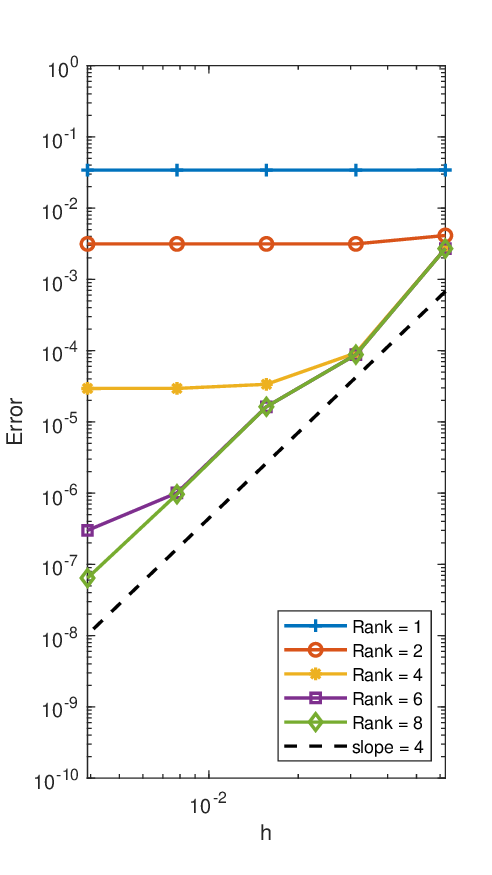}
\end{minipage}
}
\caption{Spatial accuracy tests of the DLR-MLFEM at $T=1$.}\label{SpatialOrder}
\end{figure}

The spatial errors at the final time $T$, computed using the FR-MLFEM and DLR-MLFEM for the classical AC equation, are presented in \Cref{SpatialOrder}. As expected, the FR-MLFEM achieves $(k+1)$-th order accuracy in space (\Cref{SpatialOrder}(a)). \Cref{SpatialOrder}(b)-(d) report the report the errors and spatial convergence rates for the DLR-MLFEM with different rank $r$. 
We observe that when $r$ is small, the truncation error dominates. However, for a moderately large rank ($r=8$ for this example), the same errors and convergence rates as the full-rank solution are achieved.

\begin{figure}
\centering
\subfigure[Lie-Trotter]{
\begin{minipage}[t]{0.23\linewidth}
\centering
\includegraphics[width=1\textwidth]{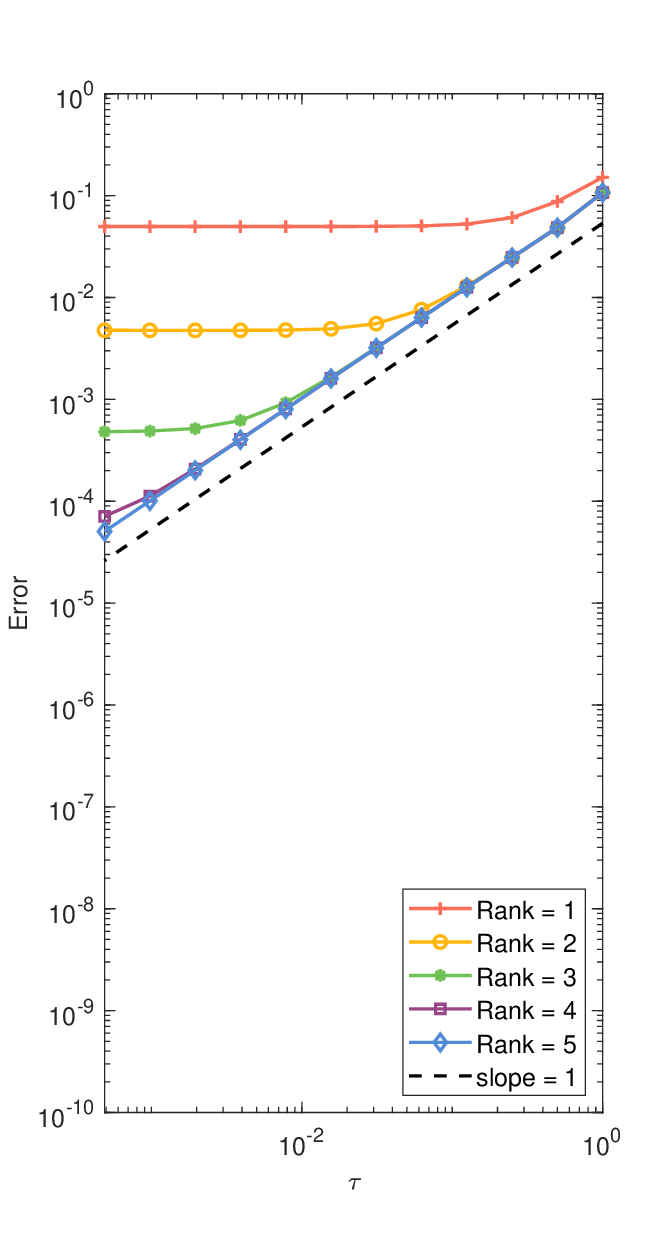}
\end{minipage}
}
\subfigure[DLR-MLFEM]{
\begin{minipage}[t]{0.23\linewidth}
\centering
\includegraphics[width=1\textwidth]{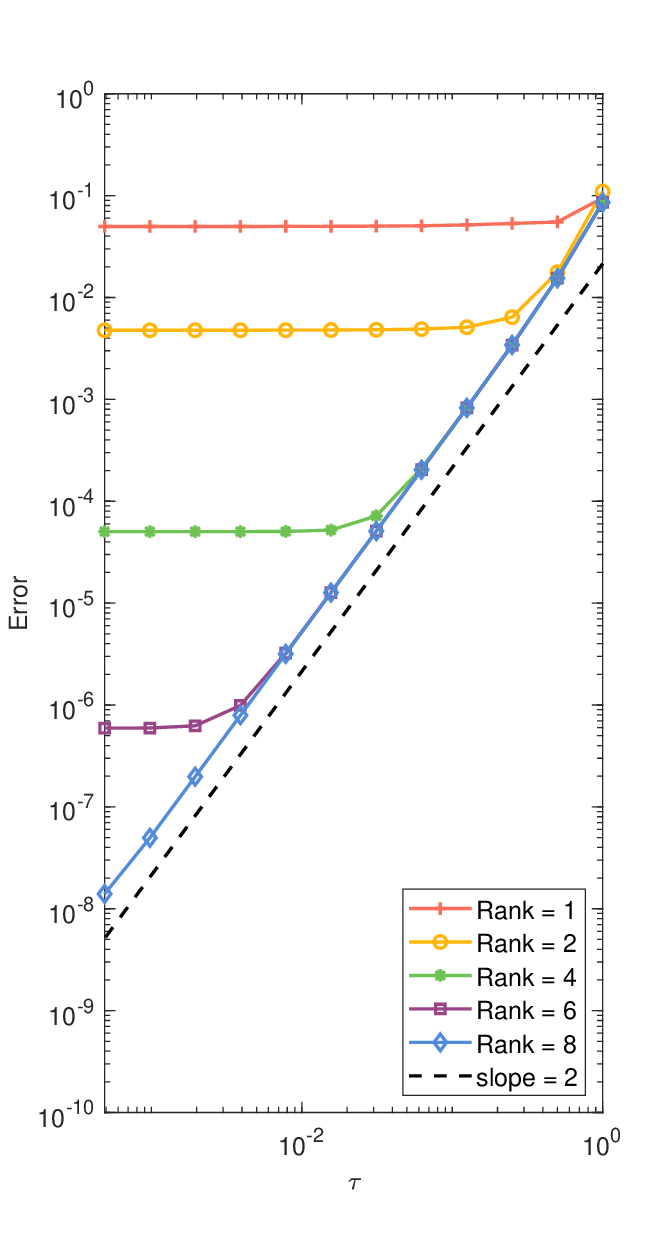}
\end{minipage}
}
\caption{Temporal accuracy tests of the DLR-MLFEM at $T=1$.}\label{TimeOrder}
\end{figure}

The temporal errors are presented in \Cref{TimeOrder}. 
 In  \Cref{TimeOrder}(a), we additionally provide the errors and convergence rates computed by the first order Lie-Trotter integrator (see \Cref{lr1st}).
The results demonstrate that the Lie-Trotter dynamical low-rank finite element solution achieves first-order accuracy in time, while the proposed \Cref{lrnonlinear2} attains second-order accuracy in time.

\subsection{Modified energy}
\begin{figure}
\centering
\subfigure[AC, FR-MLFEM]{
\begin{minipage}[t]{0.48\linewidth}
\centering
\includegraphics[height=6cm,width=6cm]{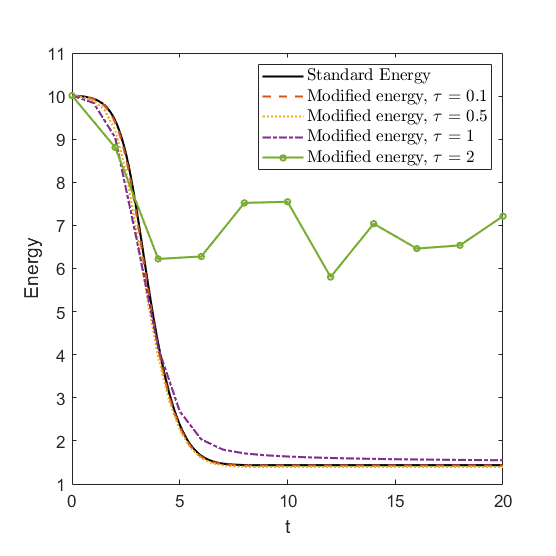}
\end{minipage}
}
\subfigure[AC, adaptive DLR-MLFEM.]{
\begin{minipage}[t]{0.48\linewidth}
\centering
\includegraphics[height=6cm,width=6cm]{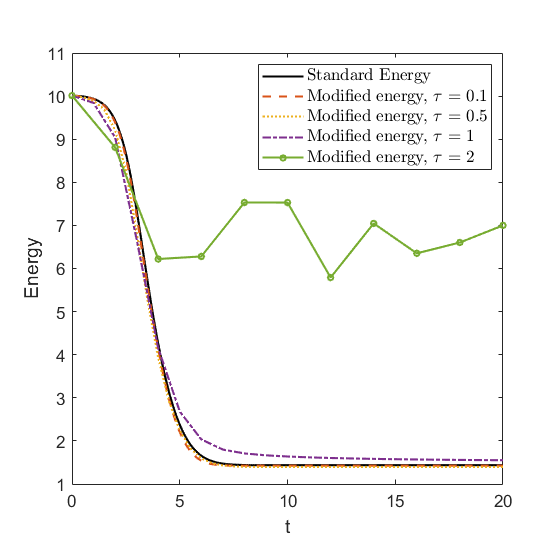}
\end{minipage}
}
\caption{Evolution of Energy and modified energy calculated by different methods; Left: FR-MLFEM. Righe: adaptive DLR-MLFEM. Parameters: $N_x = N_y = 129$, tolerance $\eta = 0.01\|\Sigma\|_2$. }\label{MEfig1}
\end{figure}

\begin{example}
In this example, we consider the classical AC equation with domain $\Omega = [0,2\pi] \times [0,2\pi]$ and the initial data
\begin{equation}
u_0 = 0.05sin(x)sin(y).
\end{equation}
\end{example}
We solve this problem using the FR-MLFEM and the adaptive DLR-MLFEM, with parameters set to $\epsilon = 0.01, k=1, N_x = N_y = 129$.
The rank adaptive tolerance is chosen as $0.01\|\Sigma\|_{2}$, where $\Sigma$ is the matrix defined in truncation step of \Cref{lrnonlinear2}.
In \Cref{MEfig1}, we plot the standard energy (computed by the FR-MLFEM with $\tau = 0.01$) and the modified energy of the numerical solutions.  
\Cref{MEfig1}(a) shows that as \(\tau\) increases from 0.1 to 2, the modified energy curve gradually deviates from the standard energy curve, eventually becoming oscillatory at \(\tau = 2\). However, when \(\tau\) is small (\(\tau = 0.1\)), there is no significant difference between the modified energy and standard energy curves.
The energy curves of the low-rank solutions are similar to those of the full-rank solution, consistent with the theoretical results in  \Cref{lem2.5} for the full-rank scheme and \Cref{energythm} for the low-rank scheme.

\subsection{Mass-conservation}
\begin{example}
Consider the classical AC equation \eqref{1.1} and the conservative AC equation \eqref{1.2}. The initial function is taken as the ``Kiss-Bubble" in $\Omega \in [-0.5, 0.5] \times [-0.5, 0.5]$,
\begin{equation}\label{4.2}
u_0(x,y)= - \sum\limits_{i=1}^{2}\tanh \left( \dfrac{\sqrt{(x-x_i)^2+(y-y_i)^2}-R}{\sqrt{2}\epsilon}\right)+1.
\end{equation}
The parameters are chosen as $\epsilon=0.01$, $R = 0.19, (x_1,y_1) = (0,-0.2)$, $(x_2,y_2) = (0,0.2)$. In addition, the parameters for both FR-MLFEM and DLR-MLFEM are chosen as $k=1, N_x = N_y = 256, \tau = 0.5$.
\end{example}

\begin{figure}
\centering
\subfigure{
\begin{minipage}[t]{0.16\linewidth}
\centering
\includegraphics[height=3.3cm,width=3.3cm]{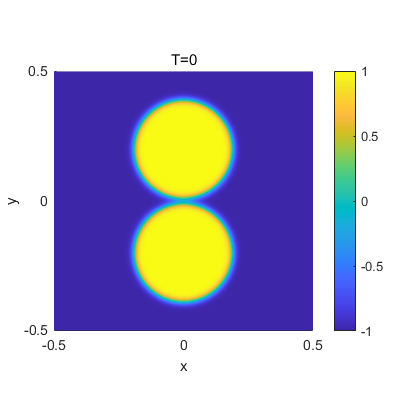}
\end{minipage}
}
\subfigure{
\begin{minipage}[t]{0.16\linewidth}
\centering
\includegraphics[height=3.3cm,width=3.3cm]{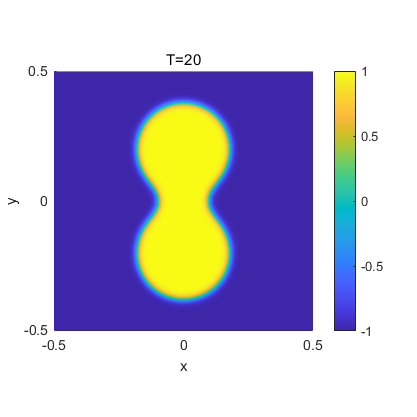}
\end{minipage}
}
\subfigure{
\begin{minipage}[t]{0.16\linewidth}
\centering
\includegraphics[height=3.3cm,width=3.3cm]{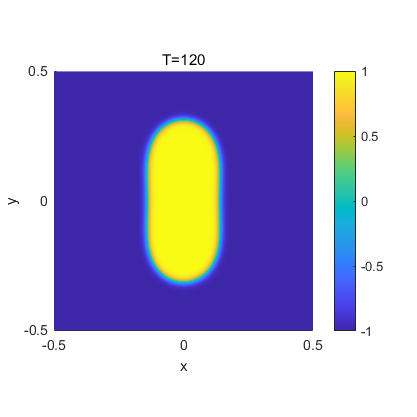}
\end{minipage}
}
\subfigure{
\begin{minipage}[t]{0.16\linewidth}
\centering
\includegraphics[height=3.3cm,width=3.3cm]{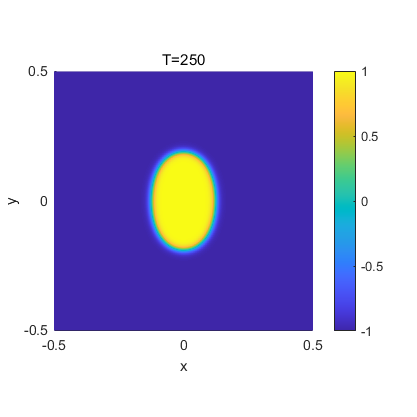}
\end{minipage}
}
\subfigure{
\begin{minipage}[t]{0.16\linewidth}
\centering
\includegraphics[height=3.3cm,width=3.3cm]{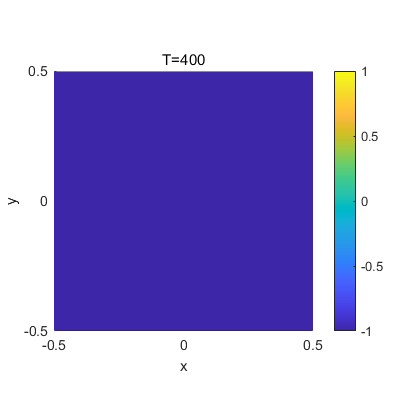}
\end{minipage}
}
\\
\vspace{-0.6cm}
\setcounter{subfigure}{0}
\subfigure[$t=0$]{
\begin{minipage}[t]{0.16\linewidth}
\centering
\includegraphics[height=3.3cm,width=3.3cm]{2D/exm3/KissBubble_T0.png}
\end{minipage}
}
\subfigure[$t=20$]{
\begin{minipage}[t]{0.16\linewidth}
\centering
\includegraphics[height=3.3cm,width=3.3cm]{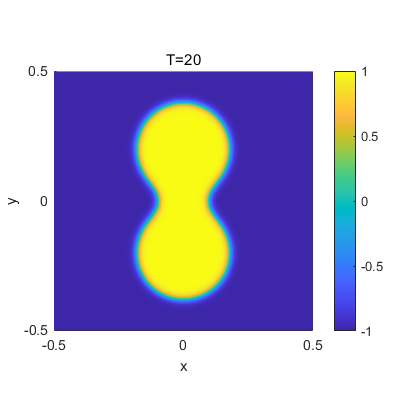}
\end{minipage}
}
\subfigure[$t=120$]{
\begin{minipage}[t]{0.16\linewidth}
\centering
\includegraphics[height=3.3cm,width=3.3cm]{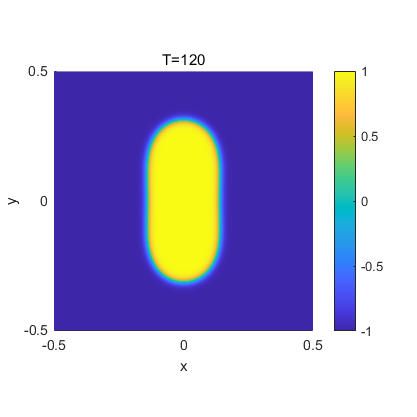}
\end{minipage}
}
\subfigure[$t=250$]{
\begin{minipage}[t]{0.16\linewidth}
\centering
\includegraphics[height=3.3cm,width=3.3cm]{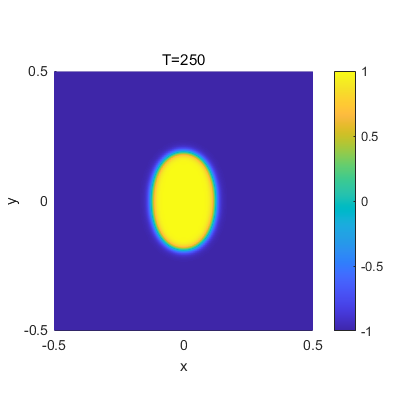}
\end{minipage}
}
\subfigure[$t=400$]{
\begin{minipage}[t]{0.16\linewidth}
\centering
\includegraphics[height=3.3cm,width=3.3cm]{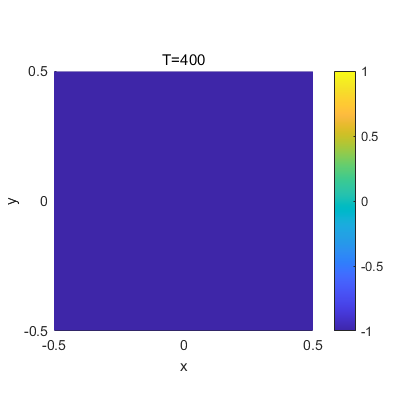}
\end{minipage}
}
\caption{Snapshots of solutions of AC equation computed using FR-MLFEM (first row) and adaptive DLR-MLFEM (second row). Parameters: $N_x = N_y = 256$, $\tau = 0.5$, $\eta = 0.01\|\Sigma\|_{2}$.}\label{2Dexm1fig1}
\end{figure}

\begin{figure}
\centering
\subfigure[Energy]{
\begin{minipage}[t]{0.31\linewidth}
\centering
\includegraphics[height=5cm,width=5cm]{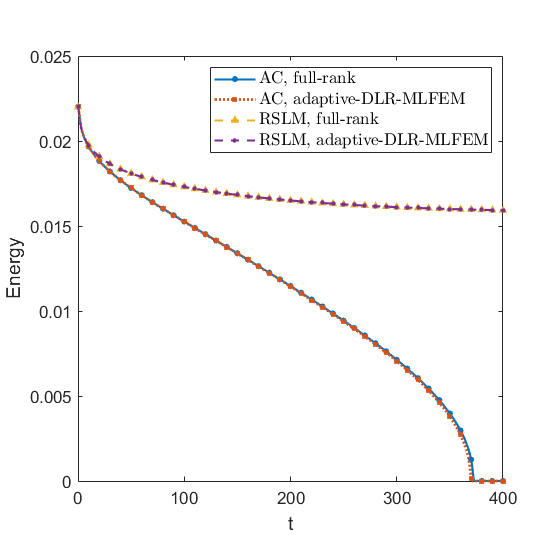}
\end{minipage}
}
\subfigure[mass error]{
\begin{minipage}[t]{0.31\linewidth}
\centering
\includegraphics[height=5cm,width=5cm]{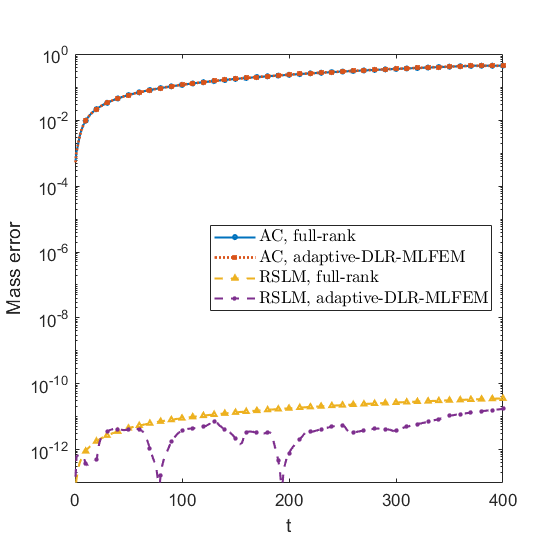}
\end{minipage}
}
\subfigure[Rank]{
\begin{minipage}[t]{0.31\linewidth}
\centering
\includegraphics[height=5cm,width=5cm]{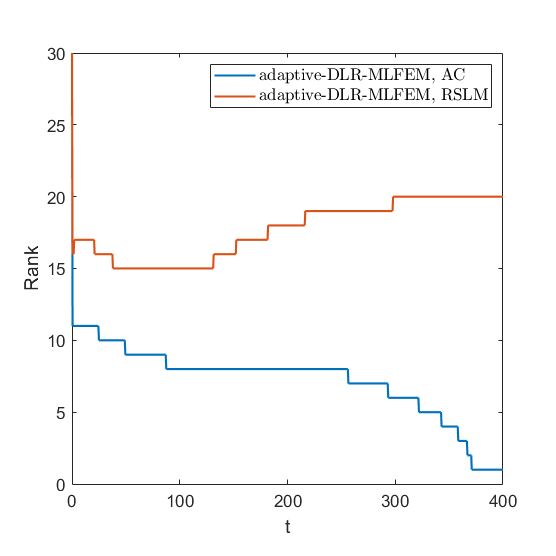}
\end{minipage}
}
\caption{Evolution of numerical solution properties. Left: energy; middle: mass error; right: Rank.}\label{2Dexm1fig3}
\end{figure}

We first solve the classical AC equation. The snapshots of the solutions computed by FR-MLFEM and adaptive DLR-MLFEM at different times are presented in \Cref{2Dexm1fig1}. As time evolves, the two bubbles begin to be absorbed by each other, coalescing into a single shrinking bubble, which eventually disappears. There is no visible difference between the full-rank and low-rank solutions. The evolutions of the mass changes, original energy, and rank are presented in \Cref{2Dexm1fig3}. From these results, we observe that: (a) the solution of the classical AC equation does not conserve mass, (b) both the full-rank and low-rank solutions satisfy the energy dissipation law, and (c) the rank in the low-rank algorithm decreases as the solution evolves toward the steady state.

\begin{figure}
\centering
\subfigure{
\begin{minipage}[t]{0.16\linewidth}
\centering
\includegraphics[height=3.3cm,width=3.3cm]{2D/exm3/KissBubble_T0.png}
\end{minipage}
}
\subfigure{
\begin{minipage}[t]{0.16\linewidth}
\centering
\includegraphics[height=3.3cm,width=3.3cm]{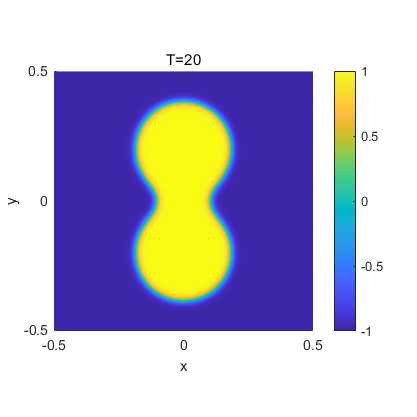}
\end{minipage}
}
\subfigure{
\begin{minipage}[t]{0.16\linewidth}
\centering
\includegraphics[height=3.3cm,width=3.3cm]{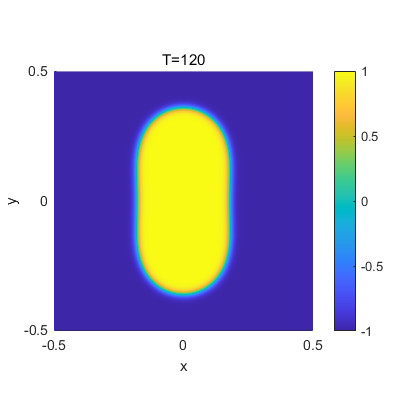}
\end{minipage}
}
\subfigure{
\begin{minipage}[t]{0.16\linewidth}
\centering
\includegraphics[height=3.3cm,width=3.3cm]{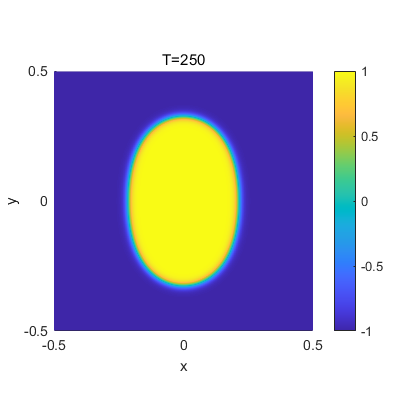}
\end{minipage}
}
\subfigure{
\begin{minipage}[t]{0.16\linewidth}
\centering
\includegraphics[height=3.3cm,width=3.3cm]{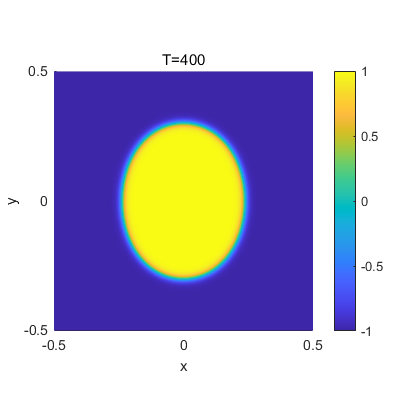}
\end{minipage}
}
\\
\vspace{-0.6cm}
\setcounter{subfigure}{0}
\subfigure[$t=0$]{
\begin{minipage}[t]{0.16\linewidth}
\centering
\includegraphics[height=3.3cm,width=3.3cm]{2D/exm3/KissBubble_T0.png}
\end{minipage}
}
\subfigure[$t=20$]{
\begin{minipage}[t]{0.16\linewidth}
\centering
\includegraphics[height=3.3cm,width=3.3cm]{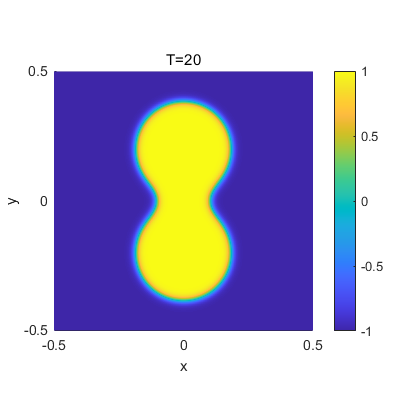}
\end{minipage}
}
\subfigure[$t=120$]{
\begin{minipage}[t]{0.16\linewidth}
\centering
\includegraphics[height=3.3cm,width=3.3cm]{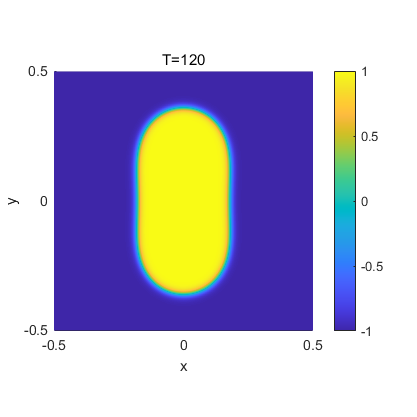}
\end{minipage}
}
\subfigure[$t=250$]{
\begin{minipage}[t]{0.16\linewidth}
\centering
\includegraphics[height=3.3cm,width=3.3cm]{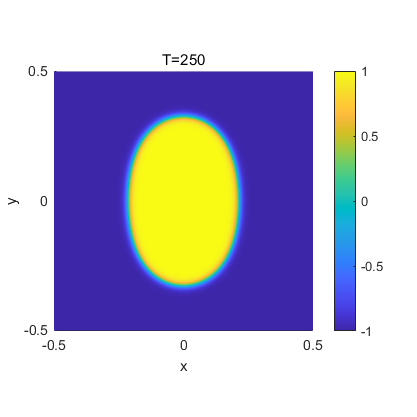}
\end{minipage}
}
\subfigure[$t=400$]{
\begin{minipage}[t]{0.16\linewidth}
\centering
\includegraphics[height=3.3cm,width=3.3cm]{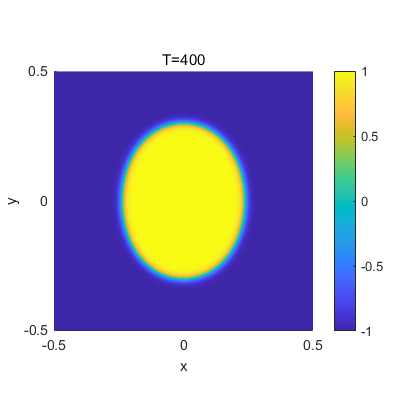}
\end{minipage}
}
\caption{Snapshots of solutions of AC equation with RSLM, computed using FR-MLFEM(first row) and adaptive DLR-MLFEM(second row). Parameters: $N_x = N_y = 256, \tau = 0.5, \eta = 0.001\|\Sigma\|_{2}.$}\label{2Dexm1fig2}
\end{figure}

Next, we solve the conservative AC equation with RSLM using the same initial condition. The snapshots of the solutions computed by FR-MLFEM and adaptive DLR-MLFEM are presented in \Cref{2Dexm1fig2}. In this case, the two bubbles coalesce into a larger, increasingly round bubble. Compared to the classical Allen-Cahn equation, the bubble does not disappear. From \Cref{2Dexm1fig3}, we can observe that (a) both the full-rank and low-rank solutions conserve mass very well, (b) both the full-rank and low-rank solutions satisfy the energy dissipation law, and (c) the rank in the low-rank algorithm is no more than $r=20$.

\begin{figure}
\centering
\subfigure[FR-MLFEM]{
\begin{minipage}[t]{0.22\linewidth}
\centering
\includegraphics[height=3.5cm,width=3.5cm]{2D/exm3/exm3_BBLM_Full1_T120.png}
\end{minipage}
}
\subfigure[adaptive DLR-MLFEM]{
\begin{minipage}[t]{0.22\linewidth}
\centering
\includegraphics[height=3.5cm,width=3.5cm]{2D/exm3/exm3_BBLM_Ada1_T120.png}
\end{minipage}
}
\subfigure[Mass error]{
\begin{minipage}[t]{0.22\linewidth}
\centering
\includegraphics[height=3.5cm,width=3.5cm]{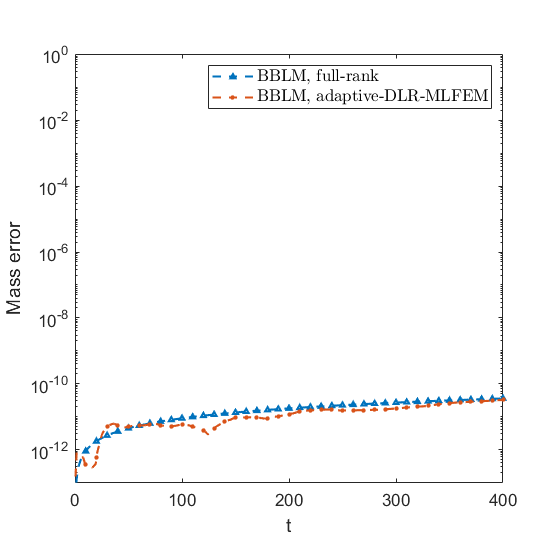}
\end{minipage}
}
\subfigure[Rank]{
\begin{minipage}[t]{0.22\linewidth}
\centering
\includegraphics[height=3.5cm,width=3.5cm]{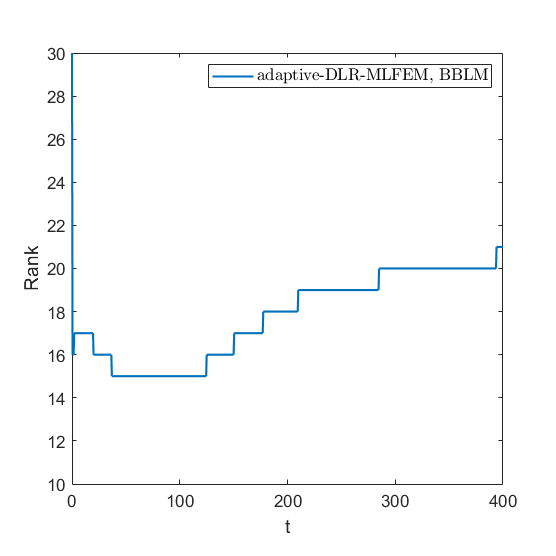}
\end{minipage}
}
\caption{The snapshots of solutions of AC equation with BBLM at $T=120$, and the evolution of mass error and rank calculated by different methods. Parameters: $N_x = N_y = 256, \tau = 0.5$, $\eta = 0.001\|\Sigma\|_2$. }\label{2Dexm2fig4}
\end{figure}

The energy of the conservative Allen-Cahn equation with BBLM is implicit. In \Cref{2Dexm2fig4}, we presented the snapshots of the solution at $T=20$, mass error and the rank history of the solution computed by FR-MLFEM and adaptive DLR-MLFEM. It can be observed that (a) There is no obvious difference between the solutions computed by FR-MLFEM and adaptive DLR-MLFEM, (b) both the full-rank and low-rank solutions conserve mass very well, and (c) the rank in the low-rank algorithm is no more than $r=21$.

In addition, for both cases, we observe that the profiles of the energy and mass errors are almost identical between the full-rank and low-rank solutions.

\subsection{Symmetry breaking}
Symmetry-breaking phenomena are commonly observed in phase field simulations. Theoretically, for a given initial odd function in a square domain centered at the origin, the exact solution should always maintain odd symmetry. In applications and numerical simulations, numerical schemes are expected to preserve this symmetry over time. However, in practice, this symmetry is often broken in long simulations \cite{li2021symmetry}.
Among existing algorithms, IMEX \cite{cheng2015fast}, operator splitting methods \cite{xu2019stability}, BDF2 \cite{xu2006stability}, and ETDRK \cite{wang2016efficient, du2021maximum} fail to preserve this symmetry, including the recently popular SAV approach \cite{shen2018scalar, shen2019new}. 

\begin{figure}
\centering
\subfigure[$u_1$]{
\begin{minipage}[t]{0.31\linewidth}
\centering
\includegraphics[height=4cm,width=4cm]{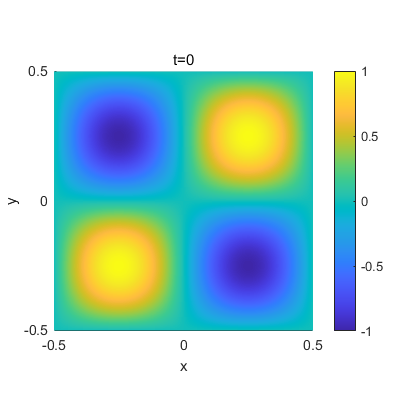}
\end{minipage}
}
\subfigure[$u_2$]{
\begin{minipage}[t]{0.31\linewidth}
\centering
\includegraphics[height=4cm,width=4cm]{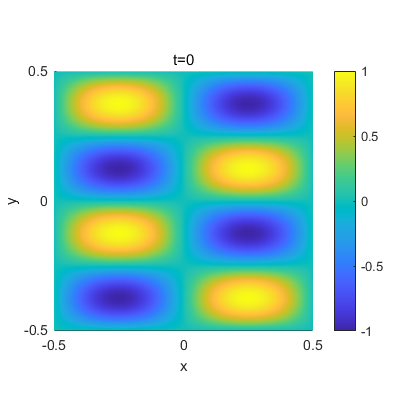}
\end{minipage}
}
\subfigure[$u_3$]{
\begin{minipage}[t]{0.31\linewidth}
\centering
\includegraphics[height=4cm,width=4cm]{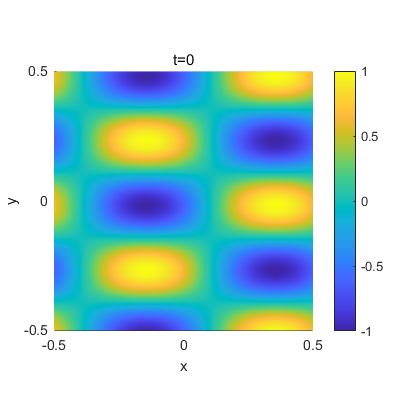}
\end{minipage}
}
\caption{Snapshots of different initial values.}\label{symfig0}
\end{figure}

\begin{figure}
\centering
\subfigure{
\begin{minipage}[t]{0.16\linewidth}
\centering
\includegraphics[height=3.3cm,width=3.3cm]{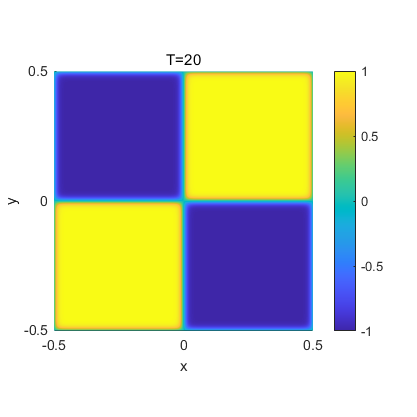}
\end{minipage}
}
\subfigure{
\begin{minipage}[t]{0.16\linewidth}
\centering
\includegraphics[height=3.3cm,width=3.3cm]{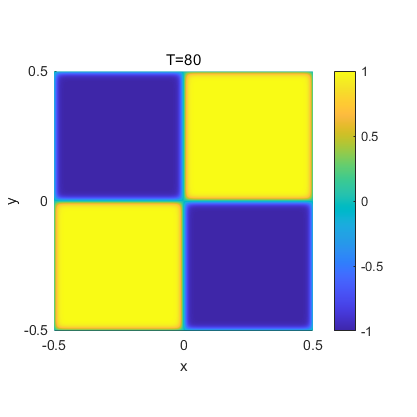}
\end{minipage}
}
\subfigure{
\begin{minipage}[t]{0.16\linewidth}
\centering
\includegraphics[height=3.3cm,width=3.3cm]{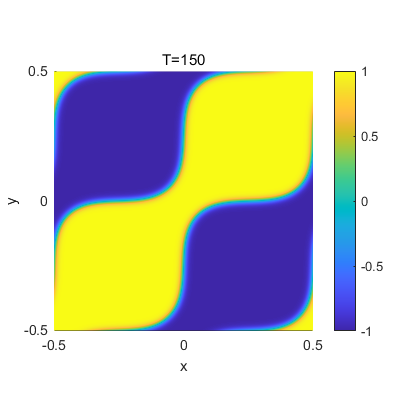}
\end{minipage}
}
\subfigure{
\begin{minipage}[t]{0.16\linewidth}
\centering
\includegraphics[height=3.3cm,width=3.3cm]{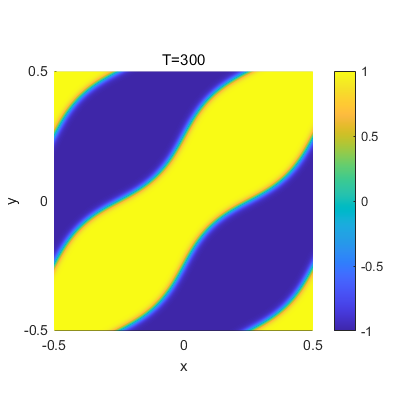}
\end{minipage}
}
\subfigure{
\begin{minipage}[t]{0.16\linewidth}
\centering
\includegraphics[height=3.3cm,width=3.3cm]{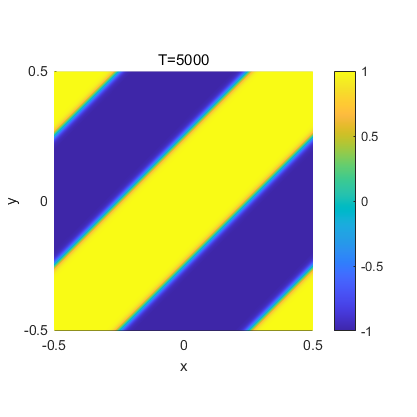}
\end{minipage}
}\\
\vspace{-0.6cm}
\setcounter{subfigure}{0}
\subfigure[$t=20$]{
\begin{minipage}[t]{0.16\linewidth}
\centering
\includegraphics[height=3.3cm,width=3.3cm]{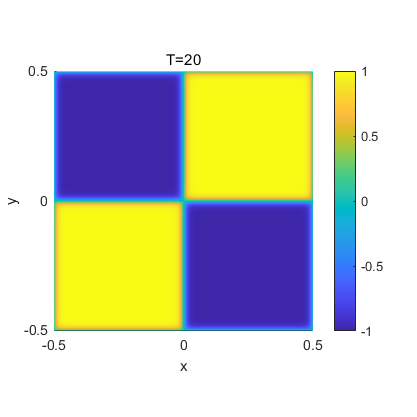}
\end{minipage}
}
\subfigure[$t=50$]{
\begin{minipage}[t]{0.16\linewidth}
\centering
\includegraphics[height=3.3cm,width=3.3cm]{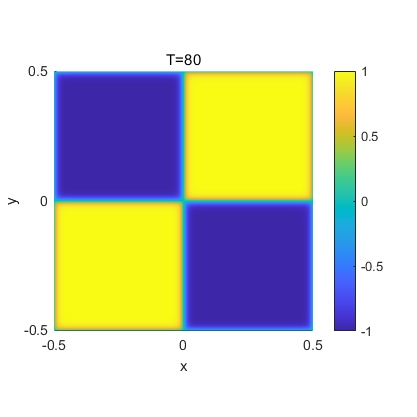}
\end{minipage}
}
\subfigure[$t=150$]{
\begin{minipage}[t]{0.16\linewidth}
\centering
\includegraphics[height=3.3cm,width=3.3cm]{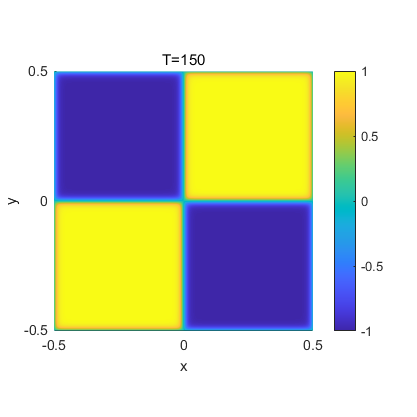}
\end{minipage}
}
\subfigure[$t=500$]{
\begin{minipage}[t]{0.16\linewidth}
\centering
\includegraphics[height=3.3cm,width=3.3cm]{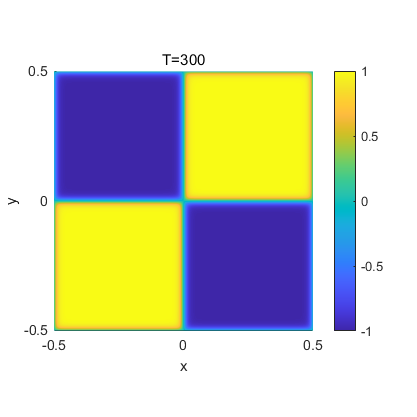}
\end{minipage}
}
\subfigure[$t=5000$]{
\begin{minipage}[t]{0.16\linewidth}
\centering
\includegraphics[height=3.3cm,width=3.3cm]{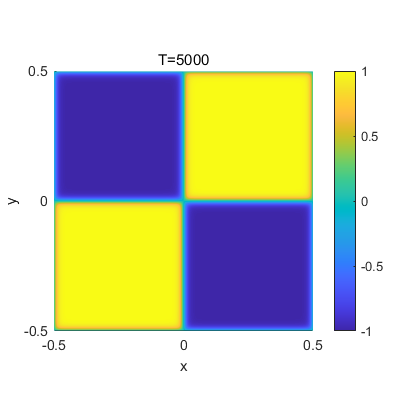}
\end{minipage}
}
\caption{Snapshots of solutions to the AC equation at different times with initial value $u_1$, computed using the FR-MLFEM (first row) and the Adaptive DLR-MLFEM (second row).}\label{symfig1}
\end{figure}

\begin{figure}
\centering
\subfigure{
\begin{minipage}[t]{0.16\linewidth}
\centering
\includegraphics[height=3.3cm,width=3.3cm]{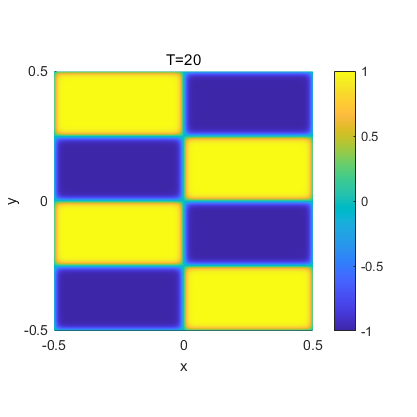}
\end{minipage}
}
\subfigure{
\begin{minipage}[t]{0.16\linewidth}
\centering
\includegraphics[height=3.3cm,width=3.3cm]{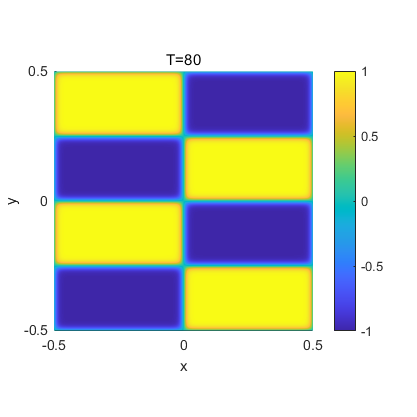}
\end{minipage}
}
\subfigure{
\begin{minipage}[t]{0.16\linewidth}
\centering
\includegraphics[height=3.3cm,width=3.3cm]{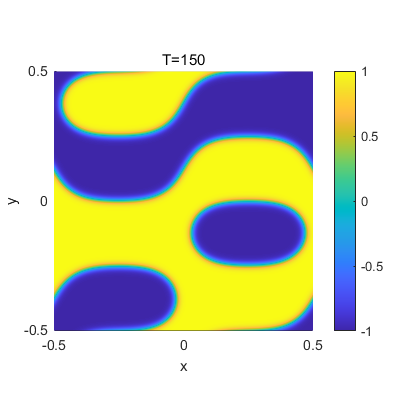}
\end{minipage}
}
\subfigure{
\begin{minipage}[t]{0.16\linewidth}
\centering
\includegraphics[height=3.3cm,width=3.3cm]{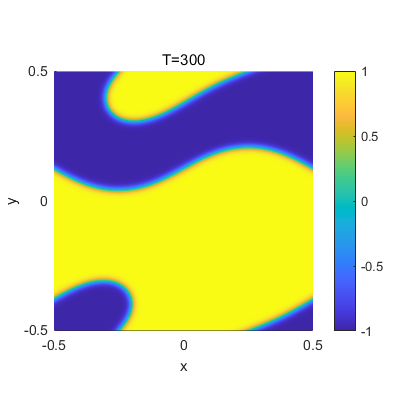}
\end{minipage}
}
\subfigure{
\begin{minipage}[t]{0.16\linewidth}
\centering
\includegraphics[height=3.3cm,width=3.3cm]{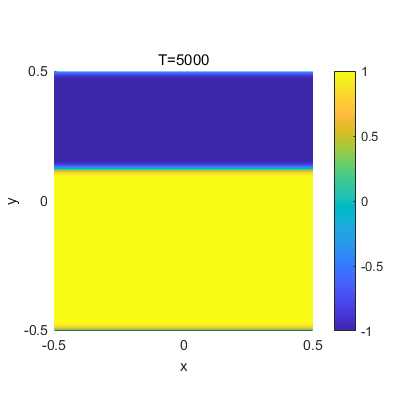}
\end{minipage}
}\\
\vspace{-0.6cm}
\setcounter{subfigure}{0}
\subfigure[$t=20$]{
\begin{minipage}[t]{0.16\linewidth}
\centering
\includegraphics[height=3.3cm,width=3.3cm]{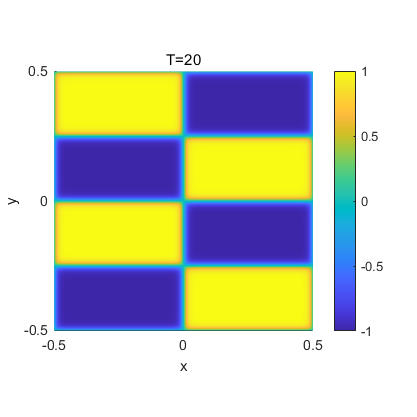}
\end{minipage}
}
\subfigure[$t=80$]{
\begin{minipage}[t]{0.16\linewidth}
\centering
\includegraphics[height=3.3cm,width=3.3cm]{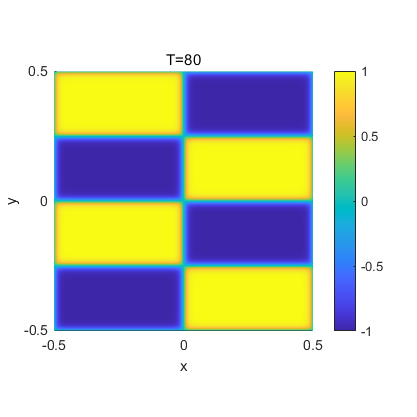}
\end{minipage}
}
\subfigure[$t=150$]{
\begin{minipage}[t]{0.16\linewidth}
\centering
\includegraphics[height=3.3cm,width=3.3cm]{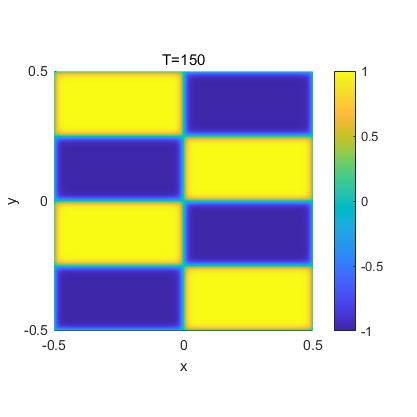}
\end{minipage}
}
\subfigure[$t=300$]{
\begin{minipage}[t]{0.16\linewidth}
\centering
\includegraphics[height=3.3cm,width=3.3cm]{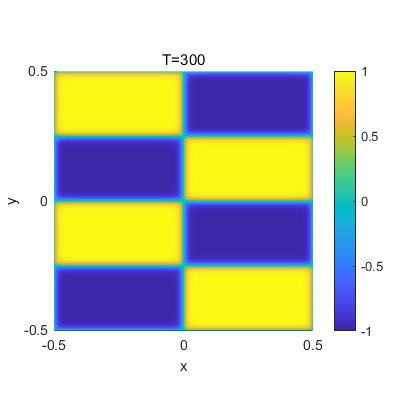}
\end{minipage}
}
\subfigure[$t=5000$]{
\begin{minipage}[t]{0.16\linewidth}
\centering
\includegraphics[height=3.3cm,width=3.3cm]{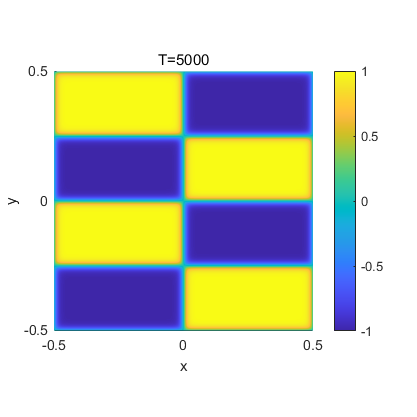}
\end{minipage}
}
\caption{Snapshots of solutions to the AC equation at different times with initial valu $u_2$, computed using the FR-MLFEM (first row) and the adaptive DLR-MLFEM (second row).}\label{symfig2}
\end{figure}

\begin{example}
We consider the classical Allen-Cahn equation in domain $\Omega \in [-0.5 , 0.5] \times [-0.5, 0.5]$ with different initial values
\begin{equation}\label{sym1}
\begin{aligned}
u_1&= sin(2 \pi x) sin(2 \pi y), \\
u_2&= sin(2 \pi x) sin(4 \pi y).
\end{aligned}
\end{equation}
We also consider a shifted initial value
\begin{equation}\label{sym2}
u_3= \sin\left(2\pi (x+\frac{\pi}{8})\right) \sin\left(4\pi (y+\frac{\pi}{8})\right).
\end{equation}
\end{example}
The snapshots of solutions for different initial values are shown in \Cref{symfig0}. Since the initial values in $\eqref{sym1}$ and $\eqref{sym2}$ are odd and smooth, the solutions should preserve odd symmetry over time, i.e., $u(-x,y,t) = -u(x,y,t)$ and $u(x,-y,t) = -u(x,y,t)$.
Figure \ref{symfig1} and Figure \ref{symfig2} compare the numerical solutions of FR-MLFEM and adaptive DLR-MLFEM. 
From the results, it is observed that the odd symmetry of the solution computed by FR-MLFEM is broken, whereas the adaptive DLR-MLFEM preserves symmetry throughout the simulation.

A similar numerical result for the case with the initial value $u_1$ was presented in \cite{li2021symmetry}, where symmetry-breaking, similar to that observed with the FR-MLFEM, was reported. To address this issue, a symmetry-preserving filter was introduced in \cite{li2021symmetry} to maintain the odd symmetry of the solution. Notably, the proposed adaptive DLR-MLFEM method achieves consistent results by inherently preserving the symmetry without requiring any such filter.

\begin{figure}
\centering
\subfigure{
\begin{minipage}[t]{0.16\linewidth}
\centering
\includegraphics[height=3.3cm,width=3.3cm]{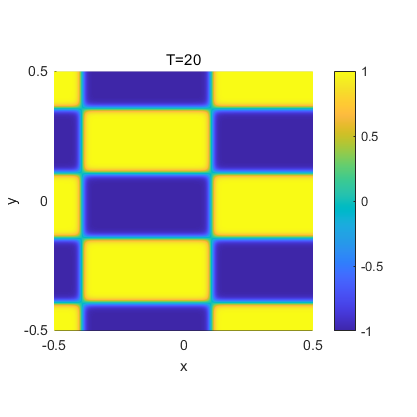}
\end{minipage}
}
\subfigure{
\begin{minipage}[t]{0.16\linewidth}
\centering
\includegraphics[height=3.3cm,width=3.3cm]{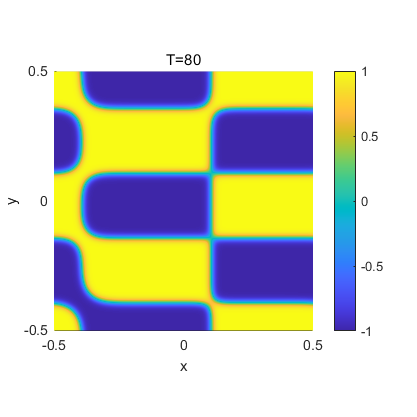}
\end{minipage}
}
\subfigure{
\begin{minipage}[t]{0.16\linewidth}
\centering
\includegraphics[height=3.3cm,width=3.3cm]{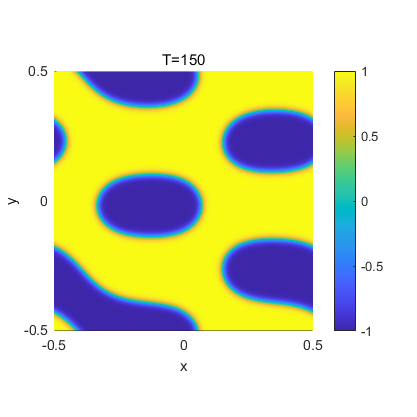}
\end{minipage}
}
\subfigure{
\begin{minipage}[t]{0.16\linewidth}
\centering
\includegraphics[height=3.3cm,width=3.3cm]{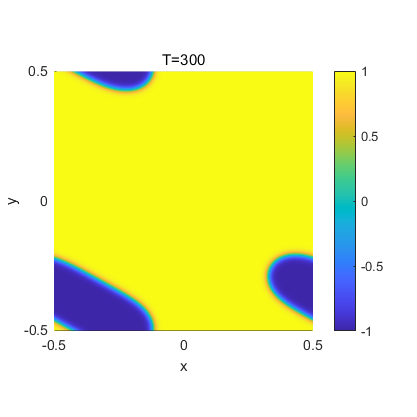}
\end{minipage}
}
\subfigure{
\begin{minipage}[t]{0.16\linewidth}
\centering
\includegraphics[height=3.3cm,width=3.3cm]{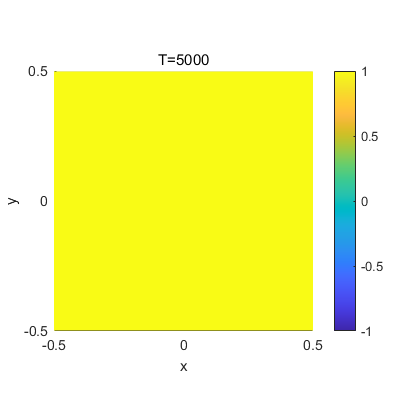}
\end{minipage}
}\\
\vspace{-0.6cm}
\setcounter{subfigure}{0}
\subfigure[$t=20$]{
\begin{minipage}[t]{0.16\linewidth}
\centering
\includegraphics[height=3.3cm,width=3.3cm]{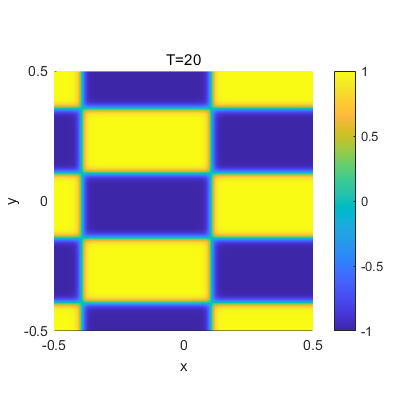}
\end{minipage}
}
\subfigure[$t=80$]{
\begin{minipage}[t]{0.16\linewidth}
\centering
\includegraphics[height=3.3cm,width=3.3cm]{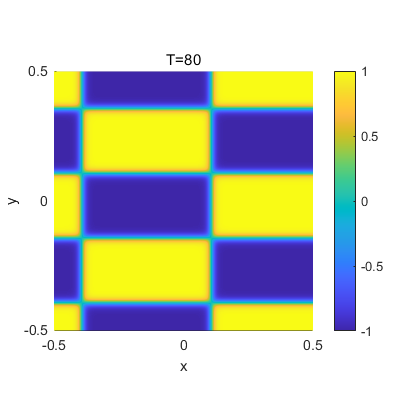}
\end{minipage}
}
\subfigure[$t=150$]{
\begin{minipage}[t]{0.16\linewidth}
\centering
\includegraphics[height=3.3cm,width=3.3cm]{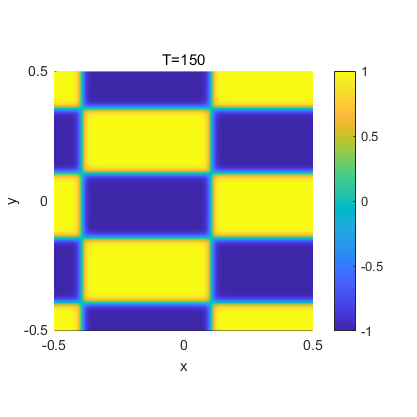}
\end{minipage}
}
\subfigure[$t=300$]{
\begin{minipage}[t]{0.16\linewidth}
\centering
\includegraphics[height=3.3cm,width=3.3cm]{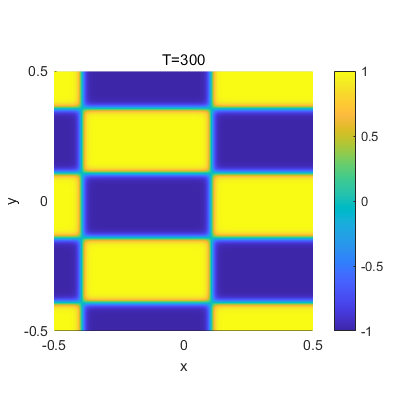}
\end{minipage}
}
\subfigure[$t=5000$]{
\begin{minipage}[t]{0.16\linewidth}
\centering
\includegraphics[height=3.3cm,width=3.3cm]{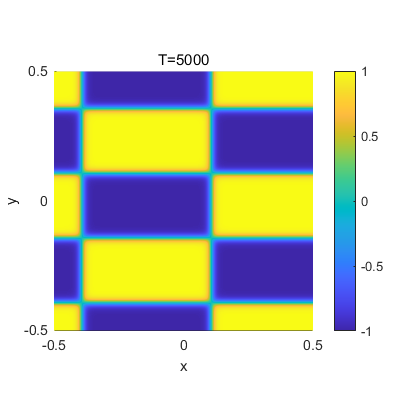}
\end{minipage}
}
\caption{Snapshots of solutions to the AC equation at different times with initial valu $u_3$, computed using the FR-MLFEM (first row) and the adaptive DLR-MLFEM (second row).}\label{symfig3}
\end{figure}

\begin{figure}
\centering
\subfigure[Energy]{
\begin{minipage}[t]{0.31\linewidth}
\centering
\includegraphics[height=5cm,width=5cm]{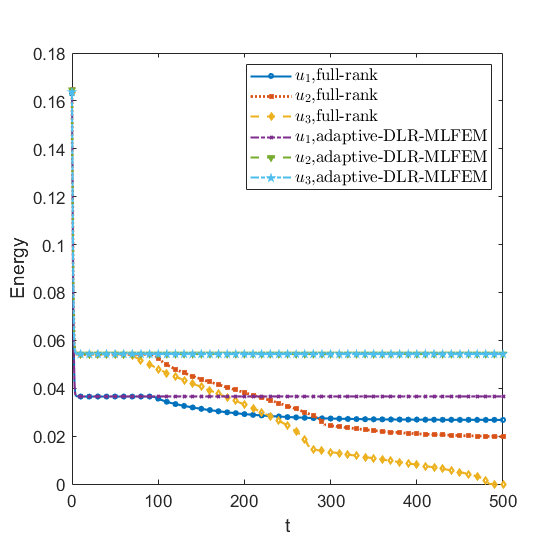}
\end{minipage}
}
\subfigure[Mass]{
\begin{minipage}[t]{0.31\linewidth}
\centering
\includegraphics[height=5cm,width=5cm]{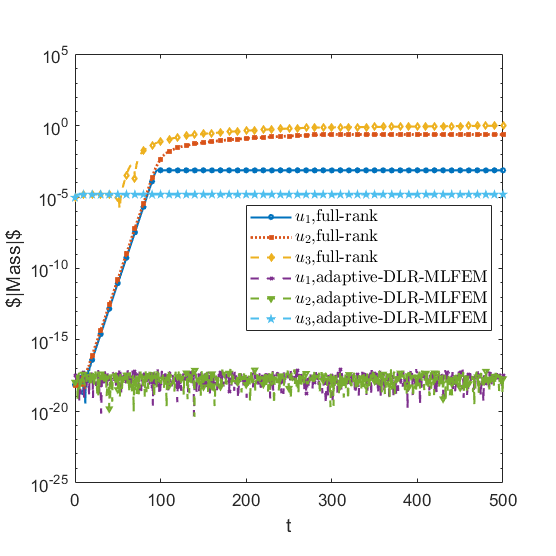}
\end{minipage}
}
\subfigure[Rank]{
\begin{minipage}[t]{0.31\linewidth}
\centering
\includegraphics[height=5cm,width=5cm]{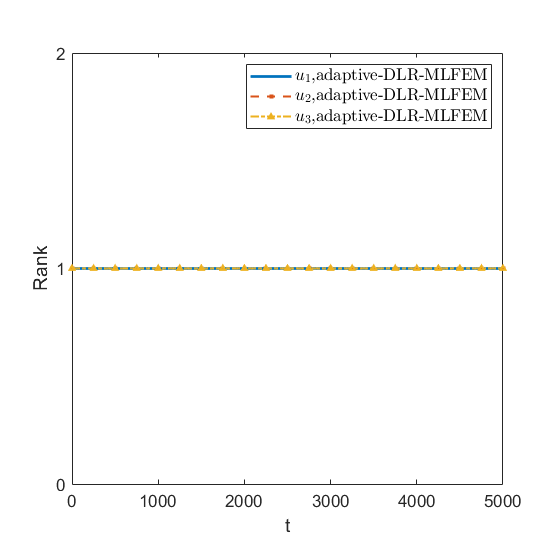}
\end{minipage}
}
\caption{Evolution of numerical solution properties. Left: Energy; middle: Mass; right: Rank.}\label{symfig4}
\end{figure}

For the case with the shifted initial value \(u_3\), snapshots of the full-rank and low-rank solutions are presented in \Cref{symfig3}, where it is evident that the FR-MLFEM fails to preserve the symmetry. The symmetry-preserving filter introduced in \cite{li2021symmetry} may not be applicable to improve it. In contrast, the adaptive DLR-MLFEM method successfully maintains the symmetry structure.
Figure \ref{symfig4} presents the evolution of energy, mass, and rank. It can be observed that the symmetry-breaking phenomenon of the FR-MLFEM method occurs at $t \approx 100$. 
In contrast, after a short period of evolution, all the properties of the solutions computed by the adaptive DLR-MLFEM method remain stable over a long period of simulation.

\begin{rem}
The reasons for this symmetry breaking may include: 1) contamination from machine rounding errors; 2) truncation errors from spatial discretization; 3) truncation errors from time discretization; or the accumulation and interaction of these factors, along with the influence of functional properties such as boundary layer effects.
Preserving the symmetry-breaking phenomena is important for numerical methods as it demonstrates their robustness in long-time simulations.
\end{rem}

\begin{rem}
The symmetry-preserving behavior of low-rank approximation integrators can be understood as follows: During the rank truncation process, the errors mentioned earlier are naturally eliminated. This effectively removes the small errors that accumulate over time and could eventually lead to symmetry breaking.
\end{rem}

\section{Concluding remarks}\label{Sec5}
A second-order DLR-MLFEM was proposed to solve both the classical and conservative AC equations. The matrix differential equation arising from the semi-discrete mass-lumped finite element method is decomposed into linear and nonlinear components. The linear component is solved analytically, while the nonlinear component is addressed using the second-order augmented BUG integrator. The computational complexity of the proposed algorithm is $\mathcal{O}((m+n)r^4)$. 
The mass is conserved up to the truncation tolerance in the conservative Allen-Cahn equation. The modified energy is dissipative up to a high-order error of $\mathcal{O}(\tau(\tau^2+h^{k+1}))$. Hence, the energy stability property remains valid.
Although some conservative BUG integrators have been developed in the literature (see, e.g., \cite{einkemmer2023conservation}), implementing them for the proposed method in the nonlinear problems is challenging, which will be addressed in future work. Numerical experiments demonstrate that the DLR-MLFEM solution preserves mass for the conservative AC equation and ensures energy dissipation. Additionally, the tests reveal that the DLR-MLFEM solution maintains symmetrical phenomena, highlighting its robustness for long-time simulations. Extending these results to three dimensions and exploring more general nonlinear terms, such as the logarithmic Flory-Huggins potential, would also be interesting directions for future research.

\section*{Acknowledgments}
Yi's research was partially supported by NSFC Project (12431014). Yin’s research was supported by the University of Texas at El Paso Startup Award.

\appendix

\section{Some proofs}
In this section, we present some necessary proofs.
\subsection{Proof of \Cref{lem2.3}}\label{lem2.2proof}
\begin{proof}
Let $q_m = [1,\cdots,1]^\top \in \mathbb{R}^{m}$ be a vector of ones.
It can be verified that  $(0, q_m)$ is an eigenpair of the matrix $ L_x = - M_x^{-1} A_x$, i.e., $L_x q_m = 0q_m$.
Then, it follows
\[
\begin{aligned}
e^{\tau \epsilon^2 L_x} q_m &= \sum\limits_{l=0}^{\infty} \frac{\tau^l}{l !} (\epsilon^2 L_x)^l q_m=\sum\limits_{l=0}^{\infty} \frac{(0\tau)^l}{l !} q_m  = I_m q_m  =q_m,
\end{aligned}
\]
which implies that $(0, q_m)$ is an eigenpair of the matrix $e^{\tau \epsilon^2 L_x}$.
Similarly, $(0, q_n)$ is an eigenpair of the matrix $e^{\tau \epsilon^2 L_y}$, i.e., $e^{\tau \epsilon^2 L_y} q_n = q_n$, where $q_n$ is also a vector of ones. 

Next, we state that the matrices $M_x e^{\tau \epsilon^2 L_x}$ and $M_y e^{\tau \epsilon^2 L_y}$ are symmetric. We will prove this for $M_x e^{\tau \epsilon^2 L_x}$; the case for $M_y e^{\tau \epsilon^2 L_y}$ can be shown similarly.
\[
\begin{aligned}
\left(M_x e^{\tau \epsilon^2 L_x} \right)^{\top} &= \left(e^{-\epsilon^2 \tau M_x^{-1} A_x} \right)^{\top} M_x^{\top} = e^{(-\epsilon^2 \tau M_x^{-1} A_x)^{\top}} M_x^{\top} =e^{-\epsilon^2 \tau A_x^{\top} (M_x^{-1})^{\top}} M_x^{\top} = e^{-\epsilon^2 \tau A_x M_x^{-1}} M_x \\
&= M_x M_x^{-1}e^{-\epsilon^2 \tau A_x M_x^{-1}} M_x = M_x e^{-\epsilon^2 \tau M_x^{-1} A_x M_x^{-1}M_x}=M_x e^{\tau \epsilon^2 L_x}.
\end{aligned}\] 
Then, using the cyclic property of the Frobenius inner product gives
\[
\begin{aligned}
(e^{\tau \epsilon^2 \mathcal{L}} Z, \mathbf{I})_{\rm M} &= (M_x e^{\tau \epsilon^2 L_x} Z e^{\tau \epsilon^2 L_y^{\top}} M_y,\mathbf{I})_{\rm F} =(Z , (M_x e^{\tau \epsilon^2 L_x})^{\top} \mathbf{I} (e^{\tau \epsilon^2 L_y^{\top}} M_y)^{\top})_{\rm F} = (Z , M_x e^{\tau \epsilon^2 L_x} q_m q_n^{\top} M_y e^{\tau \epsilon^2 L_y})_{\rm F} \\
&=( Z , M_x e^{\tau \epsilon^2 L_x} q_m ((M_y e^{\tau \epsilon^2  L_y })^{\top}q_n )^{\top})_{\rm F} = ( Z , M_x e^{\tau \epsilon^2 L_x} q_m (M_y e^{\tau \epsilon^2 L_y }q_n )^{\top})_{\rm F} \\
&=( Z  , M_x q_m ( M_y q_n)^{\top} )_{\rm F} =(M_x Z M_y  ,\mathbf{I})_{\rm F} = (Z, \mathbf{I})_{\rm M}.
\end{aligned}
\]
\end{proof}

\subsection{Proof of \Cref{lem2.5}}\label{lem2.5proof}
\begin{proof}
Recall that $W_{\mathfrak{n}+1,1} = e^{\frac{\tau}{2} \epsilon^2 \mathcal{L}} W_{\mathfrak{n}+1} = S_{\frac{\tau}{2}}^{\mathcal{L}} W_{\mathfrak{n}+1}$. Using $\eqref{2.25}$, it can be written as
\[W_{\mathfrak{n}+1,1} = S_{\frac{\tau}{2}}^{\mathcal{L}}S_{\frac{\tau}{2}}^{\mathcal{L}} S_{\tau}^{\mathcal{N}} W_{\mathfrak{n},1}=S_{\tau}^{\mathcal{L}} S_{\tau}^{\mathcal{N}} W_{\mathfrak{n},1},\]
which is equivalent to 
\[S_{-\tau}^{\mathcal{L}} W_{\mathfrak{n}+1,1} = S_{\tau}^{\mathcal{N}} W_{\mathfrak{n},1}.\]
Then it holds
\begin{equation}\label{2.30}
\frac{1}{\tau} \left( S_{-\tau}^{\mathcal{L}} W_{\mathfrak{n}+1,1} - W_{\mathfrak{n}+1,1} \right) + \frac{1}{\tau} (W_{\mathfrak{n}+1,1} - W_{\mathfrak{n},1}) = \frac{1}{\tau} (S_{\tau}^{\mathcal{N}} W_{\mathfrak{n},1} -W_{\mathfrak{n},1}).
\end{equation}
By $\eqref{2.25}$, it can show
\[S_\tau^{\mathcal{N}} W_{\mathfrak{n},1} - W_{\mathfrak{n},1} = \frac{\tau}{2} \left(\mathcal{N}(W_{\mathfrak{n},1}) + \mathcal{N} (W_{\mathfrak{n},2}) \right).\]
Thus, 
\begin{equation}\label{2.31}
\frac{1}{\tau} \left(S_{-\tau}^{\mathcal{L}}  W_{\mathfrak{n}+1,1} - W_{\mathfrak{n}+1,1} \right) + \frac{1}{\tau} (W_{\mathfrak{n}+1,1} - W_{\mathfrak{n},1}) = \frac{1}{2} \left(\mathcal{N}(W_{\mathfrak{n},1}) + \mathcal{N} (W_{\mathfrak{n},2}) \right).
\end{equation}
The equation $\eqref{2.31}$ can be reformulated in finite element form as
\[\frac{1}{\tau} \left( e^{-\tau \epsilon^2 \mathcal{L}_h} w_h^{\mathfrak{n}+1,1} - w_h^{\mathfrak{n}+1,1} \right) + \frac{1}{\tau} (w_h^{\mathfrak{n}+1,1} - w_h^{\mathfrak{n},1}) = \frac{1}{2} \left(\mathcal{N}(w_h^{\mathfrak{n},1}) + \mathcal{N} (w_h^{\mathfrak{n},2}) \right).\]
By taking the inner product with $w_h^{\mathfrak{n}+1,1} - w_h^{\mathfrak{n},1}$ and using the identity
\[(a-b,2a)_h = \|a\|^2_h - \|b\|_h^2 +\|a-b\|_h^2,\]
it follows
\begin{equation}\label{2.32}
\begin{aligned}
&\frac{1}{2\tau} \left[\left(  (e^{-\tau \epsilon^2 \mathcal{L}_h}-1) w_h^{\fn+1,1},w_h^{\fn+1,1} \right)_h - \left(  (e^{-\tau \epsilon^2\mathcal{L}_h}-1) w_h^{\fn,1} ,w_h^{\fn,1} \right)_h \right. \\
&\qquad + \left. \left(  e^{-\tau \epsilon^2\mathcal{L}_h} (w_h^{\mathfrak{n}+1,1} -w_h^{\mathfrak{n},1} ) ,w_h^{\mathfrak{n}+1,1} -w_h^{\mathfrak{n},1}\right)_h - \left\|w_h^{\mathfrak{n}+1,1} -w_h^{\mathfrak{n},1}\right\|^2_{h} \right] \\
&\qquad \qquad \qquad +\frac{1}{\tau} \left\|w_h^{\mathfrak{n}+1,1} -w_h^{\mathfrak{n},1}\right\|^2_{h} + \left( g(w_h^{\mathfrak{n},1}), w_h^{\mathfrak{n}+1,1} -w_h^{\mathfrak{n},1}\right)_h = 0.
\end{aligned}
\end{equation}
Expanding $G(w_h^{\mathfrak{n}+1,1})$ in Taylor series yields
\begin{equation}\label{2.33}
\begin{aligned}
G(w_h^{\mathfrak{n}+1,1}) - G(w_h^{\mathfrak{n},1}, 1)_h = (g(w_h^{\mathfrak{n},1}), w_h^{\mathfrak{n}+1,1} - w_h^{\mathfrak{n},1})_h+ \frac{1}{2}(g^{\prime}(\xi), (w_h^{\mathfrak{n}+1,1} - w_h^{\mathfrak{n},1})^2)_h.
\end{aligned}
\end{equation}
where $\xi$ is a function between $w_h^{\mathfrak{n}+1,1}$ and $w_h^{\mathfrak{n},1}$. Given $|g^{\prime}(\xi)| \leq C_{g}$ and $\tau \leq \frac{1}{C_g}$, it holds
\[
\begin{aligned}
\widetilde{E}(w_h^{\mathfrak{n}+1,1})- \tilde{E}(w_h^{\mathfrak{n},1})\leq -\left( \frac{1}{2\tau}  - \frac{1}{2} g^{\prime}(\xi)\right) \| w_h^{\mathfrak{n}+1,1} - w_h^{\mathfrak{n},1}\|^2_{h} \leq 0.
\end{aligned}
\]
\end{proof}

\bibliographystyle{abbrv}
\bibliography{ref.bib}

\end{document}